%% file: main.tex
\documentclass[a4paper]{article}
\usepackage[utf8]{inputenc}

\usepackage[pages=all, color=black, position={current page.south}, placement=bottom, scale=1, opacity=1, vshift=5mm]{background}
\SetBgContents{
}      %

\usepackage[margin=1in]{geometry} %
\usepackage{amsmath}
\usepackage{diagbox}
\usepackage{times}
\usepackage{bm}
\usepackage{algorithm}
\usepackage{algpseudocode}
\usepackage{tabularx}
\usepackage{array}
\usepackage{multirow}
\usepackage{booktabs}
\usepackage{comment}
\usepackage{float}
\usepackage{stmaryrd}
\usepackage{authblk}
\usepackage{amsthm}
\usepackage{pgf-umlcd}
\usepackage{amssymb}
\usepackage{graphicx}
\usepackage{pgfplots}
\usepackage{todonotes}
\usepackage{subcaption}
\usepackage[section]{placeins}
\usepackage{forest}
\usepackage{graphics}
\usepackage{graphicx}

\usepackage[normalem]{ulem}

\usepackage{xcolor}   

\tikzset{bplus/.style={rectangle split, rectangle split horizontal,text width=1em, text centered, inner xsep=2pt,draw}}
\forestset{multi edge/.style={edge path={\noexpand\path[\forestoption{edge}] (!u.#1 south)--(.child anchor)\forestoption{edge label};}}}
\ExplSyntaxOn
\seq_new:N \l_tree_node_seq
\seq_set_from_clist:Nn \c_node_names_seq {one,two,three,four}
\NewDocumentCommand{\mpn}{m}{
    \seq_set_from_clist:Nn \l_tree_node_seq {#1}
    \seq_map_indexed_inline:Nn \l_tree_node_seq { 
        \nodepart{\seq_item:Nn \c_node_names_seq {##1}} {##2}
    }
}
\ExplSyntaxOff  

\newcommand{\norm}[1]{\left\lVert#1\right\rVert}
\newcommand{\ddd}{\text{\rm D}}

\newcommand{\uu}[1]{\mathbf{#1}}
\newcommand{\mean}[1]{\{\!\!\{#1\}\!\!\}}                %
\newcommand{\jump}[1]{[\![#1]\!]}                        %
\newcommand{\ud}{\,\mathrm{d}}
\newcommand{\mbf}[1]{\mbox{\boldmath$\rm{#1}$}}

\newcommand{\diam}{\operatorname{diam}}
\newcommand{\dint}{\text{\rm int}}

\usepackage{pifont}%
\usepackage{siunitx}
\usepackage[export]{adjustbox}
\usepackage{pgfplotstable}

\theoremstyle{definition}
\newtheorem{definition}{Definition}[section]

\usepgfplotslibrary{colorbrewer}
\pgfplotsset{
    cycle list/Dark2-7,
    cycle multiindex* list={
        mark list*\nextlist
        Dark2-7\nextlist
        },
}
\pgfplotsset{every axis plot/.append style={thick}}
\usepgfplotslibrary{groupplots,external}
\usetikzlibrary{er,positioning,external,shapes}
\tikzset{external/only named=true}
\usepackage{multicol}
\usepackage{mathtools}
\usepackage{caption} %

\renewcommand*\Call[2]{\textproc{#1}(#2)}
\usepackage{hyperref}
\hypersetup{
    colorlinks = true,
}

\usepackage[sort&compress,numbers,square]{natbib}
\usepackage{graphicx}
\graphicspath{{}}
\usepackage{mathrsfs} %

\pgfplotsset{compat=1.9}

\begin{document}
\captionsetup[table]{skip=6pt}

\title{R3MG: R-tree based agglomeration of polytopal grids with applications to multilevel methods}

\author[1]{Marco Feder}
\author[1]{Andrea Cangiani}
\author[2]{Luca Heltai}
\affil[1]{MathLab, Mathematics Area, International School for Advanced Studies (SISSA), Via Bonomea 265, Trieste, 34136, Italy}
\affil[2]{Department of Mathematics, University of Pisa, Largo Bruno Pontecorvo, 5, Pisa, 56127, Italy}

\maketitle

\begin{abstract}
    We present a novel approach to perform agglomeration of polygonal and polyhedral grids based on spatial indices. Agglomeration strategies are a key ingredient in polytopal methods for PDEs as they are used to generate (hierarchies of) computational grids from an initial grid. Spatial indices are specialized data structures that significantly accelerate queries involving spatial relationships in arbitrary space dimensions. We show how the construction of the R-tree spatial database of an arbitrary fine grid offers a natural and efficient agglomeration strategy with the following characteristics: i) the process is fully automated, robust, and dimension-independent, ii)  it automatically produces a balanced and nested hierarchy of agglomerates, and iii) the shape of the agglomerates is tightly close to the respective axis aligned bounding boxes. 
    Moreover, the R-tree approach provides a full hierarchy of nested agglomerates which permits fast query and allows for efficient geometric multigrid methods to be applied also to those cases where a hierarchy of grids is not present at construction time. 
    We present several examples based on polygonal discontinuous Galerkin methods, confirming the effectiveness of our approach in the context of challenging three-dimensional geometries and the design of geometric multigrid preconditioners.
        \newline
        \newline
\end{abstract}

\section{Introduction}\label{sec:intro}

Traditional Finite Element Methods (FEMs) often struggle to accurately represent complex geometries coming from real-world applications, where grid generation -- a crucial step in FEMs -- is often a bottleneck due to the limitations of using only tetrahedral, hexahedral, or prismatic elements.

To overcome these limitations, a natural possibility is to use methods that leverage general polygons and polyhedra as mesh elements. These methods, known as polytopal methods, have witnessed tremendous development in the last two decades and have shown to be particularly effective in representing complex geometries. We refer, for instance,  to polygonal FEM \cite{Tabarraei:Sukumar:2004,Fries:Belytschko:2009},
Mimetic Finite Differences \cite{MimeticBook2014,HYMAN1997130,Mimetic2005}, Virtual Element Method (VEM) \cite{VEM6,VEMbook}, Polygonal Discontinuous Galerkin  \cite{dgp,BASSI201245,AntoniettiGeo}, Hybridized Discontinuous Galerkin (HDG) \cite{HDG2009,HDG_proj,HDG2008}, and Hybrid High-Order
method (HHO) \cite{DIPIETRO20151,DIPIETRO201531}.

Regardless of the underlying method, generating hierarchies of computational grids is a key step in the efficient numerical solution of partial differential equations. For simple geometries and traditional FEMs, the generation of these hierarchies is often straightforward, and it is performed in a bottom-up approach, by refinement of an initially coarse grid. For complex geometries, however, the construction of a hierarchy of grids may be challenging or impractical. One may only have access to non-nested sequences, which requires
the construction of challenging transfer operators with a possibly high computational cost due to the non-matching nature of consecutive levels. For instance, this may occur when CAD (computer-aided design) models are meshed with external software and one is given a very fine geometry for which no hierarchical information is available.

The use of polytopal methods is attractive in this respect since coarse grids can be simply generated by merging polygonal and polyhedral elements \cite{Chan1998AnAM,antonietti2020agglomeration,BASSI201245,PAN2022110775}. However, providing automated and good-quality agglomeration strategies for polygonal and polyhedral elements remains an open and challenging task. 
It is indeed crucial to preserve the original mesh quality, as any deterioration could
potentially affect the overall performance of the method in terms of stability and accuracy. 

During the last few years, there has been a significant advancement in Machine Learning (ML) algorithms, providing a framework to automatically extract information from data to enhance and accelerate numerical methods for scientific computing~\cite{RAISSI2019686,RAY2018166}, including for effective polygonal mesh agglomeration~\cite{ANTONIETTI2022111531,antonietti2024polytopalmeshagglomerationgeometrical}, adaptivity~\cite{ANTONIETTI2022110900}, and enhancement of linear solvers~\cite{caldana2023deep}.

In this work, we propose an alternative approach to produce efficient and high-quality agglomerates based on R-trees~\cite{Guttman,RtreeBook}, a spatial indexing data structure.
R-trees excel at organizing spatial data using bounding boxes, particularly in contexts where performing spatial queries for large sets of geometric objects in a fast way is required. Such geometric predicates are ubiquitous within the efficient implementations of
non-matching finite element techniques, where two different geometries, usually representing different physical components, overlap arbitrarily such as in immersed methods \cite{KrauseZulian,Boffi23} or in particle methods \cite{JOACHIM2023112189}.
By leveraging the R-trees structure, we develop an efficient and automatic mesh agglomeration algorithm specifically tailored for polygonal and polyhedral grids. This approach is fully automated, robust, and dimension-independent. It is designed to
be independent of the specific shape of the underlying fine elements, which can indeed be polytopic; hexahedral and simplicial elements are particular cases of this. Moreover, we show that it preserves mesh quality while significantly reducing the computational cost associated with the agglomeration process. The effectiveness of the proposed approach as a polytopic finite element mesh generator and agglomerator is validated within a polytopic Discontinuous Galerkin framework applied to the Poisson problem. However, we stress that
the procedure depends only on the initial grid and not on the PDE at hand.

To validate the methodology, we assess the quality metrics of the polygonal elements produced from a set of representative fine grids. We compare our findings with agglomerations obtained with \textsc{METIS} \cite{METIS}, a state-of-the-art graph partitioning tool, which is designed to process only the graph information regarding the mesh, and is often employed as an agglomeration algorithm. Our results show that the R-trees approach produces meshes which are either similar or superior in quality to those produced by \textsc{METIS} using a fraction of the computational time.
In particular, R-tree based agglomeration preserves structured meshes, a property not shared by \textsc{METIS}. Thus, for instance, when starting from a square grid, the repeated application of R-tree agglomeration produces a sequence of nested square grids. Moreover, even in the case of unstructured meshes, subsequent agglomerates tend to align to the corresponding bounding boxes, thus the resulting grid appears logically rectangular.

A crucial feature of our algorithm is that, by exploiting the tree structure and appropriate parent-child relationships, it becomes possible to automatically extract \emph{nested} sequences of agglomerated meshes, which can then be directly used within multigrid solvers. %
We propose R-tree based MultiGrid (\textsc{R3MG}) preconditioning with discontinuous Galerkin methods. 
We test \textsc{R3MG} on a second-order elliptic model problem discretized using the polygonal interior penalty discontinuous Galerkin method analysed in~\cite{Cangiani_PolyDG,dgp,dgease}.
In the context of polytopic methods, for instance with VEM discretizations, one has to fallback to non-nested approaches: in consequence, ad-hoc definitions of the
transfer operators have to be carried out, see e.g.~\cite{VEMMG}. However, for DG methodologies both nested (~\cite{hpMGAntonietti}) and non-nested approaches, cf.~\cite{AntoniettiPennesi} are possible.
Exploiting the nested sequences of polytopal meshes produced by R-tree agglomeration, we construct multigrid preconditioners and report on the number of iterations required by the preconditioned conjugate gradient iterative solver. Our results show good convergence properties of both the two- and three-level preconditioned iterative solver in two- and three-dimensions.
The experiments are carried out using a newly developed C++ library based on the well-established \textsc{deal.II} Finite Element Library~\cite{dealIIdesign}. Our implementation allows for the use of polygonal discontinuous Galerkin
methods on agglomerated grids both in 2D and 3D, in parallel, and with different agglomeration strategies.

The paper is organized as follows. In Section~\ref{sec:rtree} we recall the R-tree data structure and the algorithmic realization of our approach, which is then used to build hierarchies of agglomerated grids in Section~\ref{sec:rtree-agglomeration}. We validate the R-tree based agglomeration strategy in Section~\ref{sec:validation}, while in Section~\ref{sec:polydg} we
introduce the polygonal discontinuous Galerkin method for the discretization of second-order elliptic problems on general meshes and showcase the application of R-tree based agglomeration in the solution of benchmark PDEs as well as
its application to multigrid methodology. Finally, Section~\ref{sec:conc} summarizes our conclusions and points to further research directions.

\input{rtree_data_structure}

\input{agglomeration_rtree}

\input{validation}
\input{poly_dg}
\input{tests}

\input{applications}

\section{Conclusions and outlook}\label{sec:conc}

We presented a novel method to generate  agglomerated grids based on a spatial structure that can be easily built out of an underlying fine mesh. The method  relies on the so-called R-tree data structure and generates sequences of nested agglomerated grids, independently on the shape of the underlying fine mesh.  We showed its robustness and effectiveness on a sequence of two- and three-dimensional benchmarks and
non-trivial geometries. 
The new method is also considerably faster than graph-partitioning-based approaches such as \textsc{METIS}. Hinging upon the reported computing times, we claim
that the computational cost relative to the agglomeration procedure is negligible with respect to all the components involved in a full pipeline of a finite element simulation. 

Further, we proposed the use of the sequence of nested grids produced via R-tree agglomeration to construct multigrid preconditioners, coining the new approach as the R-tree based  Geometric Multigrid (R3MG) algorithm.
In practice R3MG  shows the good convergence properties typical of multilevel solvers, also in cases where a hierarchy of grids was not available to start with. Several research directions can be explored.
The theoretical analysis of the convergence properties of the multigrid method covering an arbitrary number of agglomeration levels is under investigation. %
From the implementation standpoint, a careful study of the performance of the method
to assess its scalability on larger problems will be the subject of a dedicated work. In this regard, it is well known that the time to assemble the linear system is the most computationally intensive operation. However, this cost can be dramatically reduced, for instance,  by employing (if possible) the so-called quadrature free algorithm as shown in \cite{quadraturefree} and more recently in \cite{qudaraturefreehouston}, which allows
the computation of volume integrals without resorting to the sub-tessellation of the background mesh $\Omega$.

\section*{Acknowledgments}
LH acknowledges the MIUR Excellence Department Project awarded to the Department of Mathematics, University of Pisa, CUP I57G22000700001, and partial support from grant MUR PRIN 2022 No. 2022WKWZA8 “Immersed methods for multiscale and multiphysics problems (IMMEDIATE)”.
AC acknowledges partial support from PNRR NGE iNEST “Interconnected Nord-Est Innovation Ecosystem” project.
LH and AC acknowledge the support of the European Research Council (ERC) under the European Union's Horizon 2020 research and innovation programme (call HORIZON-EUROHPC-JU-2023-COE-03, grant agreement No. 101172493 ``dealii-X''). 
The authors are members of Gruppo Nazionale per il Calcolo Scientifico (GNCS) of Istituto Nazionale di Alta Matematica (INdAM).
\bibliography{refs}
\end{document}

%% file: rtree_data_structure.tex
\section{R-trees}
\label{sec:rtree}
We briefly list the basic properties of the R-tree data structure proposed by Guttman in the seminal paper~\cite{Guttman} and discuss its variations following~\cite{RtreeBook}.

R-trees are hierarchical data structures used for the dynamic organization of collections of $d$-dimensional geometric objects, representing them by their \emph{minimum bounding -- axis aligned -- d-dimensional rectangle}, also denoted by MBR. In this
context, dynamic means that no global reorganization is required upon insertion or deletion of new elements of the tree. 
We will use the terms \emph{MBR} or \emph{bounding box} interchangeably throughout the work. An internal node of an R-tree consists of
the MBR that bounds all its children. In particular, each internal node stores two pieces of data: a way of identifying a child node
and the MBR of all entries within this child node. The actual data is stored in the \emph{leaves} of the tree, i.e., the terminal nodes of the tree data strcuture. We summarize these aspects in the following definition of R-tree. 

\begin{figure}[htb]

  \begin{minipage}[b]{0.5\textwidth}
  \centering
  
  \includegraphics{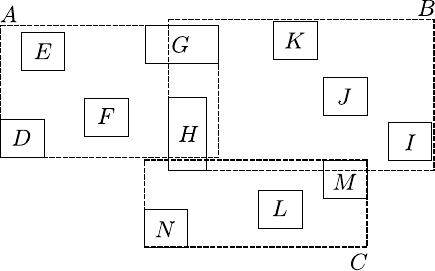}
  
  \caption{Example of MBRs holding geometric data and their MBRs.}\label{fig:rtree_example}
  
  \end{minipage}\hfill %
  \begin{minipage}[b]{0.5\textwidth}
  \centering
  \begin{forest}for tree={bplus, edge={->},l sep=1cm}
    [\mpn{A,B,C,\phantom{M}}
        [\mpn{D,E,F,G},multi edge=one]
        [\mpn{H,I,J,K},multi edge=two]
        [\mpn{L,M,N,\phantom{M}},multi edge=three]
    ]
  \end{forest}
  
  \caption{Corresponding R-tree data structure.}\label{fig:rtree_diagram}
  \end{minipage}
  
\end{figure}

\begin{definition}{R-tree of order $(m,M)$.}
  \label{def:definition_rtree}
  An R-tree of order $(m,M)$ has the following characteristics:
  \begin{itemize}
    \item Each \emph{leaf} node %
    is a container that can host between a minimum of $m\leq\frac{M}{2}$ and a maximum of $M$ entries (except for the root node, which is allowed to contain fewer elements if the number of objects to classify is less than $m$). Each entry is of type $(\text{MBR},\text{id})$ where $\text{id}$ is the object's identifier and $\text{MBR}$ is the minimum bounding rectangle that covers the object.
    \item Each \emph{internal} node is a container that can store again between a minimum of $m\leq\frac{M}{2}$ and a maximum of $M$ entries. Each entry is of the form $(\text{MBR},\text{p})$  where $\text{p}$ is a pointer to a child node (that can be either an \emph{internal} or a \emph{leaf} node) and $\text{MBR}$ is the minimum bounding rectangle that covers all the MBRs contained in this child.
    \item All leaves of the R-tree are at the same level.
    \item The minimum allowed number of entries in the root node is $2$ unless it is a
    leaf node, in which case it may contain zero or one entry (and this is the only exception to the rule that all leaves must contain a minimum of $m$ entries).
  \end{itemize}
\end{definition}
An example of the minimum bounding rectangles of some geometric objects (not shown) is given in Figure \ref{fig:rtree_example}, while the associated R-tree of
order $(2,4)$ is shown in Figure \ref{fig:rtree_diagram}. In this example, the leaf level stores the minimum bounding rectangles $D,E,F,G,H,I,J,K,L,M,N$. Conversely, the internal node comprises the three MBRs $A,B,C$. Notice the empty box in the internal node and the rightmost leaf node, meaning another entry could be stored there since the order is $(2,4)$.

The original R-tree is based solely on the minimization of the measure of each MBR. Several variants have been proposed, aimed at either improving performance or flexibility, generally depending on the domain of application.  

In view of optimizing quality and handling of grid generation, it is desirable to minimize the \emph{overlap} between MBRs. Indeed, the larger the overlap, the larger the number of paths to be processed during queries. Moreover, the smaller the overlap the closer the agglomerated element will conform to the corresponding bounding boxes, thus making the resulting agglomerated grid qualitatively close to rectangular.  Among the several available variants, we adopt the R$^*$-tree data structure designed in \cite{Beckmann1990TheRA}.
The criteria that R$^*$-trees aim to achieve are the following: 
\begin{itemize}
  \item \emph{Minimization of the area covered by each MBR.} This criterion aims at minimizing the dead space (area covered by MBRs but not by the enclosed elements) to reduce the number of paths pursued during query processing.
  \item \emph{Minimization of overlap between MBRs.} The larger the overlap, the larger the expected number of paths followed for a query. As such, this criterion has the same objective as the previous one.
  \item \emph{Minimization of MBRs perimeters.} Shaping more square bounding boxes results in reduced query time as this suffers in the presence of large overlaps and/or heterogeneous shapes. Moreover, since square objects are packed more easily, the corresponding MBRs at upper levels are expected to be smaller.
  \item \emph{Maximization of storage utilization.} Increasing the storage utilization per node will generally reduce the cost of queries since the height of the tree will be low. This holds especially true for larger queries, where a significant portion of the entries satisfies the query. Conversely, when storage utilization is low, more nodes tend to be invoked during query processing.
\end{itemize}
 It should be noted that such requirements can easily be orthogonal. As an example, keeping the area and the overlap between MBRs low could imply a lower number of entries packed within each node of the tree, leading to higher storage usage. Therefore, the R$^*$-tree does some heuristics to find the best possible compromise of these criteria. We refer to the original paper \cite{Beckmann1990TheRA} for the relevant algorithmic details.

In our implementation, we rely on the \textsc{Boost.Geometry} module supplied by the Boost C++ Libraries \cite{BoostLibraries} for the construction and manipulation of R$^*$-trees. \textsc{Boost} is a generic C++ library providing concepts,
geometry types, and general algorithms for solving, among others, problems in computational geometry. Its kernels are designed to be agnostic with respect to the number of space dimensions, coordinate systems, and types. Moreover, \textsc{Boost} provides the capability to perform spatial queries with polygonal shapes, therefore we envision its usage also when the underlying grid is already polytopal, such as with Voronoi tessellations. 

In the remainder of the paper, we will not make any distinction between R-trees and R$^*$-trees, as we always employ the latter.

\subsection{R-trees and finite element meshes}
\label{subsec:rtreefem}
\begin{figure}[ht]
  \centering
  \includegraphics[width=5cm]{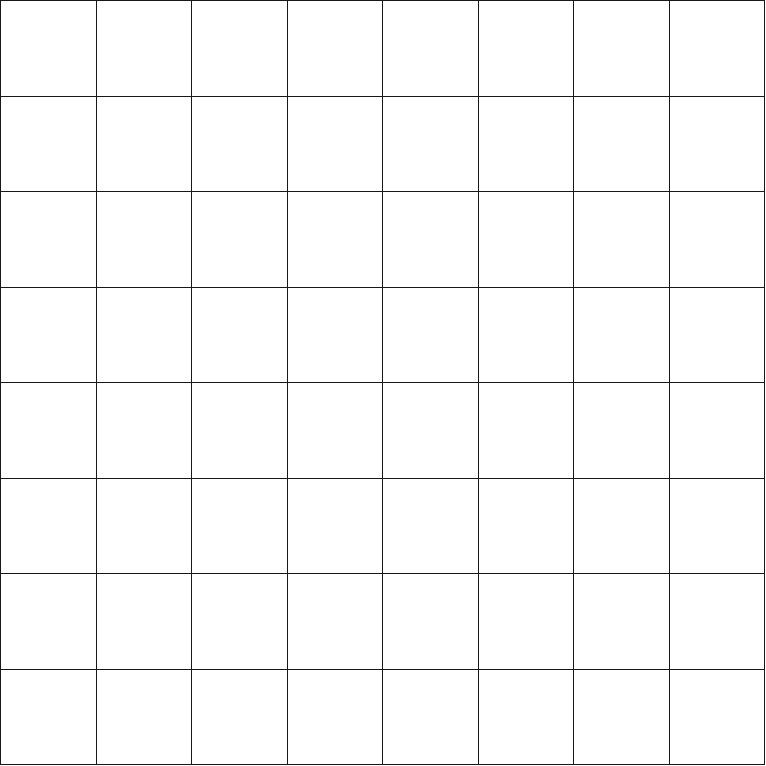}
  \hfill
  \includegraphics[width=5.7cm]{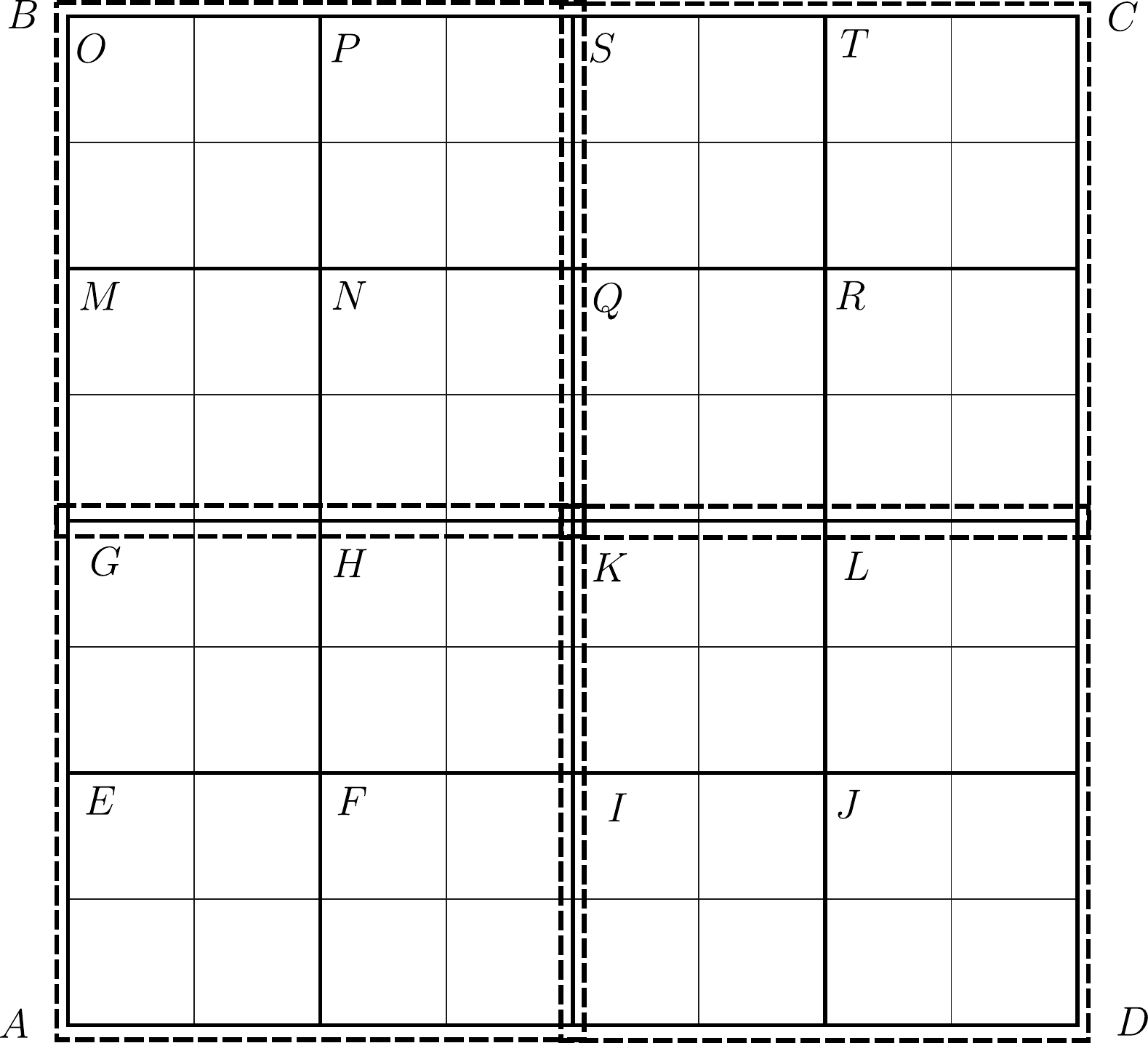}
  \caption{Left: $8 \times 8$ square grid obtained with $3$ uniform refinements of $\Omega=[0,1]^2$. Right: Minimal bounding rectangles generated on top of the grid elements of $\Omega$.}
  \label{fig:raw_mesh}
\end{figure}

\begin{figure}[ht]
  \centering
  \begin{forest}for tree={bplus, edge={->},l sep=1cm}
    [\mpn{A,B,C,D}
    [\mpn{E,F,G,H}[\mpn{,,,},multi edge=one][\mpn{,,,},multi edge=two][\mpn{,,,},multi edge=three][\mpn{,,,},multi edge=four],multi edge=one]
    [\mpn{I,J,K,L}[\mpn{,,,},multi edge=one][\mpn{,,,},multi edge=two][\mpn{,,,},multi edge=three][\mpn{,,,},multi edge=four],multi edge=four]
    ]
  \end{forest}
\caption{R-tree data structure for the $8 \times 8$ mesh example. For the sake of readability, only two internal nodes are shown. Each entry at the leaf
level is one mesh element $T$. Notice how each node stores exactly $4$ entries.}\label{fig:rtree_mesh_diagram}
\end{figure}
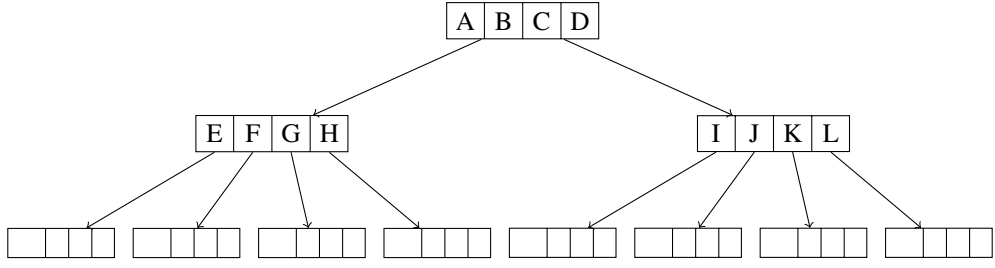

Here we describe the construction of the R-tree data structure associated with a given finite element mesh.
Let a domain $\Omega\subset{\mathbb R}^d$, $d\ge 1$ be given, which we identify with its partition into non-overlapping mesh elements $T \in \Omega$  covering $\Omega$.
The construction of the R-tree data structure for  $\Omega$ is initialised with the construction of the minimal bounding rectangle for each mesh element. 
We denote by $\{\text{MBR}(T_i)\}_{i=1}^N$, with $N$ the cardinality of $\Omega$, the resulting collection of MBRs. To give a simple yet practical example, consider the discretization of the unit square $[0,1]^2$ with the
$8 \times 8$ square grid in Figure \ref{fig:raw_mesh} (left). In this particular instance, for every $T \in \Omega$, it holds that $\text{MBR}(T) \equiv T$. Once the container with all MBRs is stored, the R-tree is ready to be built. The generated
hierarchical structure is depicted in Figure \ref{fig:raw_mesh} (right). The root node, which constitutes the first level of the R-tree, has four entries with MBRs $A,B,C,D$ (dashed lines). Each one of the entries points to an internal node, which in
turn is composed of other four MBRs. For instance, looking at the entry associated with $A$, it can be seen that it has as a child the internal node with entries $E,F,G$ and $H$, each one composed of other four elements, which are leaves
of the R-tree (and coincide with the mesh elements). The same pattern applies to the other three blocks of the mesh. The resulting R-tree is a $(2,4)$ R-tree, according to Definition \ref{def:definition_rtree}. A schematic view of the
tree hierarchy is shown in Figure \ref{fig:rtree_mesh_diagram}.

%% file: agglomeration_rtree.tex
\section{Agglomeration based on R-trees}
\label{sec:rtree-agglomeration}

Here we describe the construction of an agglomerated mesh starting from a given finite element mesh and the associated  R-tree data structure.
Let the mesh $\Omega$ of cardinality $N$ be given and assume the associated R-tree of order $(m,M)$ has been constructed. We denote with $L$ the total number of levels of the R-tree and with
$N_l$ the set of nodes in level $l$, for $l \in \{1,\ldots,L\}$. Our agglomeration strategy depends on an input parameter $l \in \{1,\ldots,L\}$ which describes
the level to be employed to generate the final agglomerates. The basic idea consists of looping through the nodes $N_l$ and, for each node $n \in N_l$, descending recursively its children
until leaf nodes are reached. These leaves share the same ancestor node $n$ (on level $l$) and thus are \emph{agglomerated}. We store such elements in a vector $v[n]$, indexed by the
node $n$. This procedure is outlined in Algorithm~\ref{alg:r_tree_algo}.

\begin{algorithm}[H]
  \caption{Creation of agglomerates.}
  \begin{algorithmic}[1]
      \Require{R-tree $R$ of order $(m,M)$}
      \Require{$l \in \{1, \ldots,L\}$ target level.}
      \Ensure{$v$ vector s.t. $v[n]$ stores leaves associated to node $n$}
      \Statex
    \Function{ComputeAgglomerates}{$R,l$}
      \For{node $n$ in $N_l$}
        \State $v[n] \gets $ \Call{Extractleaves}{$l,n$}
      \EndFor
      \State \Return $v$
    \EndFunction
  \end{algorithmic}
  \label{alg:r_tree_algo}
\end{algorithm}

\begin{algorithm}[H]
  \caption{Recursive extraction of leaves from a node $n$ on level $l$.}
  \begin{algorithmic}[1]
      \Require{$l \in \{1, \ldots,L\}$ target level}
      \Require{$n \in N_l$}
      \Ensure{vector $v$ containing leaves which share the ancestor node $n$.}
      \Statex
      \Function{Extractleaves}{$l,n$}
      \If{$l=1$}
      \State $v.\Call{add}{n}$
      \State \Return $v$
      \EndIf
      \For{node $ch\_n$ in \Call{children}{$n$}}
        \State  $v$.\Call{add}{\Call{Extractleaves}{$l-1,ch\_n$}}
      \EndFor
      \State \Return $v$
      \EndFunction
    \end{algorithmic}
  \label{alg:recursive_extract}
\end{algorithm}
After the recursive visit of the children of a node $n$, a list of mesh elements can be stored in $v[n]$
and flagged appropriately for agglomeration. Since the R-tree data structure provides a spatial partition of mesh elements, each of which is
uniquely associated with a node, the traversal of all nodes on a level $l$ provides a partition of mesh elements into agglomerates. The overall procedure
can be summarized as follows:
\begin{enumerate}
  \item Compute $\{\text{MBR}(T_i)\}_i$ for $i=1,\ldots,N$,
  \item Build the R-tree data structure using $\{\text{MBR}(T_i)\}_i$ as described in Section \ref{sec:rtree},
  \item Choose one level $l \in \{1, \ldots,L \}$ and apply Algorithm \ref{alg:r_tree_algo},
  \item For each node $n$, flag together elements of $v[n]$.
\end{enumerate}
We point out that elements in $v[n]$ are usually mesh-like iterators, i.e. lightweight objects such as pointers that uniquely identify
elements of $\Omega$. As we will show in Section \ref{sec:performance}, the overall construction of the R-tree and the actual identification
of agglomerates turns out to be quite efficient also for large 3D meshes. Moreover, having in mind multilevel methods, we notice that Algorithm \ref{alg:r_tree_algo} can be employed
to generate sequences of \emph{nested} agglomerated meshes. As a consequence, such meshes can be used as a hierarchy in a multigrid algorithm, allowing the usage of simpler and much cheaper intergrid transfer operators, when compared to the non-nested case. Building efficient intergrid transfer operators for non-matching meshes
is far from trivial even for simple geometries and their construction constitutes a severe bottleneck in terms of computational efforts, becoming critical in 3D \cite{AntoniettiPennesi,HHOnonnestedMG}. We notice, however, that in the context of
Lagrangian Finite Elements on standard hexahedral or quadrilateral grids, a completely parallel and matrix-free implementation of the non-nested geometric multigrid method has been recently addressed in \cite{MGNonNested}.

%% file: validation.tex
\section{Validation}
\label{sec:validation}

To show the effectiveness of the R-tree based grid agglomeration strategy, we perform several experiments with different grid types and compare it with \textsc{METIS}~\cite{METIS}, a standard graph partitioning algorithm often used for grid agglomeration. Other approaches, not considered here, include
Machine Learning-enhanced agglomeration strategies such as the ones developed in~\cite{antonietti2024agglomeration,antonietti2024polytopalmeshagglomerationgeometrical}.

More in detail, we use the multilevel k-way partitioning algorithm implemented in METIS to perform a partition of the graph associated with the given grid. Each partition corresponds to one agglomerated element and the cardinality of the resulting grid is thus fixed a priori as the number of partitions of the graph. 
We recall that, instead,  the input parameter in R-tree-based agglomeration is the level of the tree used to extract the agglomerates, as described in Algorithm~\ref{alg:r_tree_algo}.
Therefore, we will write $\mathtt{extraction\_level}$ and $\mathtt{n\_partitions}$ for the parameters required by the
R-tree and METIS strategy, respectively. 
The value of parameter $M$ in \ref{def:definition_rtree} is set to $2^d$, being $d$ the space dimension, while $m$ is set to $\frac{M}{2}=2^{d-1}$.

In view of a fair comparison, we always employ METIS setting $\mathtt{n\_partitions}$ as the number of elements produced
by the R-tree-based algorithm with a given $\mathtt{extraction\_level}$. In this way, the grids we compare, although effectively different, are guaranteed to be made of the same number of elements.

An important part of this work has been the  development of the \textsc{C++} library \textsc{PolyDEAL}, using the well-established Finite Element library \textsc{deal.ii}~\cite{dealII95,dealIIdesign} as a third-party
library. \textsc{PolyDEAL} provides the building blocks for solving partial differential equations with polytopal discontinuous Galerkin methods. It is distributed and builds on the Message Passing Interface (MPI) communication model~\cite{mpi40}. 
Providing an implementation within an existent finite element library has several advantages. Most importantly, we seamlessly inherit many robust features readily available and well-tested. Among them, we mention here \textsc{p4est}~\cite{p4est} for mesh partitioning across several processors, and \textsc{Trilinos}~\cite{Trilinos} or \textsc{PETSc}~\cite{balay1998petsc}
as parallel linear algebra libraries which provide a large variety of solvers. The library can be compiled by following the instructions available at the maintained GitHub repository~\cite{PolyDeal}. %
The results presented 
in this section, as well as the grids used in all the examples, are also made available in the same repository.

\subsection{Validation on a set of grids}
\label{subsec:validation_rtree}
\begin{figure}[!htb]
  \centering
  \subfloat[$\Omega_1$]{%
  \label{fig:srtuctured_square_view}%
  \includegraphics[width=0.28\textwidth]{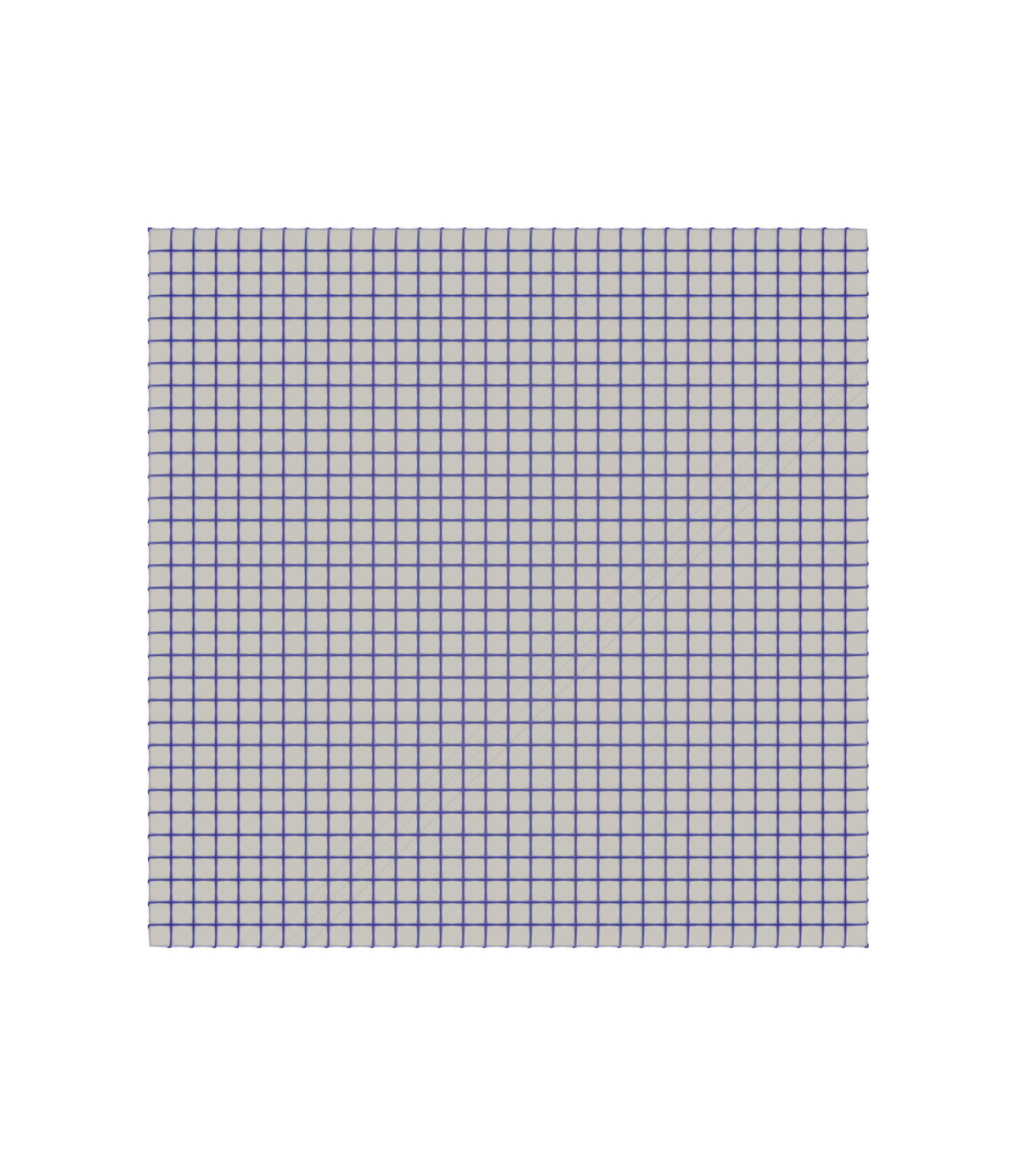}%
  }\hfill
  \subfloat[$\Omega_2$]{%
  \label{fig:structured_ball_view}%
  \includegraphics[width=0.28\textwidth]{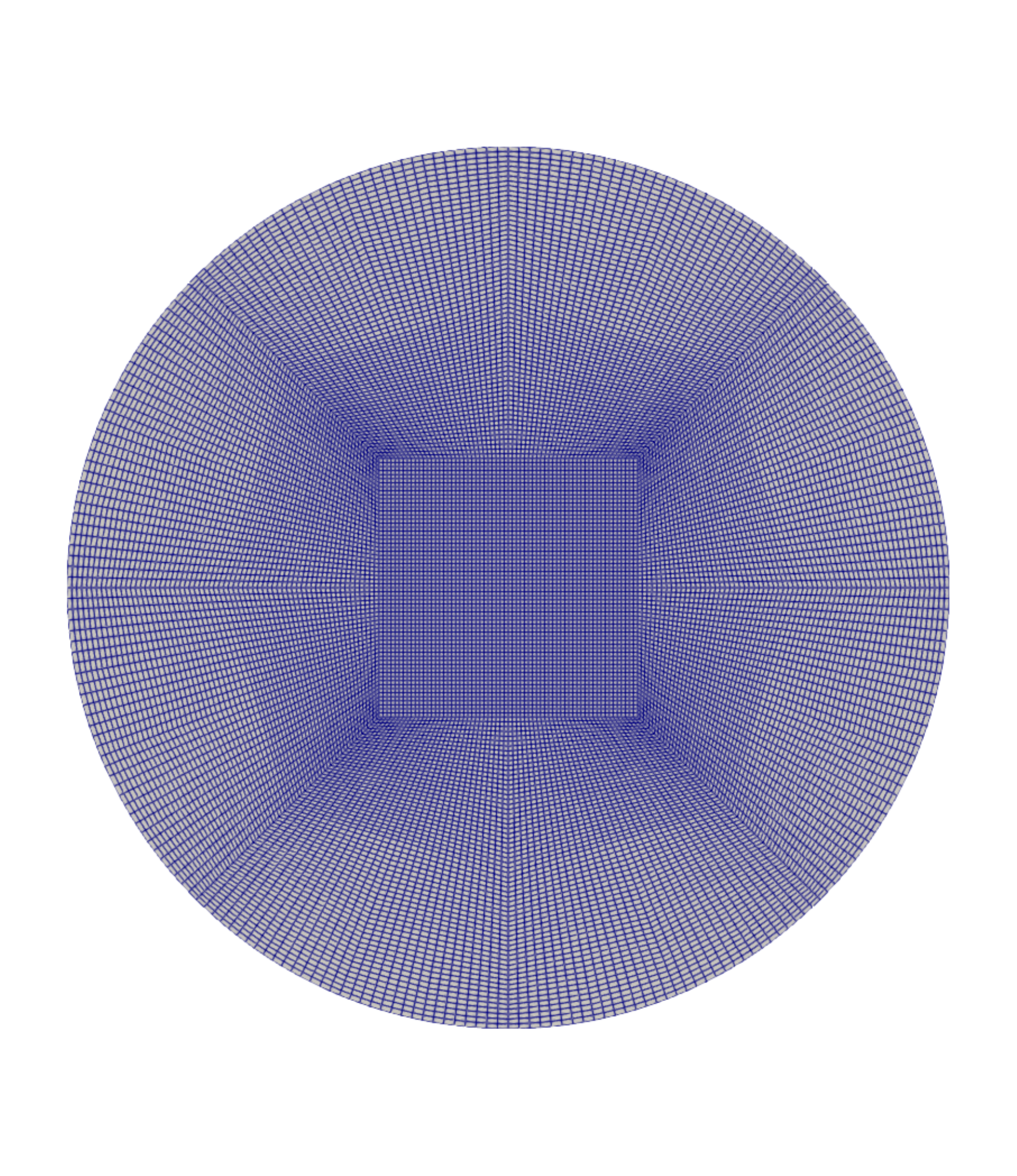}%
  }\hfill
  \subfloat[$\Omega_3$]{%
  \label{fig:unstructured_square_view}%
  \includegraphics[width=0.28\textwidth]{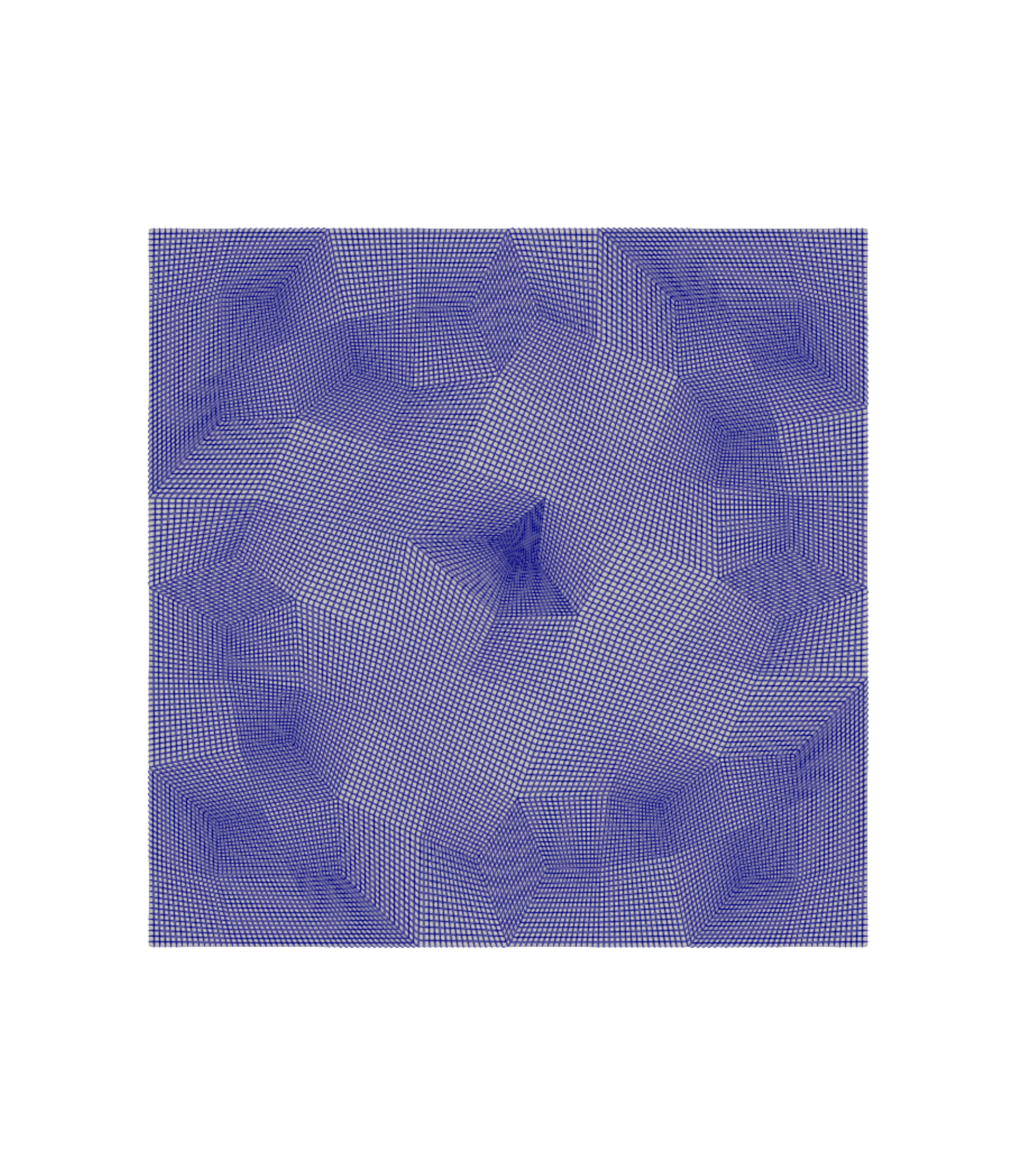}%
  }
  
  \medskip
  \subfloat[$\Omega_4$]{%
  \label{fig:voro_view}%
  \hspace{.4cm}
  \includegraphics[width=0.2\textwidth]{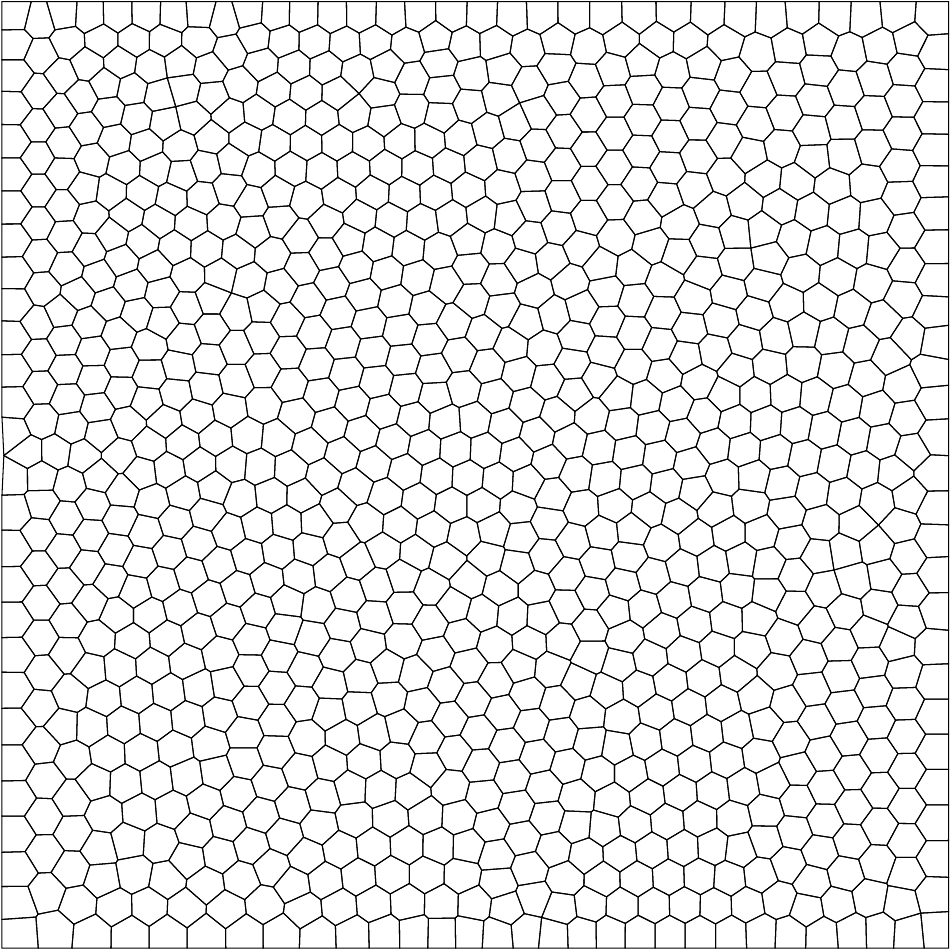}%
  }\hspace{1cm}%
  \subfloat[$\Omega_5$]{%
  \label{fig:cube_view}%
  \includegraphics[width=0.37\textwidth]{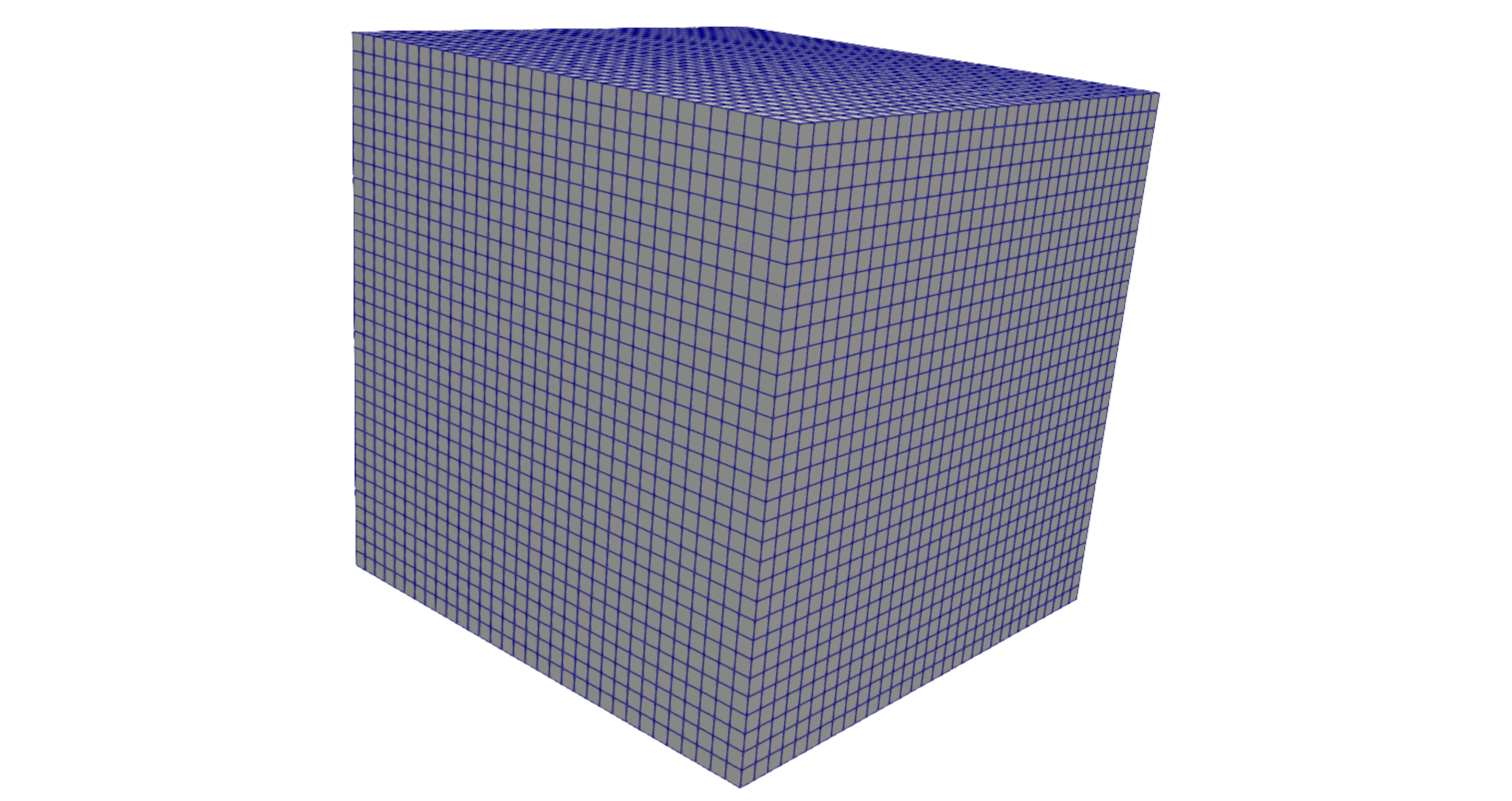}%
  }\hspace{3mm}%
  \subfloat[$\Omega_6$]{%
  \label{fig:piston_view}%
  \includegraphics[width=0.3\textwidth]{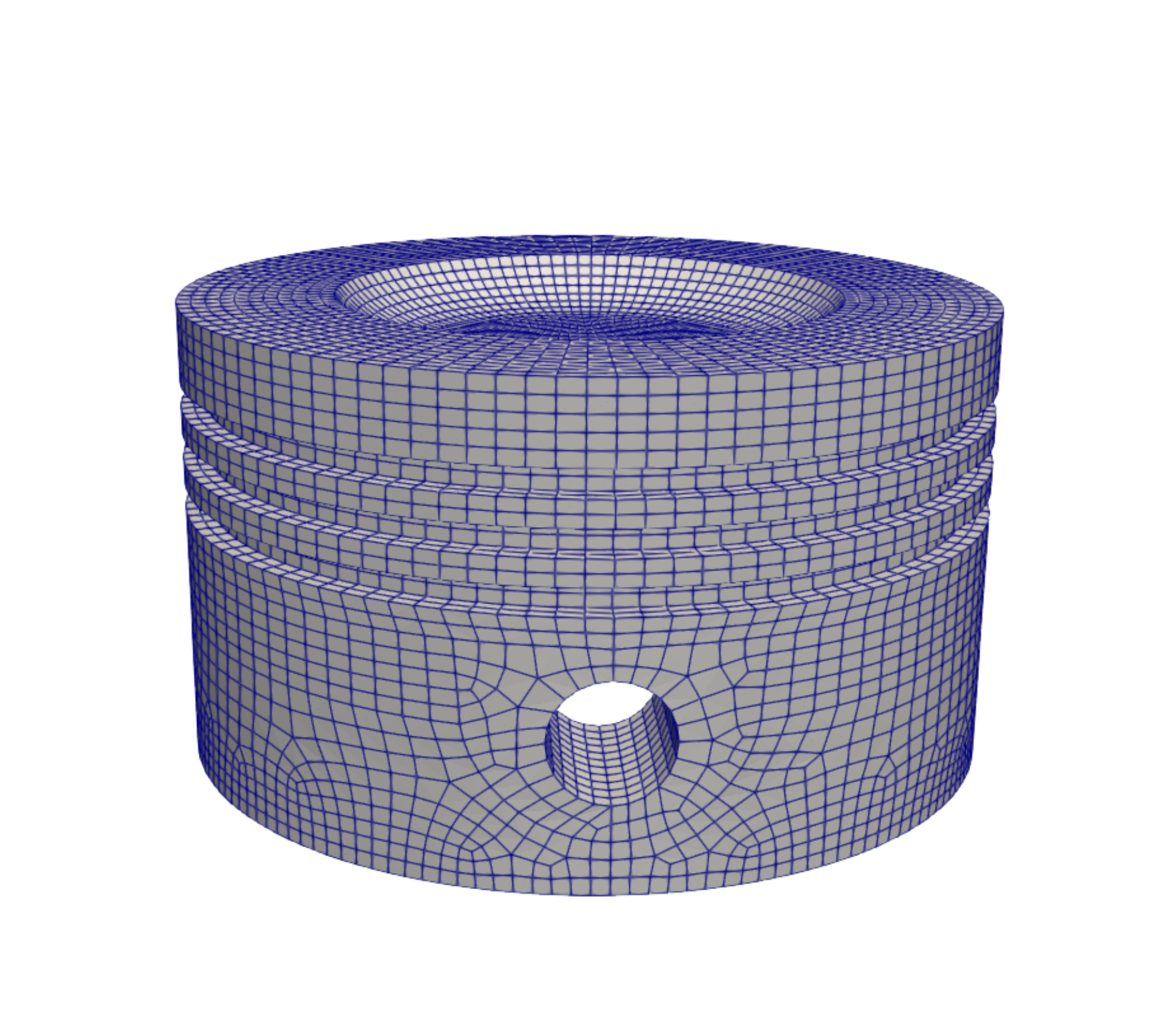}%
  }

  \medskip
  \hspace*{\fill} 
  \subfloat[$\Omega_7$]{%
  \label{fig:brain_view}%
  \includegraphics[width=0.29\textwidth]{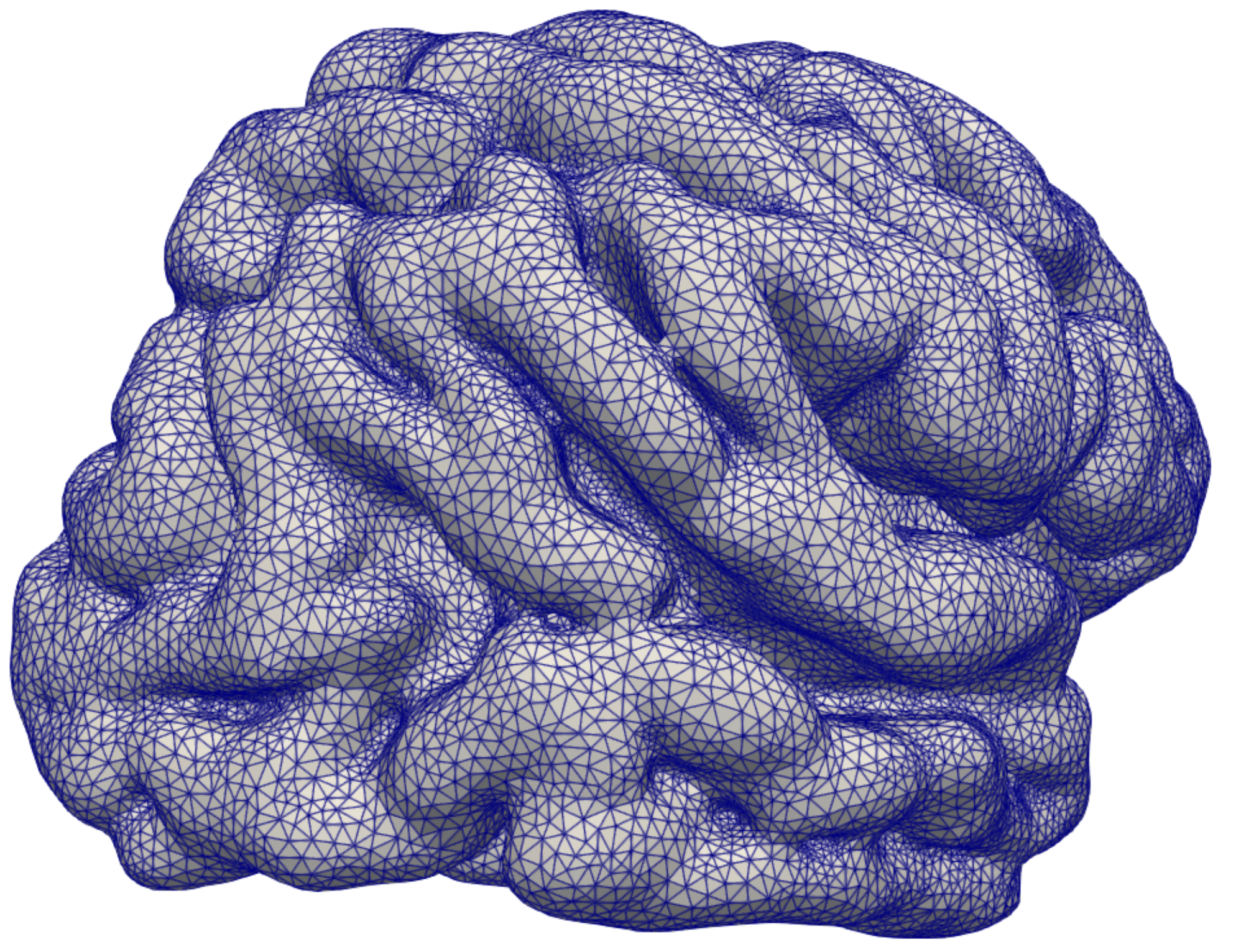}%
  }%
  \hspace{0.05\textwidth} 
  \subfloat[$\Omega_8$]{%
  \label{fig:liver_view}%
  \includegraphics[width=0.26\textwidth]{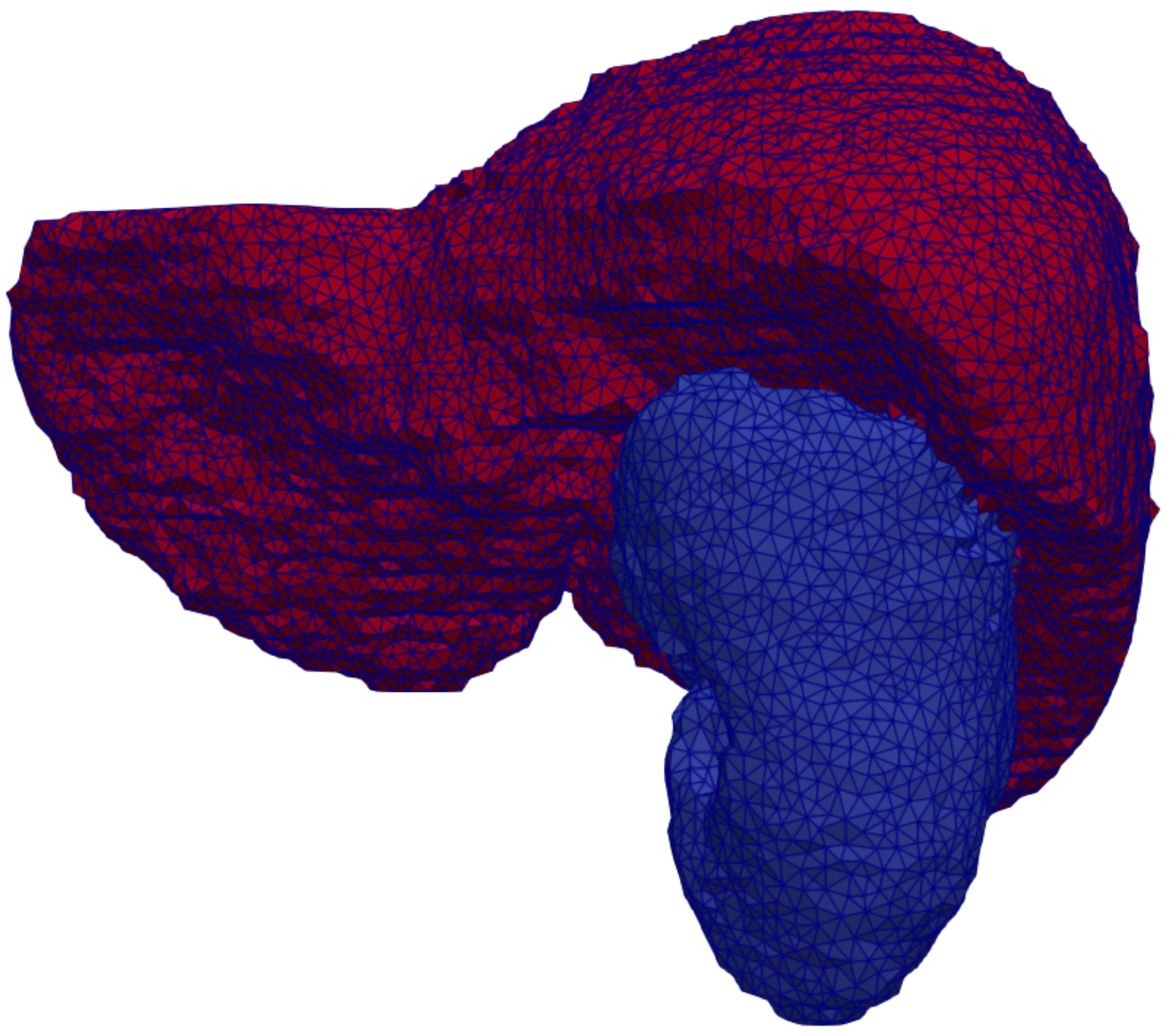}%
  }%
  \hspace*{\fill}

  \caption{The set of meshes used in the numerical experiments.}\label{fig:meshes}
  
  \end{figure}

To validate our methodology, we consider the following set of grids, sampled in Figure~\ref{fig:meshes}:
\begin{itemize}
  \item $\Omega_1$ structured partition of $(0,1)^2$,
  \item $\Omega_2$ structured partition of the open unit ball centered at the origin,
  \item $\Omega_3$ unstructured partition of $(0,1)^2$,
  \item $\Omega_4$: Centroidal Voronoi Tessellation (CVT) of $(0,1)^2$.
  \item $\Omega_5$ structured partition of  $(0,1)^3$,
  \item $\Omega_6$: CAD-modelled mesh of a piston,
  \item $\Omega_7$: mesh of a human brain,
  \item $\Omega_8$: mesh of a human liver,
\end{itemize}
We expect the circular domain and grid $\Omega_2$ to provide a less favorable case, for which higher overlap is to be expected due to the fixed, axis-aligned, orientation of the bounding boxes.
The Centroidal Voronoi Tessellation $\Omega_4$ was generated with the Polymesher package~\cite{polymesher} applying 100 iterations of Lloyd's algorithm to a random seeds Voronoi grid.
The grid $\Omega_6$ was generated after repairing and meshing the associated CAD $3$D model with the commercial
mesh generator \textsc{CUBIT}~\cite{Cubit}. $\Omega_7$ is a mesh of a brain created using the \textsc{SVMTK} library~\cite{bookMRI2FEM}, resulting in 634,472 tetrahedral elements. $\Omega_8$ is a mesh of a human liver obtained from a segmented medical image available in the \text{CGAL} library~\cite{cgal:rty-m3-24b}, with each voxel labeled according to its tissue type, resulting in a mesh composed of different domains where each subdomain corresponds to a specific tissue, totaling 284,201 tetrahedral elements.
The last two grids represent challenging scenarios and serve as excellent three-dimensional test cases for evaluating the performance and applicability of our algorithm. Furthermore, together with $\Omega_4$, they show that our approach seamlessly extends to elements other than quadrilaterals or hexahedra.

  \subsubsection{\emph{Test: Structured square}}
\begin{figure}[h]
    \centering
    \subfloat[METIS, $\mathtt{n\_partitions}=16$.]{%
    \includegraphics[width=.4\linewidth]{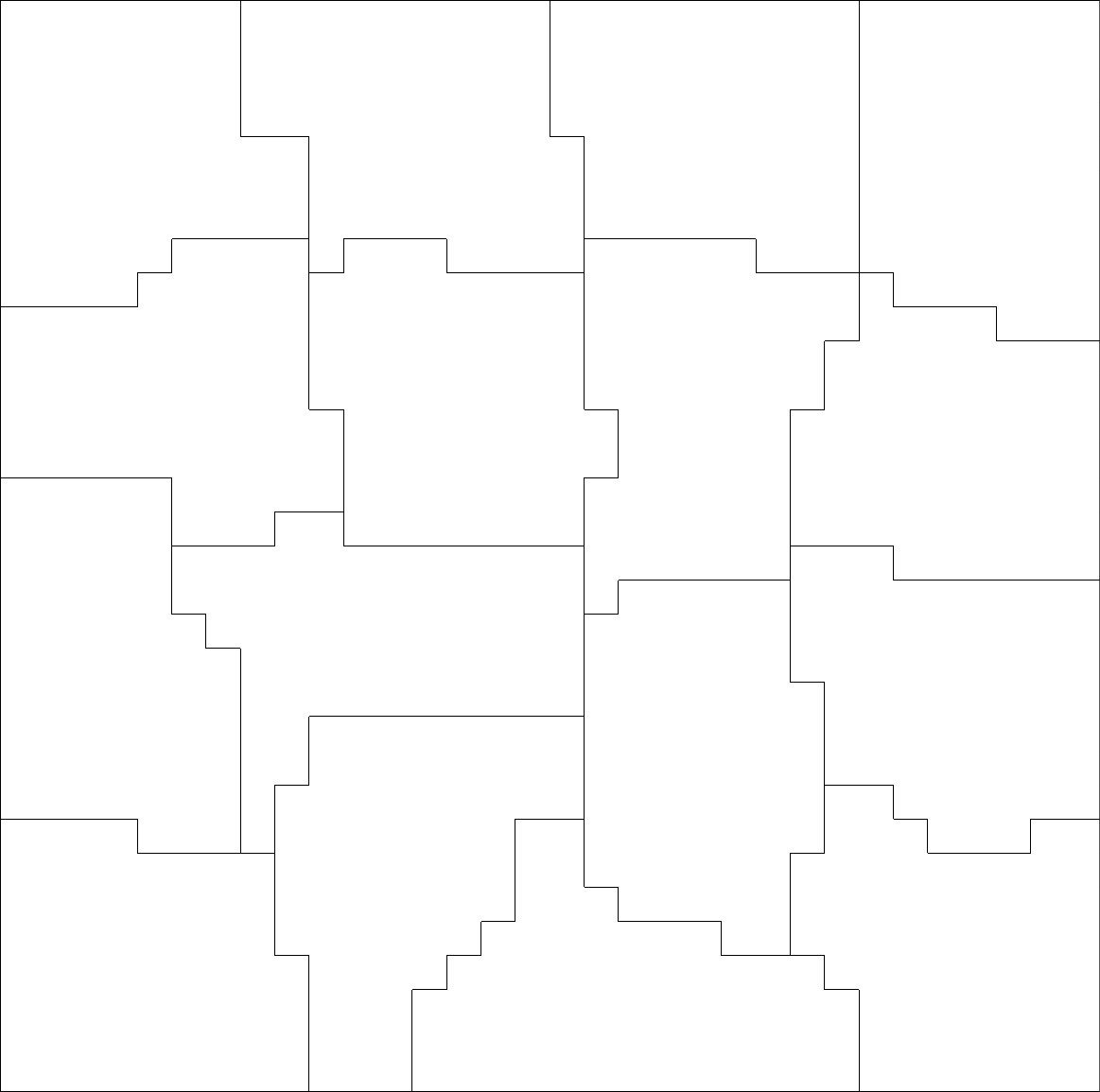}%
    \label{subfig:structured_square_METIS}%
    }\qquad\qquad
    \subfloat[R-tree, $\mathtt{extraction\_level}=2$.]{%
    \includegraphics[width=.4\linewidth]{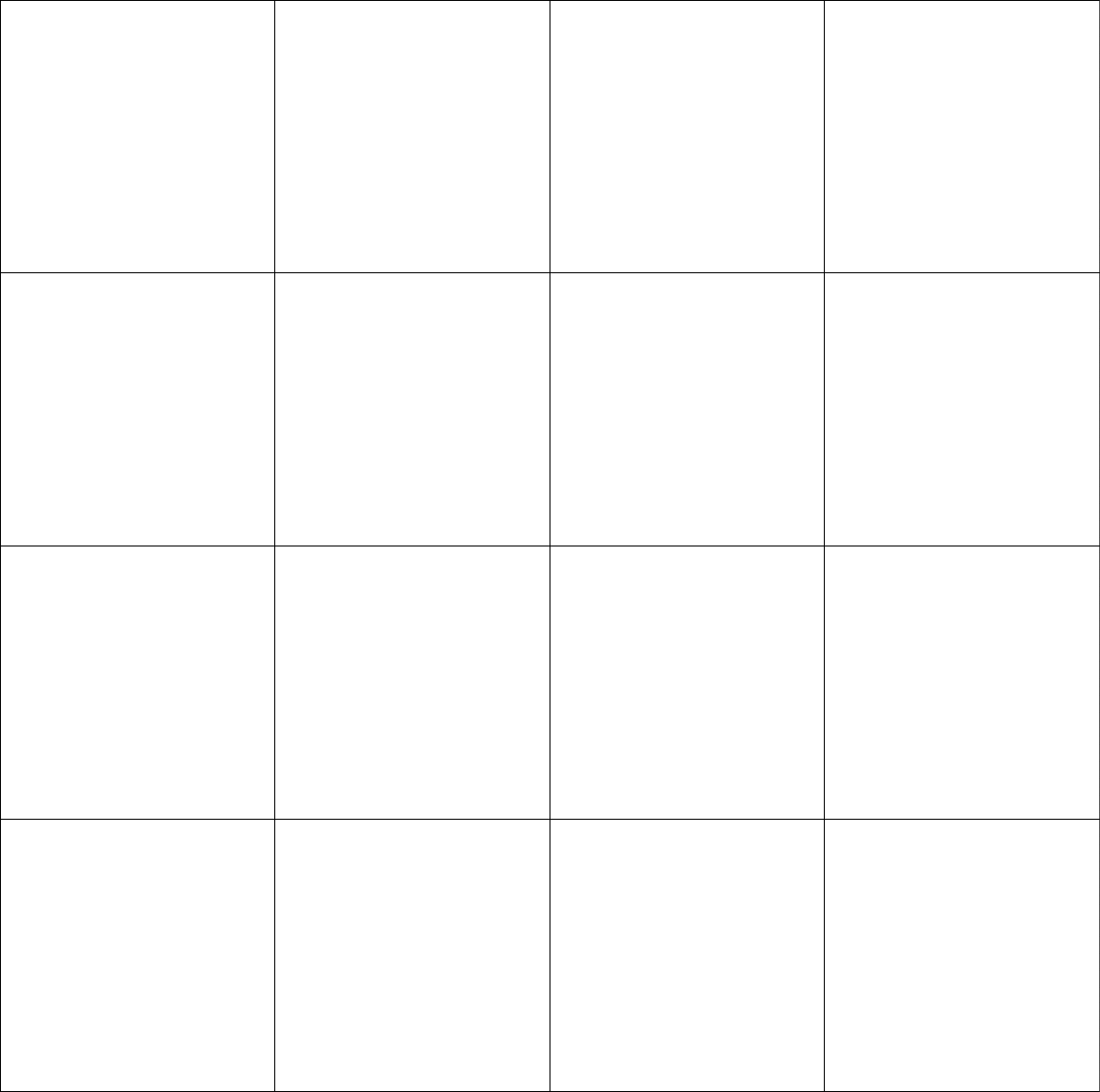}%
    \label{subfig:structured_square_Rtree}%
    }\vspace{1.00mm} 
    \centering
    \subfloat[METIS, $\mathtt{n\_partitions}=64$.]{%
    \includegraphics[width=.4\linewidth]{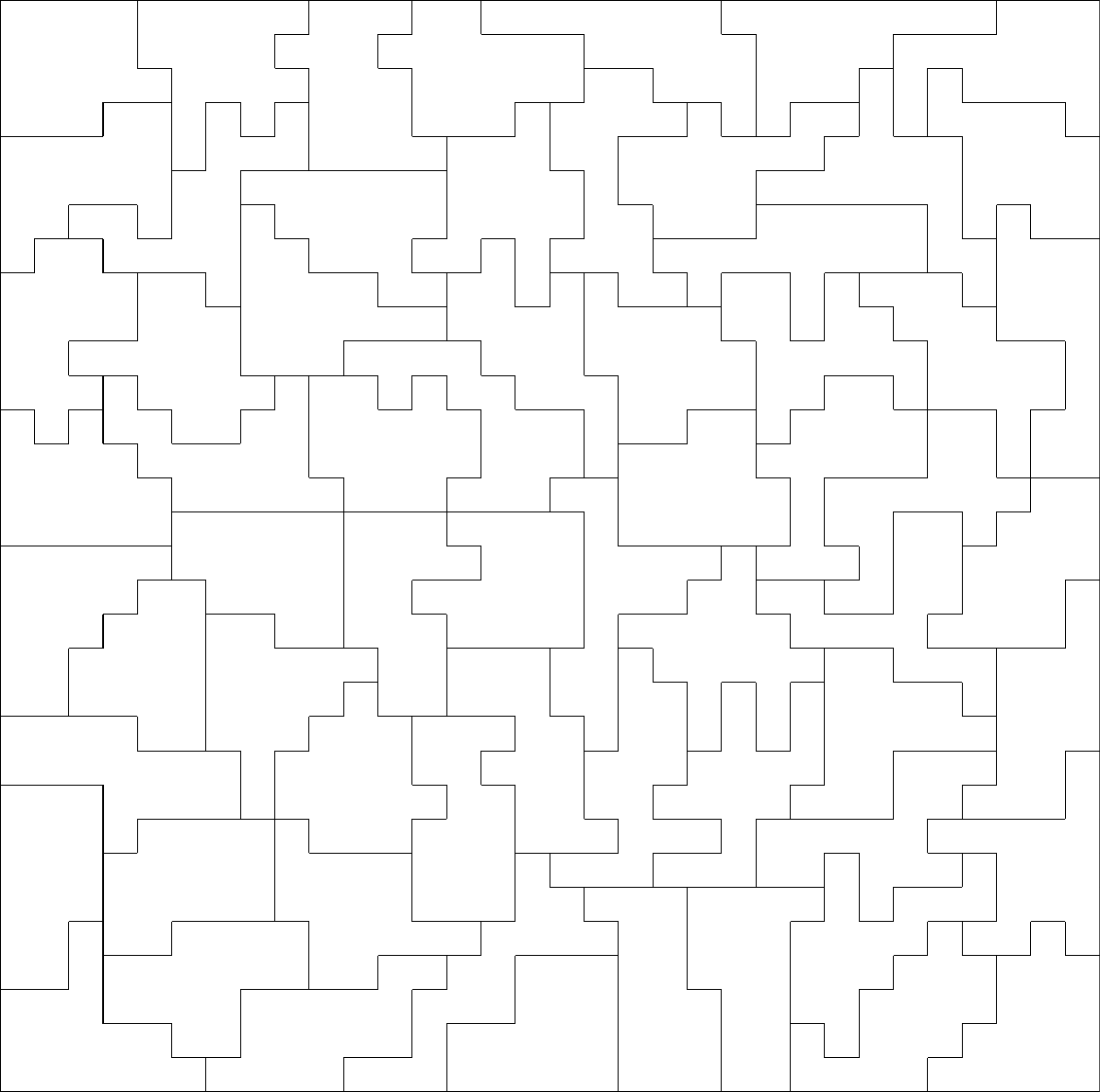}%
    \label{subfig:structured_square_level3_METIS}%
    }\qquad\qquad
    \subfloat[R-tree, $\mathtt{extraction\_level}=3$.]{%
    \includegraphics[width=.4\linewidth]{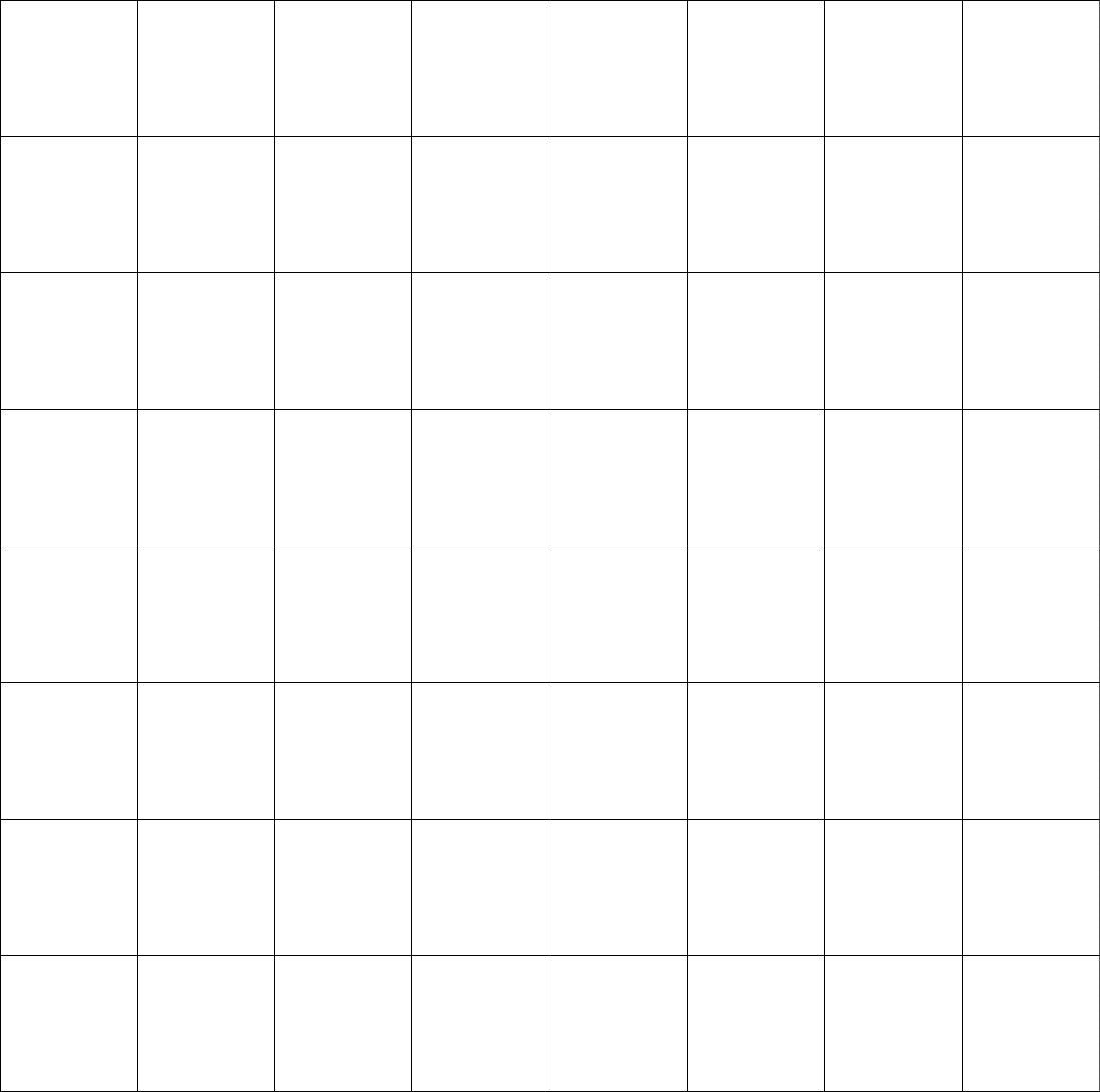}%
    \label{subfig:structured_square_level3_Rtree}%
    }
    \caption{Comparison between METIS and R-tree based agglomeration starting from the grid $\Omega_1$ displayed in Figure~\ref{fig:srtuctured_square_view}. Grids displayed in the same row always comprise the same number of elements.}
    \label{fig:structured_square}
  \end{figure}

  We fix as underlying mesh $\Omega_1$ the  
  $32 \times 32$ structured grid of squares. We report in Figure~\ref{fig:structured_square} the grids obtained by agglomeration of $\Omega_1$ using either METIS or the R-tree.
  With the extraction of the R-tree second level, we obtain the $4 \times 4$ square mesh reported in Figure~\ref{subfig:structured_square_Rtree}. The target number of mesh elements required by METIS is hence set to $16$. The resulting agglomerated elements are jagged, cf. Figure \ref{subfig:structured_square_METIS}. This is to be expected: METIS only processes the information coming from the graph topology of the mesh, hence the agglomerates are not supposed to preserve in any way the initial geometry, despite the fact that each polygon is made by $16$ sub-elements (as it is also the case for the R-tree). 
  
  We repeat the procedure by setting $\mathtt{extraction\_level}=3$. The R-tree naturally produces the $8 \times 8$ Cartesian mesh, exactly the subdivision one would get by globally refining the coarser mesh in Figure \ref{subfig:structured_square_Rtree}. We thus set $64$ as the number of target elements for METIS. Results are shown in the bottom plots in Figure \ref{subfig:structured_square_level3_METIS}. We observe that METIS produces polygons with even more jagged shapes, with many skinny and elongated elements, a particularly poor result considering that the underlying mesh is Cartesian. 
  
\subsubsection{\emph{Test: structured ball}}

\begin{figure}[h]
  \centering
  \subfloat[METIS, $\mathtt{n\_partitions}=20$.]{%
      \includegraphics[width=.4\linewidth]{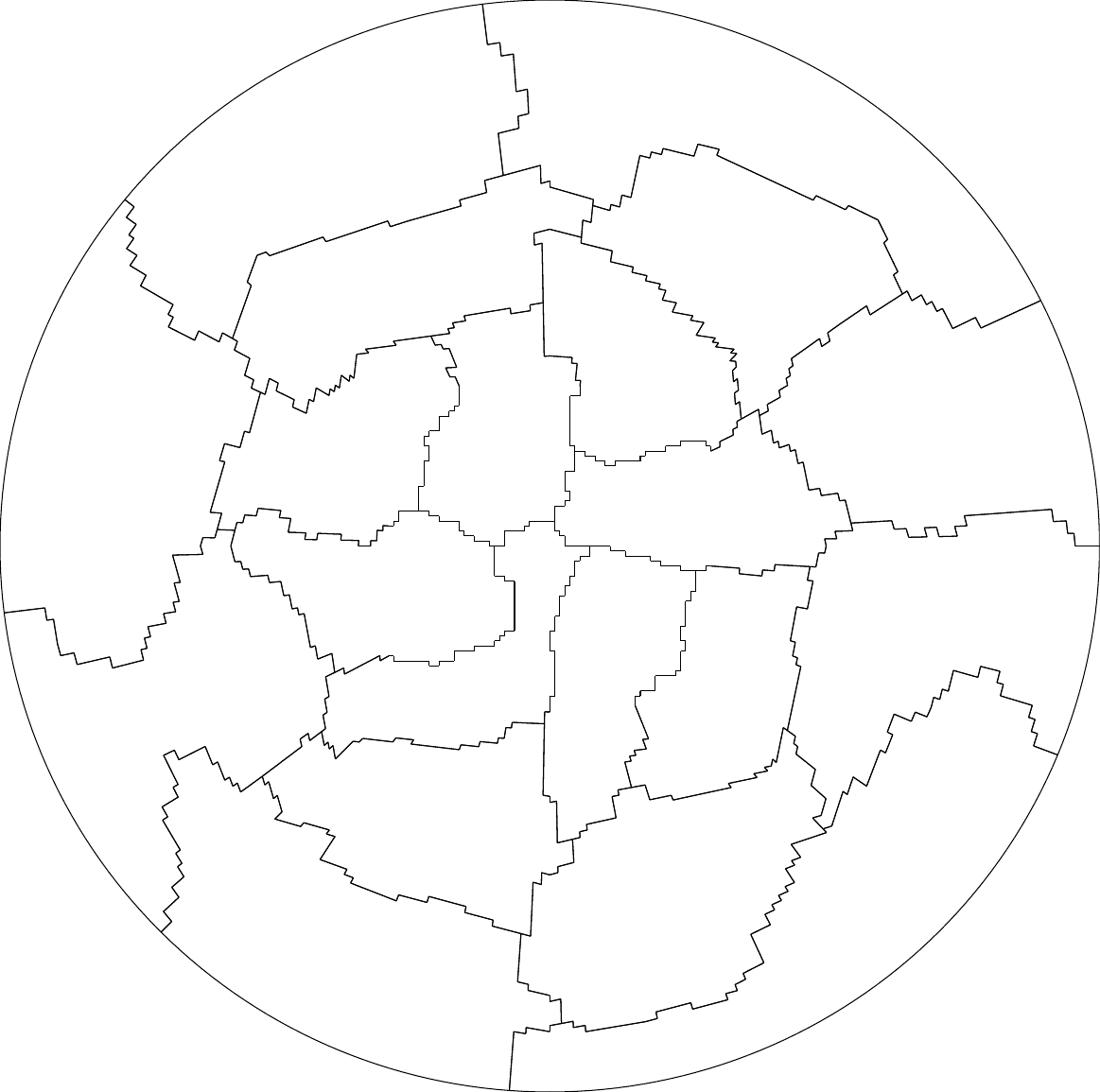}%
      \label{subfig:structured_ball_METIS}%
  }\qquad\qquad
  \subfloat[R-tree, $\mathtt{extraction\_level}=3$.]{%
      \includegraphics[width=.4\linewidth]{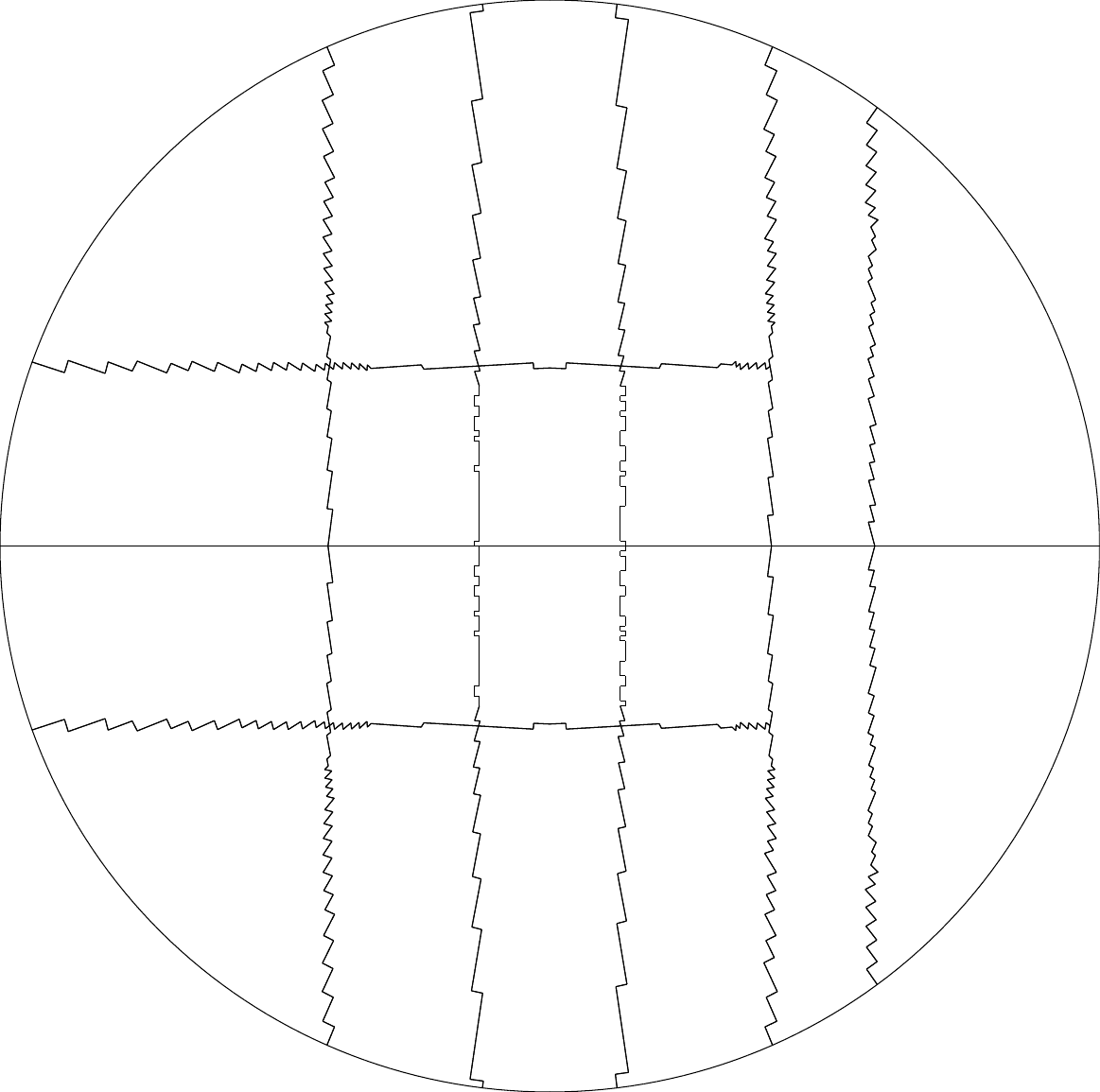}%
      \label{subfig:structured_ball_Rtree}%
  }\vspace{1.00mm} 
  \centering
  \subfloat[METIS, $\mathtt{n\_partitions}=80$.]{%
      \includegraphics[width=.4\linewidth]{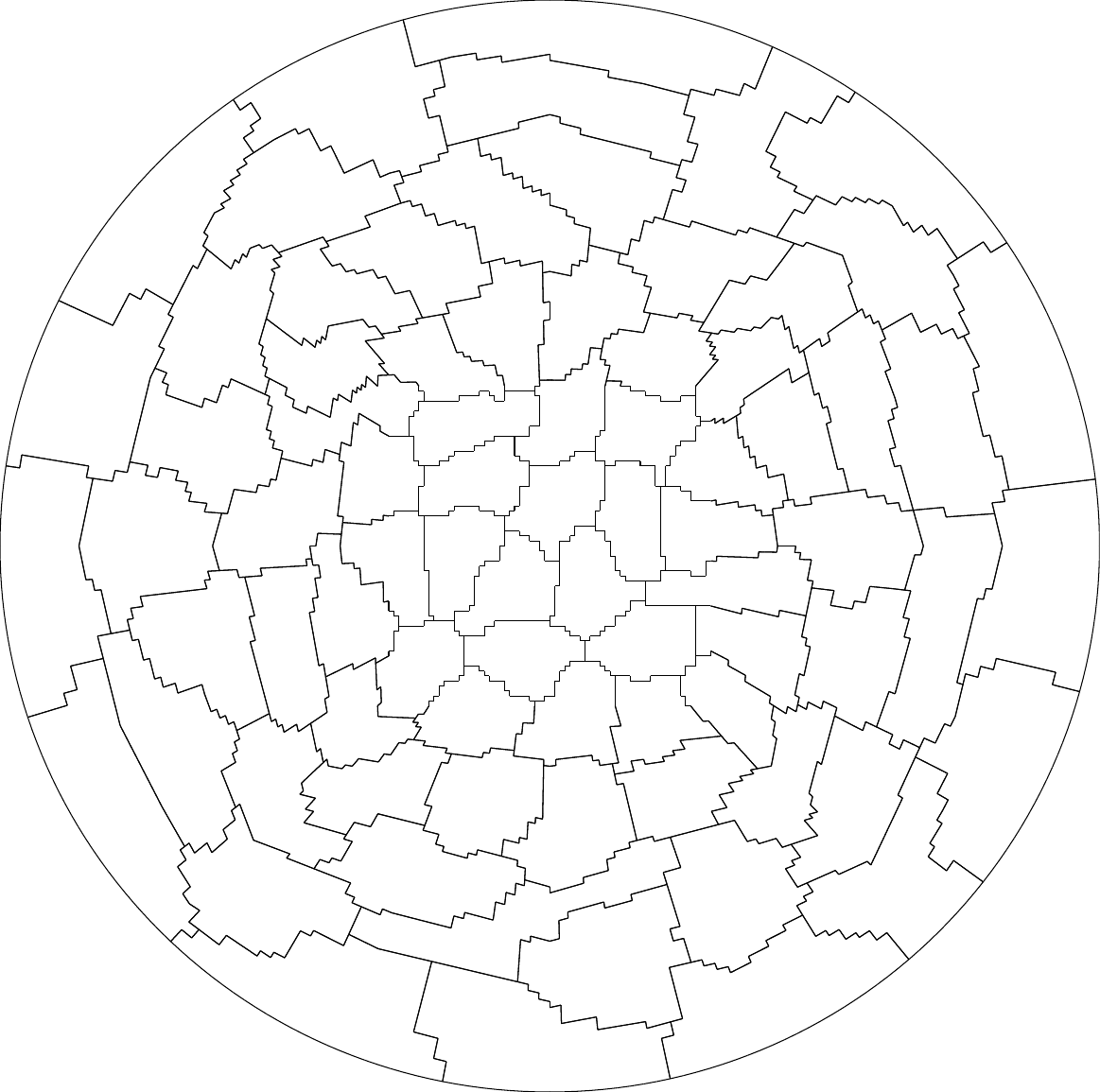}%
      \label{subfig:structured_ball_level4_METIS}%
  }\qquad\qquad
  \subfloat[R-tree, $\mathtt{extraction\_level}=4$.]{%
      \includegraphics[width=.4\linewidth]{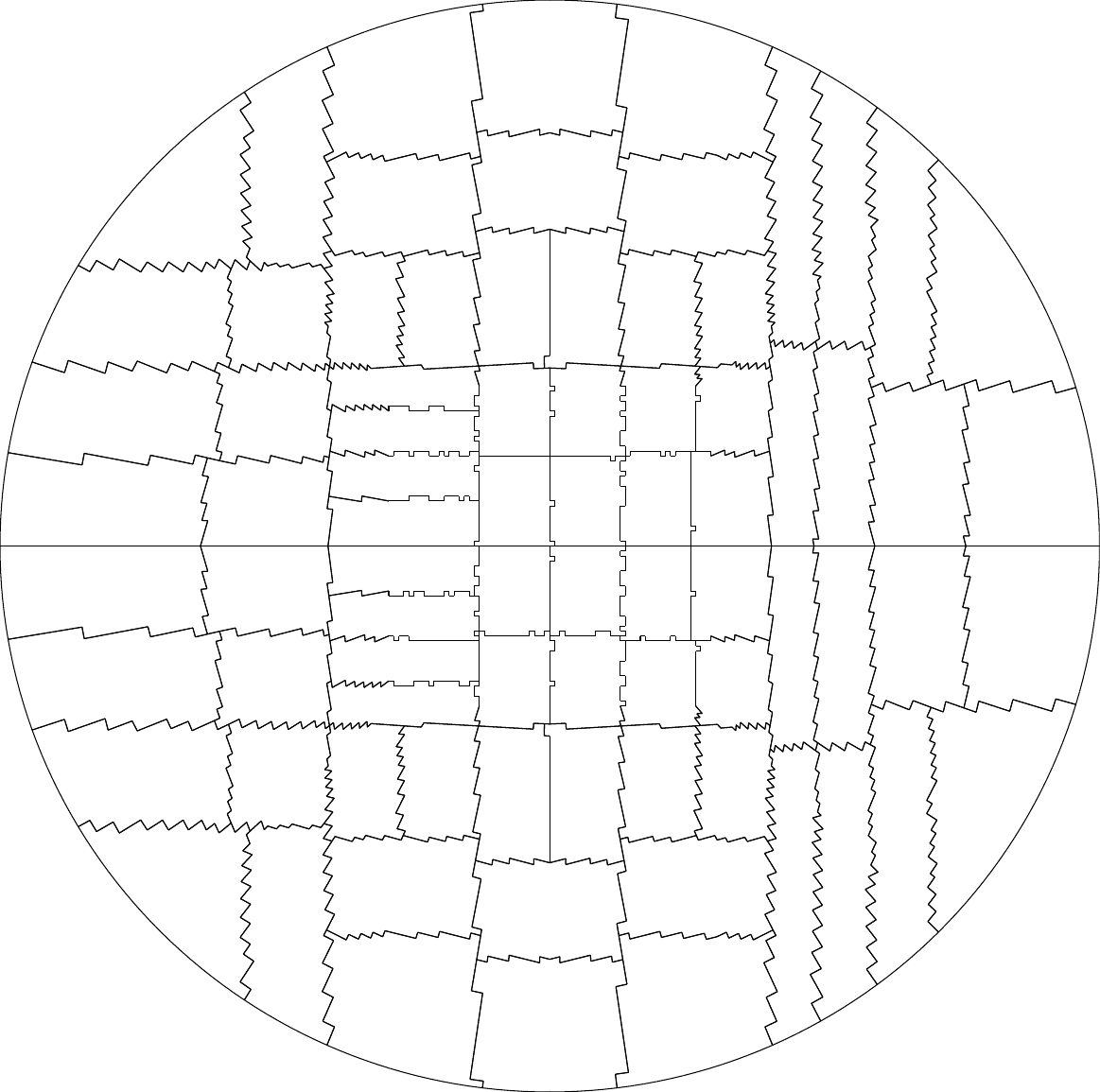}%
      \label{subfig:structured_ball_level4_Rtree}%
  }
  \caption{Comparison between METIS and R-tree based agglomeration of the grid $\Omega_2$ shown in Figure~\ref{fig:structured_ball_view}. Grids displayed in the same row always comprise the same number of elements.}
  \label{fig:structured_ball}
\end{figure}

We consider next as $\Omega_2$ the structured partition of a circle shown in Figure \ref{fig:structured_ball_view}. We remark that this grid is much finer than $\Omega_1$ as it is made of 20,480 elements. We start with extracting the third level of the R-tree, which gives $20$ agglomerates and accordingly set the target number of elements for METIS as $\mathtt{n\_partitions}=20$.
We observe from Figure \ref{subfig:structured_ball_Rtree} that the R-tree agglomerates conform to the rectangular shape of the respective bounding boxes, even though in this case the underlying quadrilateral mesh is not axis-aligned. METIS partitions, on the other hand, produce elements
with general shapes. Extracting the next level (i.e. setting $\mathtt{extraction\_level}=4$) of the hierarchy is equivalent to partitioning each agglomerate of $\mathtt{extraction\_level}=3$ into balanced sub-agglomerates, as shown in Figure \ref{subfig:structured_ball_level4_Rtree}, similarly to what happens with the square case shown previously. On the contrary, the two corresponding grids produced by METIS and shown in Figure~\ref{subfig:structured_ball_METIS} and Figure~\ref{subfig:structured_ball_level4_METIS} are completely unrelated.

\subsubsection{\emph{Test: unstructured square}}
\begin{figure}[h]
  \centering
  \subfloat[METIS, $\mathtt{n\_partitions}=91$.]{%
  \includegraphics[width=.4\linewidth]{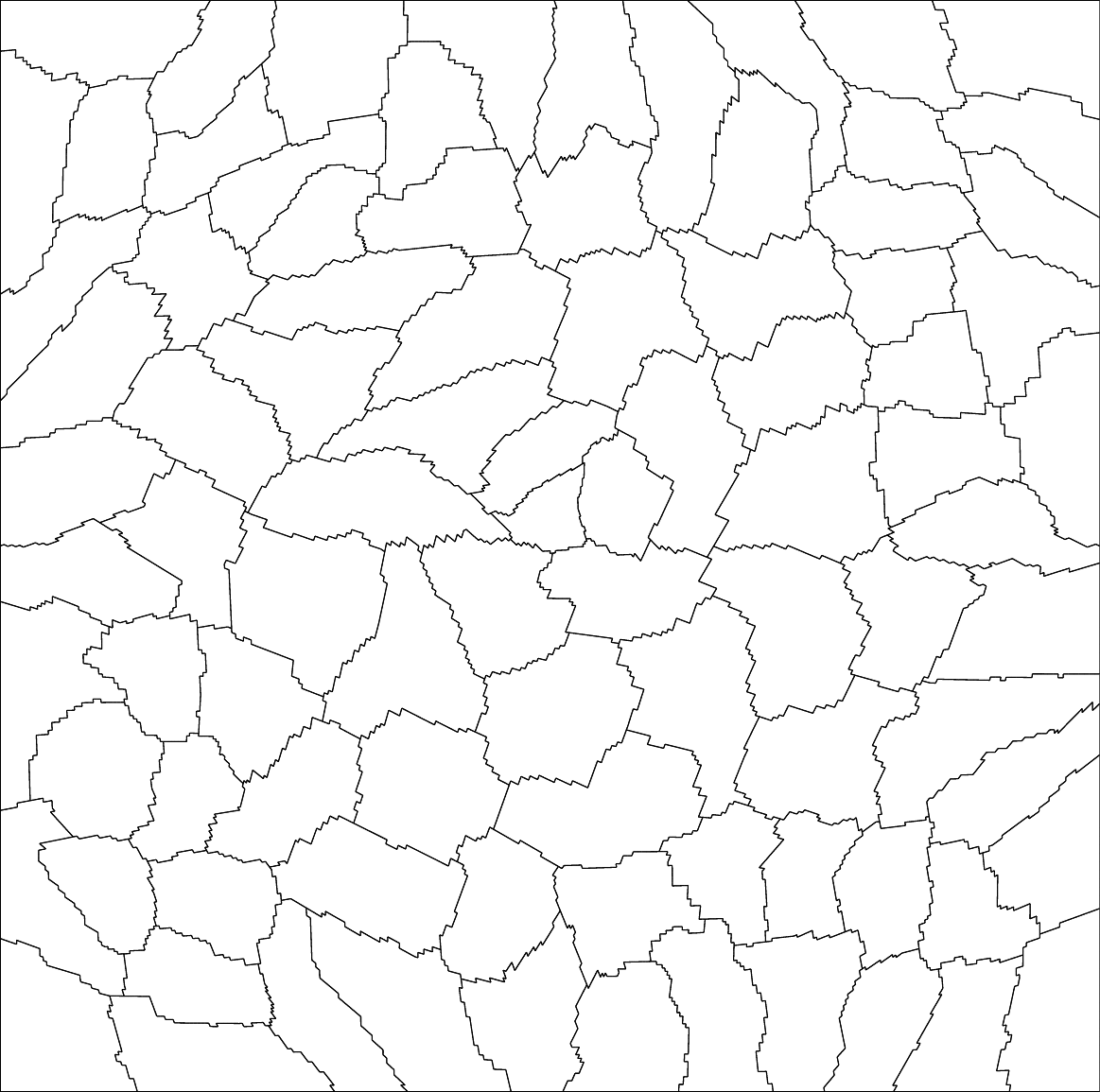}%
  \label{subfig:unstructured_square_level4_METIS}%
  }\qquad\qquad
  \subfloat[R-tree, $\mathtt{extraction\_level}=4$.]{%
      \includegraphics[width=.4\linewidth]{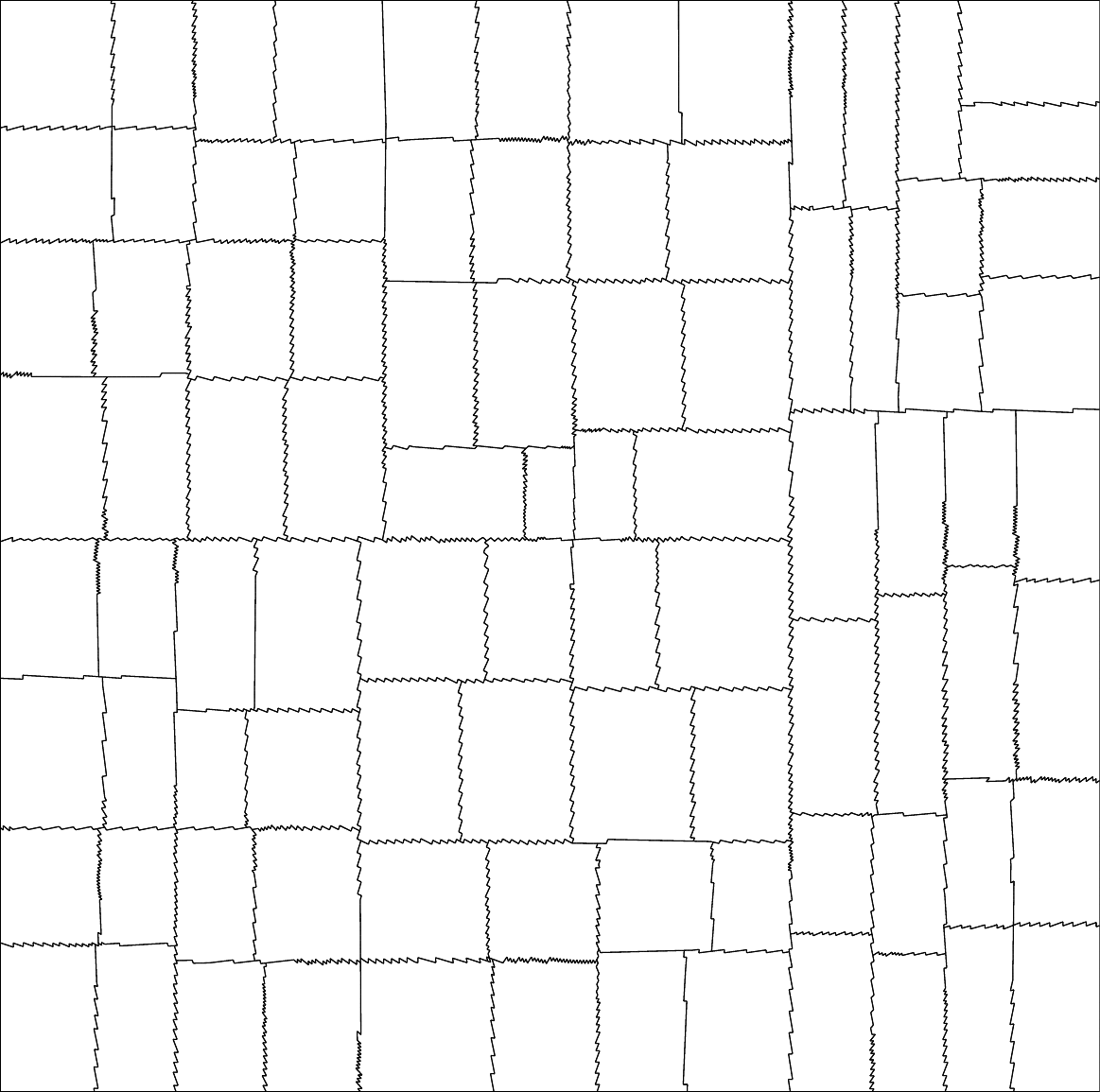}%
      \label{subfig:unstructured_square_level4_Rtree}%
  }\vspace{1.00mm} 
  \centering
  \subfloat[METIS, $\mathtt{n\_partitions}=364$.]{%
      \includegraphics[width=.4\linewidth]{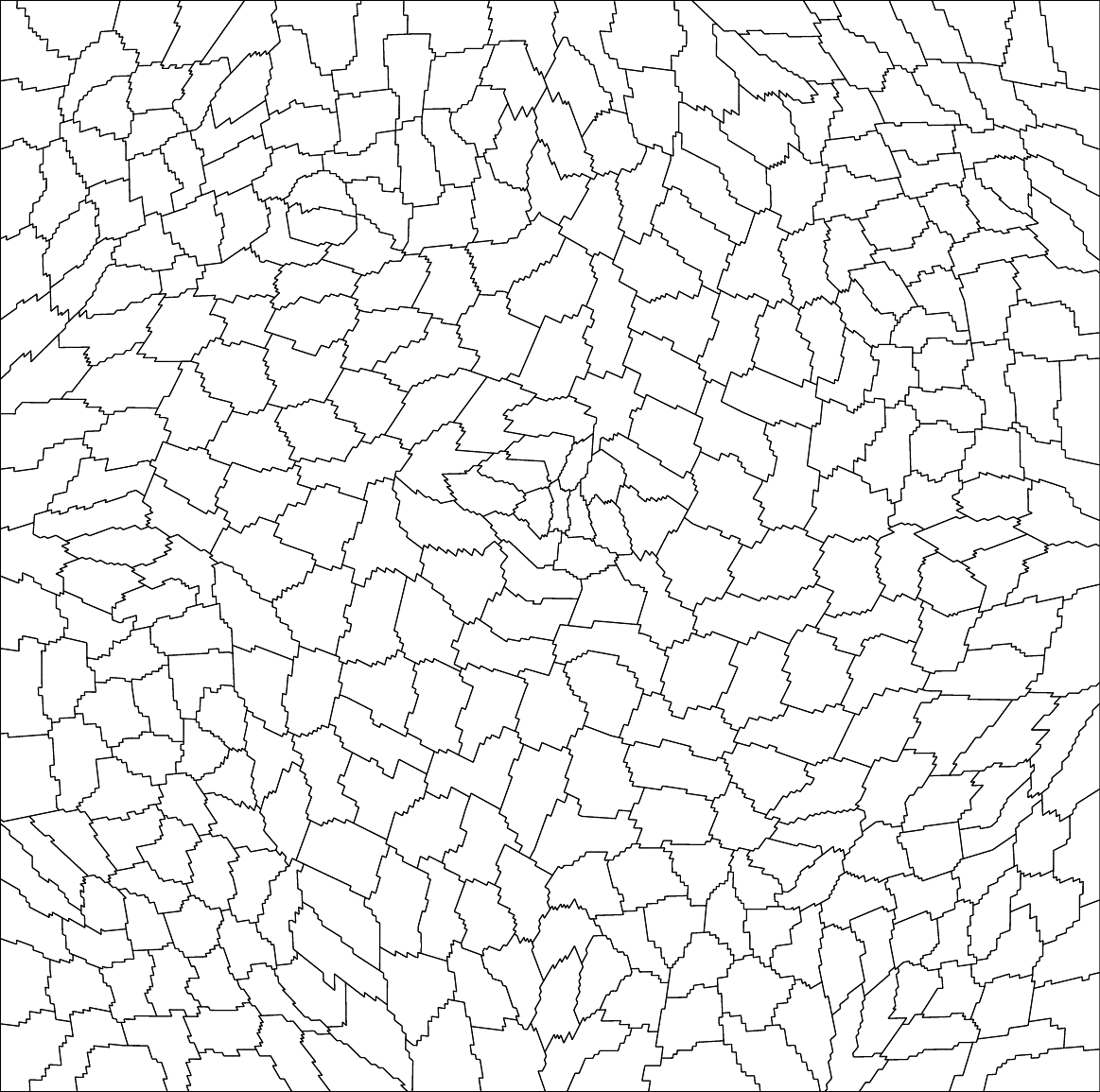}%
      \label{subfig:unstructured_square_level5_METIS}%
  }\qquad\qquad
  \subfloat[R-tree, $\mathtt{extraction\_level}=5$.]{%
      \includegraphics[width=.4\linewidth]{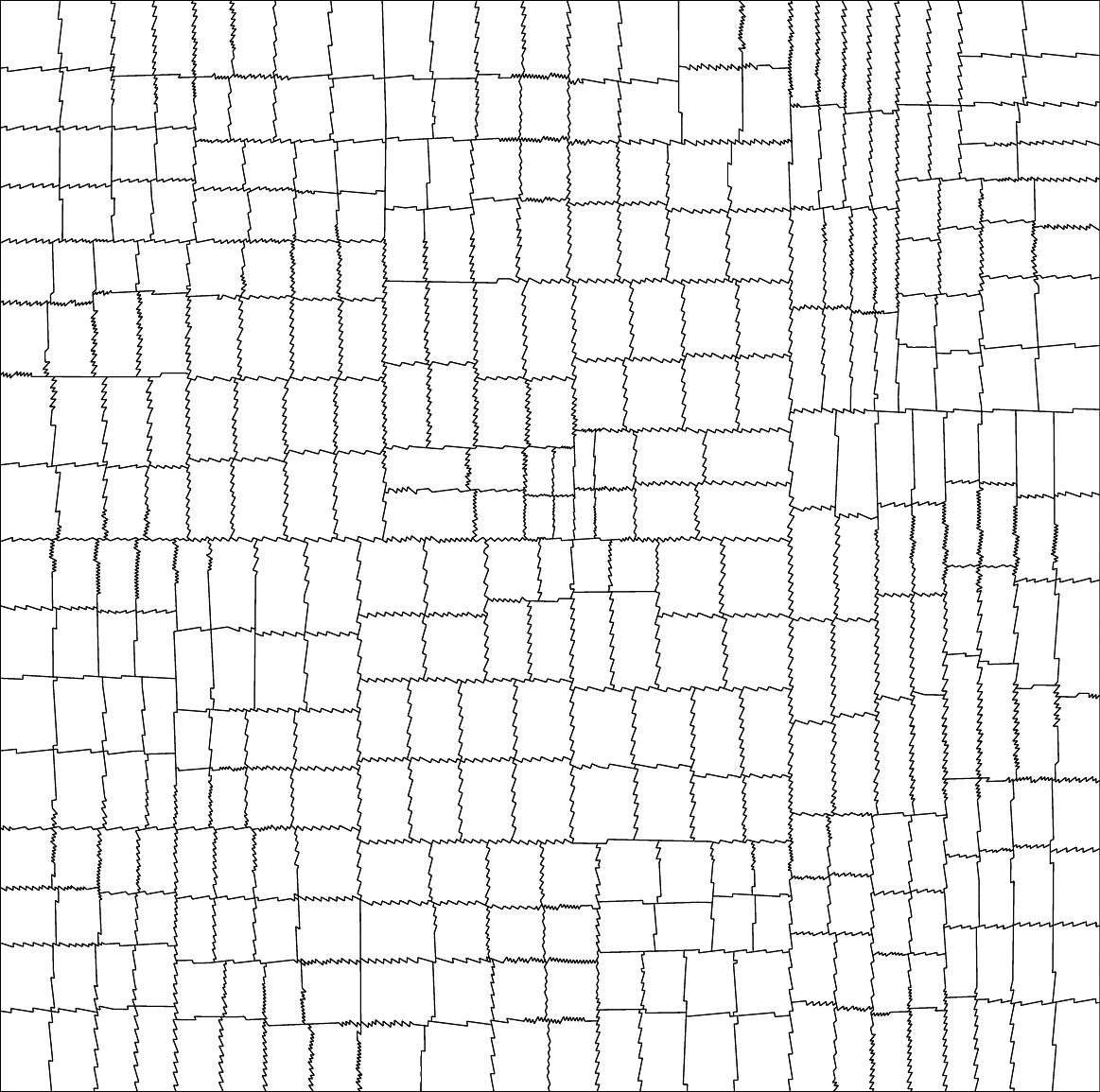}%
      \label{subfig:unstructured_square_level5_Rtree}%
  }
  \caption{Comparison between METIS and R-tree agglomeration starting from the grid $\Omega_3$ seen in Figure~\ref{fig:unstructured_square_view}. Grids displayed in the same row always comprise the same number of elements.}
  \label{fig:unstructured_square}
\end{figure}

The same procedure is repeated with the grid $\Omega_3$, composed of 93,184 non-uniform quadrilaterals. The main difference between the two examples
above relies on the fact that this grid is fully unstructured. With the R-tree, we extract first level $4$, which gives $91$ elements. Employing METIS with $\mathtt{n\_partitions}=91$ as input value, we
obtain the grid in Figure~\ref{subfig:unstructured_square_level4_METIS}. Here it is even more evident how the R-tree approach inherently gives rectangular-like shapes, in contrast to \textsc{METIS}. This is confirmed by the grids corresponding to $\mathtt{extraction\_level}=5$ and shown in Figures \ref{subfig:unstructured_square_level5_METIS} and~\ref{subfig:unstructured_square_level5_Rtree}.

\subsubsection{\emph{Test: Voronoi tessellation}}
\begin{figure}[h!]
  \centering
  \subfloat[\centering Coarse level R-tree, $\mathtt{extraction\_level}=2$.]{\includegraphics[width=.4\linewidth]{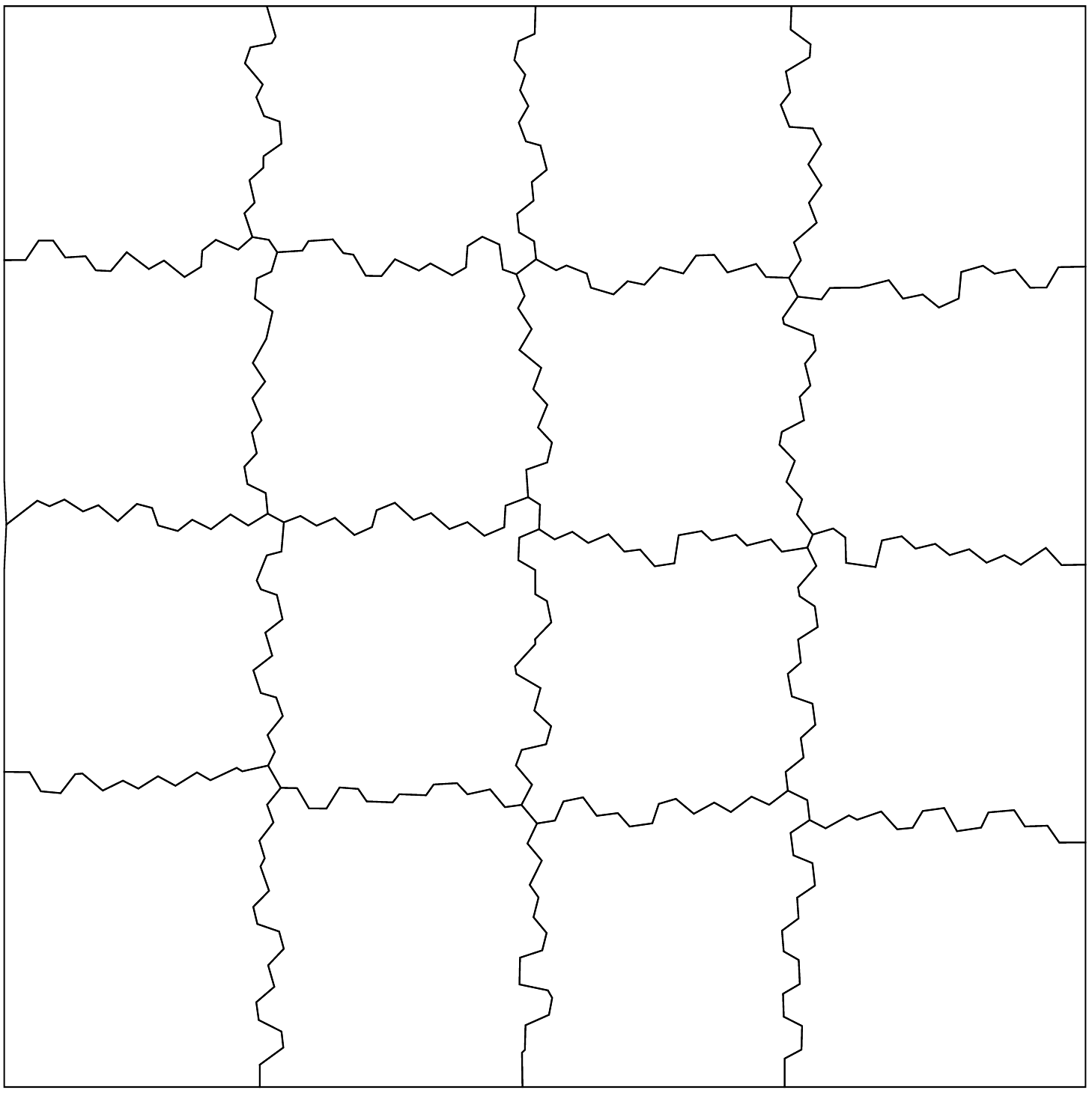}%
  \label{subfig:coarse_level}}
  \qquad\qquad
  \subfloat[\centering Fine level R-tree, $\mathtt{extraction\_level}=3$.]{\includegraphics[width=.4\linewidth]{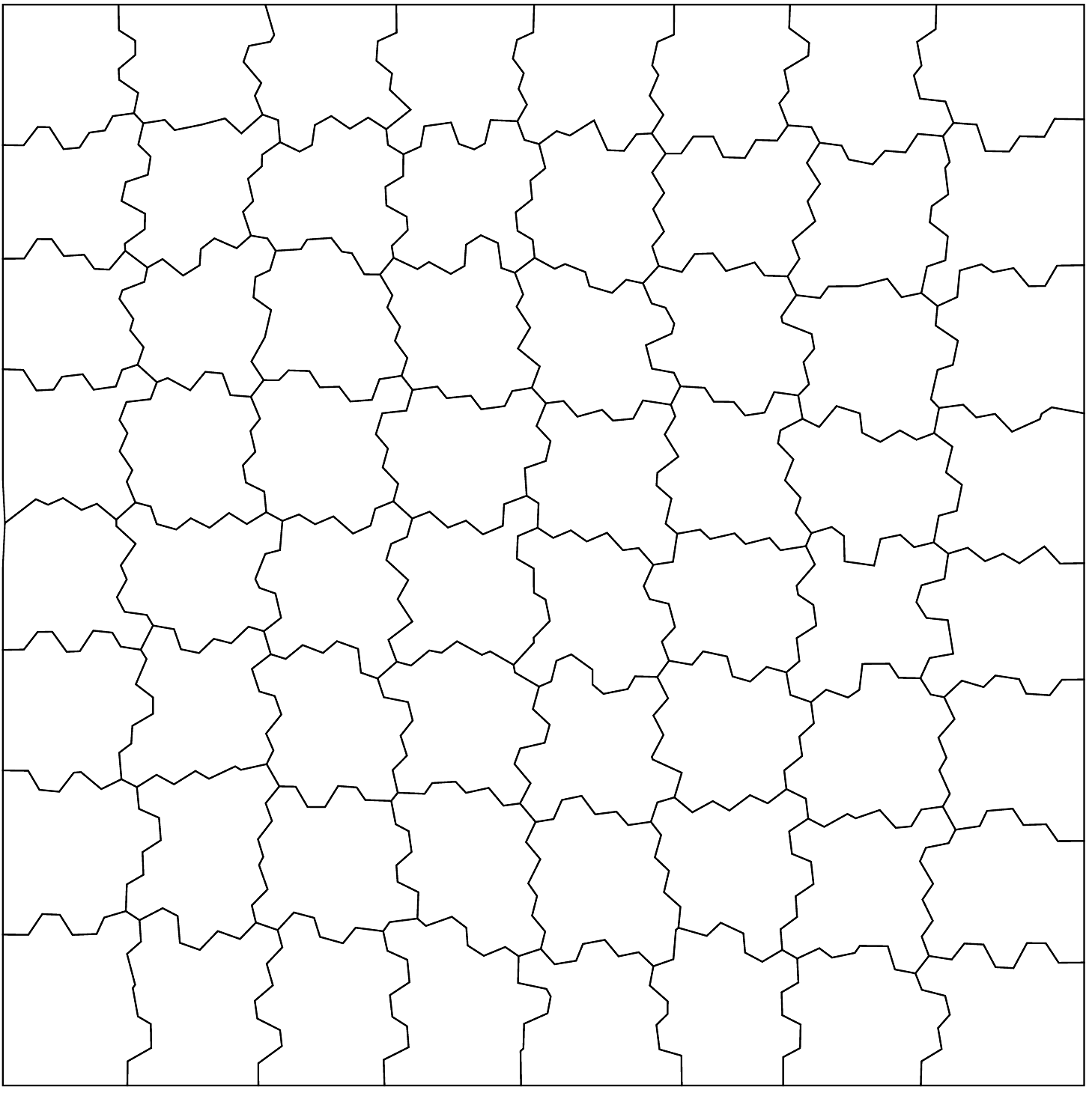}%
  \label{subfig:fine_level}}
  \caption{Two successive extraction levels for the CVT in Figure~\ref{fig:voro_view}.}%
  \label{fig:voronoi_agglo}%
\end{figure}

We consider the case of the grid $\Omega_4$, shown in Figure~\ref{fig:voro_view}, consisting of $1024$ polygonal elements. The application of the R-tree methodology to this grid produces a tree with $4$ levels, two of which are shown in Figure~\ref{fig:voronoi_agglo}. Once again, we observe that the resulting meshes are nested by construction. Moreover, the sequence of grids follows the classical square mesh pattern, with the level $l$ grid made of $2^l \times 2^l$ elements.  %
Based on the previous set of experiments, as well as from other comparisons such as the ones in~\cite{antonietti2024agglomeration}, we do not report here the results with METIS.

\subsubsection{\emph{Test: brain mesh}}
\begin{figure}[h]
  \centering
  \subfloat[\centering Agglomerates on the boundary of the brain.]{\includegraphics[width=.4\linewidth]{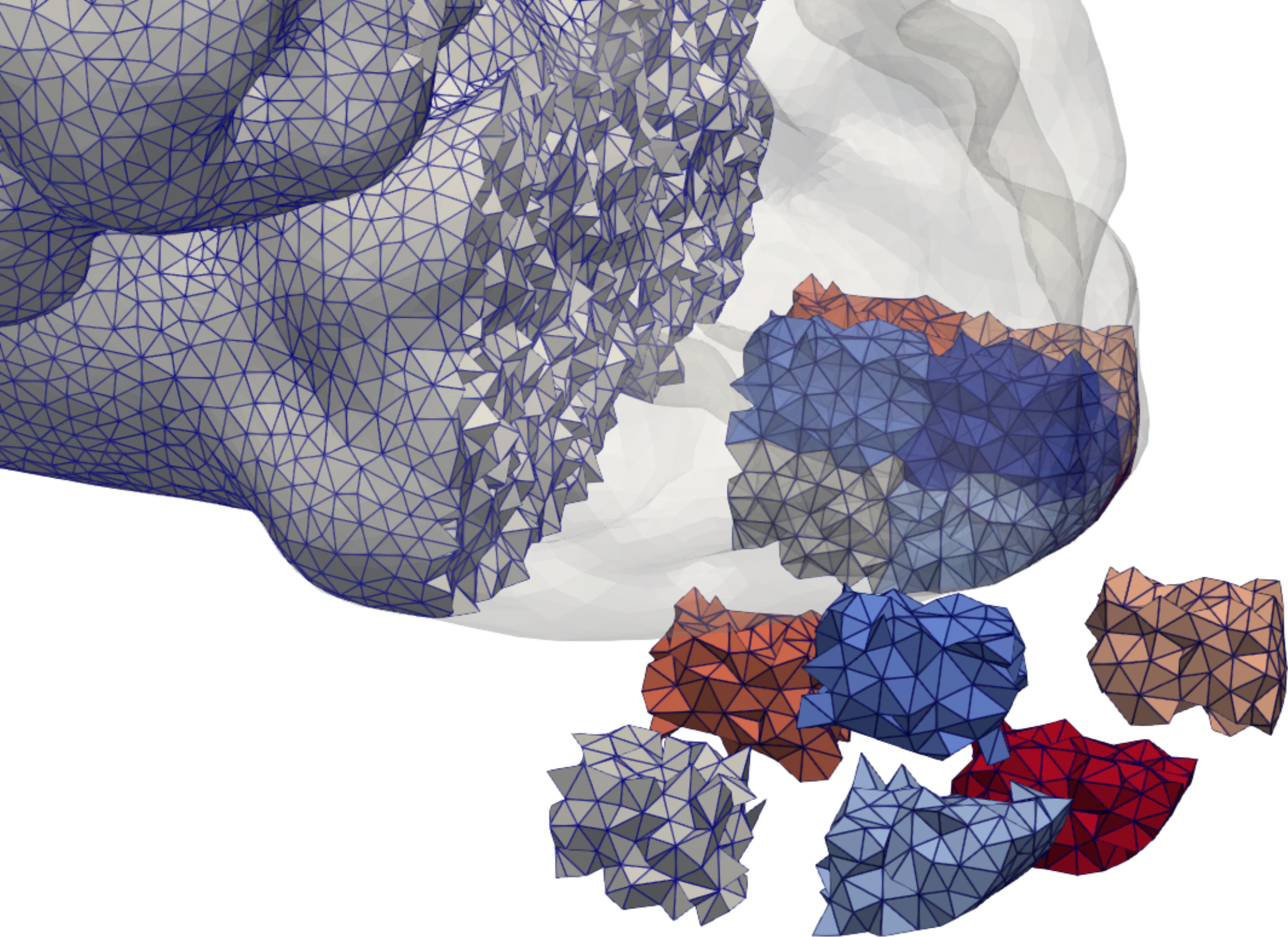}%
  \label{subfig:boundary_view}}
  \qquad
  \subfloat[\centering Agglomerates on the internal part of the brain.]{\includegraphics[width=.45\linewidth]{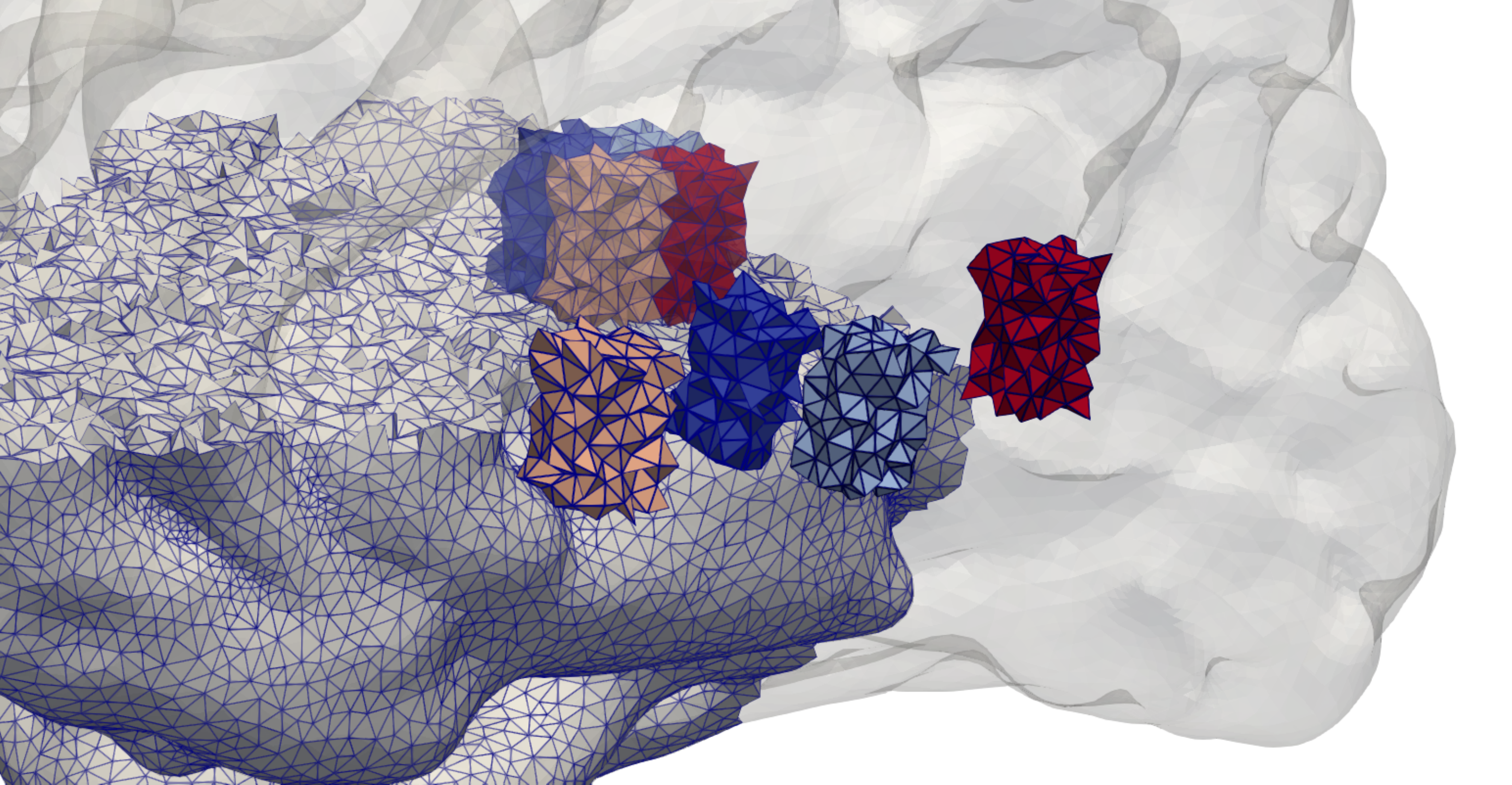}%
  \label{subfig:center_view}}
  \caption{Exploded view of agglomerates on different regions of the brain mesh $\Omega_7$. Each color maps to a different agglomerated element.}%
  \label{fig:brain_agglomerates_parent_child}%
\end{figure}

\begin{figure}[h]
  \centering
  \subfloat[\centering Parent agglomerate (in green).]{\includegraphics[width=.6\linewidth]{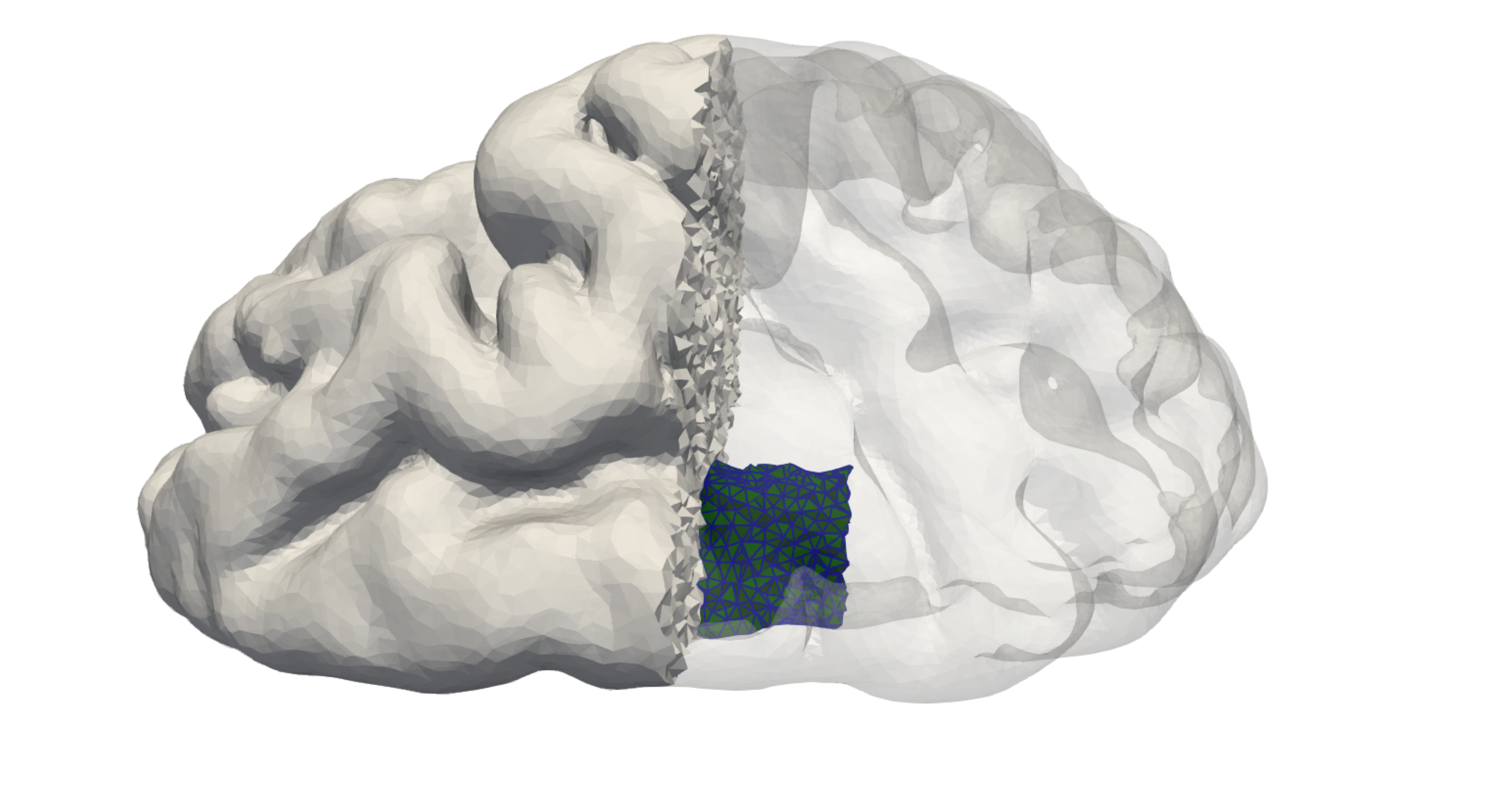}%
  \label{subfig:parent}}
  \qquad
  \subfloat[\centering Children agglomerates of the green element on the left.]{\includegraphics[width=.3\linewidth]{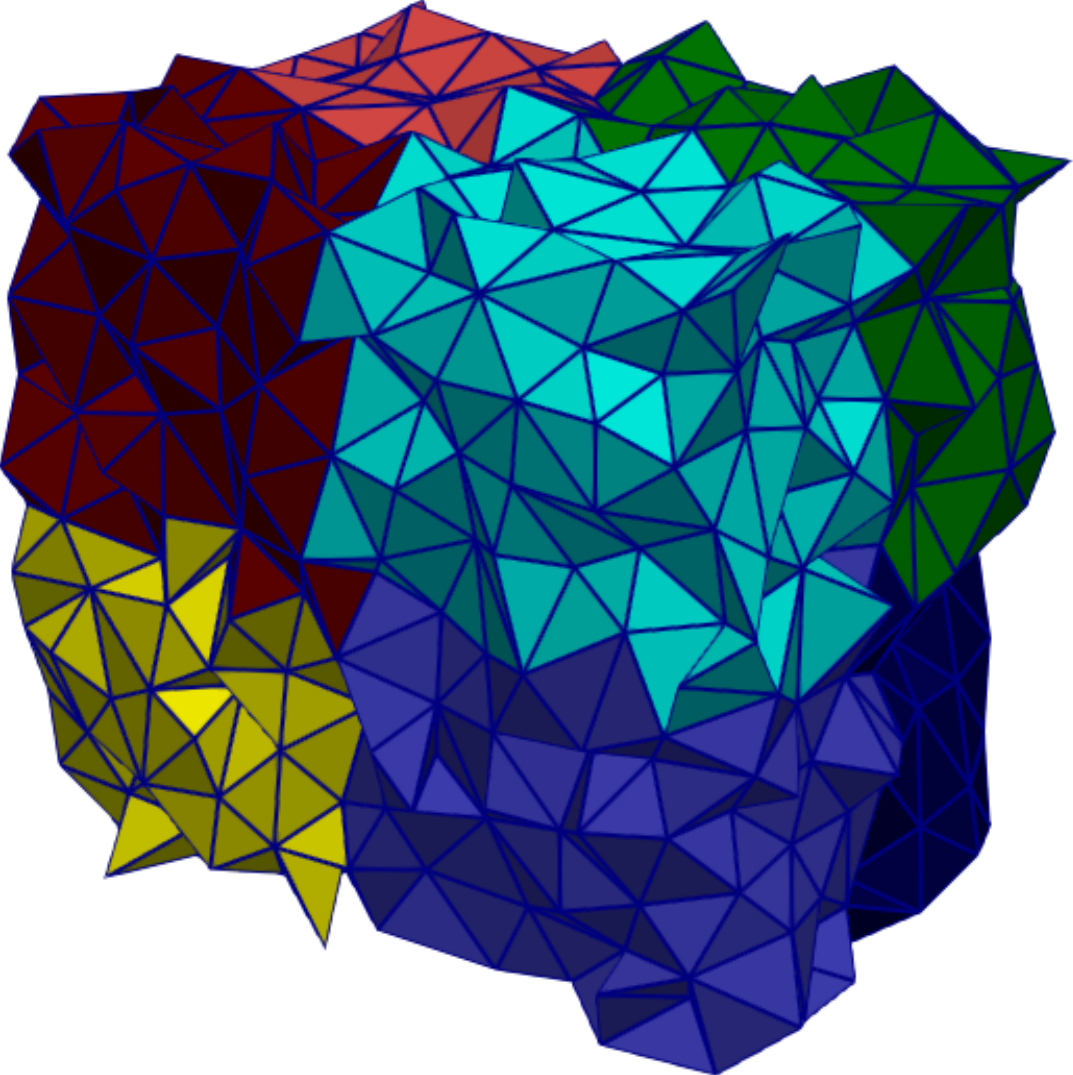}%
  \label{subfig:children}}
  \caption{Parent-child relationship between agglomerates at different levels. Each color maps to a different agglomerated element.}%
  \label{fig:brain_agglomerates}%
\end{figure}

To highlight the generality of our approach we apply it to a complex geometry representing a brain. We show in Figures~\ref{subfig:boundary_view} and~\ref{subfig:center_view} exploded
views of some agglomerates generated with the R-tree procedure for mesh $\Omega_7$. %
We extract one of the $7$ levels generated by the R-tree algorithm, obtaining
a computational polytopic mesh comprising $1,088$ agglomerates. In Figures~\ref{subfig:parent} and~\ref{subfig:children}, we display the hierarchical structure of some selected agglomerates. With a solid representation (in dark green), we display the
parent agglomerate, while in different colors we show its children. Notably, the shape of the resulting polytopes is quite close to the
associated bounding boxes, confirming the findings observed in the two-dimensional case.

\subsubsection{\emph{Test: liver mesh}}
We consider the case of meshes with different material properties attached to mesh elements. This case exemplifies a practical situation where the PDE model involves regions with different
coefficients associated to different parts of the geometry. In this scenario, it is desirable to have agglomerates that are not simultaneously living on both material regions. To automate this process, it is sufficient to build (in an automated and independent way) different trees \emph{for each} region of the mesh, which are then glued \emph{by construction} once the procedure is completed
for every region. The adaption of this approach to METIS is also possible, although the generation of a hierarchy is not automated, just as in the standard case. 

When regions differ substantially in size, shape, or number of local elements (as is the case for the liver mesh), the proposed methodology will still generate well-balanced trees, albeit of different depths. This is not an issue for the automation of the process, since the algorithm works independently on each region, resulting in an agglomerated coarsest grid which respects the material regions. However, the size of the coarsest grid will be  dictated by the shallowest tree.  
An alternative agglomeration strategy that addresses the automatic construction of region-preserving agglomerates was recently proposed in~\cite{antonietti2024polytopalmeshagglomerationgeometrical} where the authors present an agglomeration procedure based on deep learning techniques.

\begin{figure}[h]
  \centering
  \includegraphics[width=.6\linewidth]{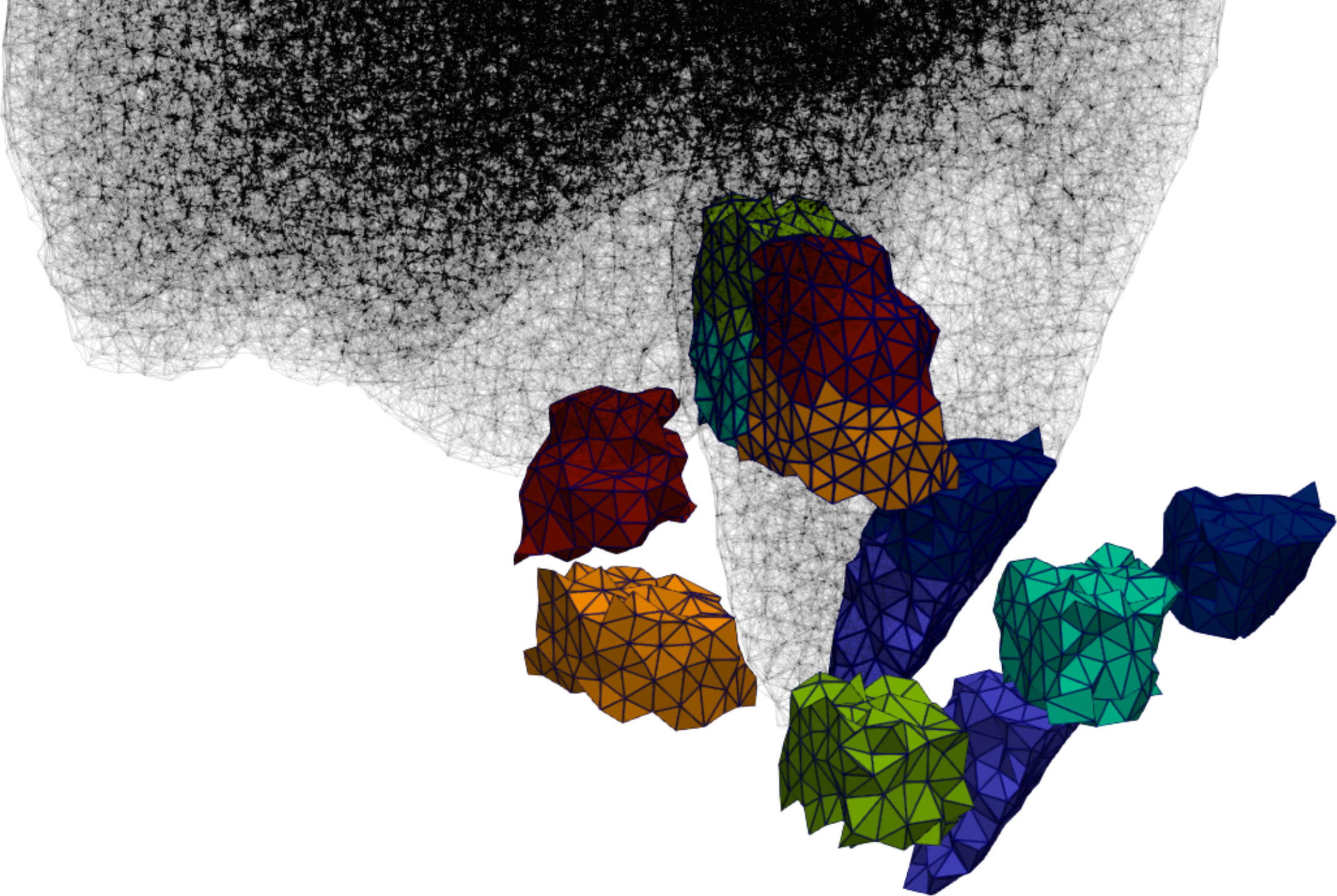}
  \caption{Exploded view of agglomerates (with different colors) generated for the bottom part of the liver mesh.}
  \label{fig:liver_test_view}
\end{figure}

A view of the agglomerates generated on the bottom part of the liver
can be visualized in Figure~\ref{fig:liver_test_view}. As it will be shown in detail in Section~\ref{sec:performance}, the construction of the R-tree data structure, as well as the definition of
the associated polytopic grid, does not constitute a bottleneck in terms of wall-clock time.

\subsection{Quality of resulting elements}
To assess the quality of the proposed methodology, we follow the same approach used in~\cite{antonietti2024agglomeration} and compute some of the quality metrics devised in \cite{ATTENE20211392}, which are reported hereafter for completeness:
\begin{itemize}
  \item \emph{Uniformity Factor}. Ratio between the diameter of an element $K$ and the mesh size $h$ defined as the maximum overall diameter of the polygons in the mesh: $$\text{UF}(K)=\frac{\text{diam}(K)}{h}.$$
  \item \emph{Circle Ratio}. Ratio between the radius of the inscribed circle and the radius of the circumscribed circle of an element $K$: $$\text{CR}(K)=\frac{\max_{\{B(r) \subset K \}}{r}}{\min_{\{K \subset B(r)\}}{r}}.$$
  \item \emph{Box Ratio}. Ratio between the measure of an element $K$ and the measure of its  bounding box: $$\text{BR}(K)=\frac{|K|}{|\text{MBR}(K)|}.$$
\end{itemize}
In addition to the previous metrics, we consider also the following:
\begin{itemize}
\item \emph{Overlap Factor}. Ratio between the global measure of the domain and the sum of the measures of all the bounding boxes in the polytopal mesh: $$\text{OF}(\Omega)=\frac{|\Omega|}{\sum_{i=1}^{N} |\text{MBR}(K_i)|},$$ where $N$ is the cardinality of the mesh $\Omega$.
\end{itemize}
All the metrics take values in $[0,1]$. In particular, higher average values of $\text{UF}$ imply that elements of the grid have comparable sizes, while $\text{CR}$ is a quality measure
related to the elongation or distortion of elements: the further away from $1$, the more elements tend to be elongated. To estimate this particular metric, we compute the radius of the inscribed circle of general polygons using the \textsc{C++} library \textsc{CGAL} \cite{CGALKernel} and approximate the radius of the inscribed circle with $\frac{\text{diam}(K)}{2}$. The metric \text{BR} takes into account how much an agglomerate $K$ and its bounding box are coinciding. Therefore, values close to $1$ indicate that the bounding box is tightly close to the agglomerate. 

The metric $\text{OF}(\Omega)$ is a global version of the metric \text{BR}: it estimates the sum of the overlaps between all bounding boxes relative to the measure of $\Omega$. In addition, it indicates how tightly the computational domain can be covered by boxes. For instance, we obtain the optimal $\text{BR}=1$ for the axis-aligned grid of Figure~\ref{fig:srtuctured_square_view}, cf. Table~\ref{tab:metrics1}.

In the table, we report the average values of $\text{UF}$, $\text{CR}$, and $\text{BR}$ for each geometry when we perform the first round of agglomeration of the grids displayed in Figure~\ref{fig:meshes} (top), namely the structure square, structure ball, and unstructured square grids. The last two columns of Table \ref{tab:metrics1} show the global overlap factor $\text{OF}(\Omega)$. The R-tree approach produces elements with shapes close to the associated Cartesian bounding box, resulting in values of $\text{OF}(\Omega)$ close to $1$ even in the case of underlying grids whose elements are \emph{not} axis-aligned. With reference to
the grid $\Omega_3$, which is fully unstructured, we are nevertheless able to obtain an almost optimal value for the overlap factor. In the case of $\Omega_2$, we have a slightly larger value of $1.2$  due to the curved geometry. In this particular case,  the contributions coming from agglomerates on the boundary exit from the domain due to the axis-aligned nature of the bounding boxes. 

Regarding the procedure employing \textsc{METIS}, we note that this is not meant to conserve any geometric information. This fact can be visually deduced by examining the Figures in Section \ref{subsec:validation_rtree}. The results reported in Table \ref{tab:metrics1} show that the R-tree grids are superior to the corresponding METIS grids in all four metrics.
In Figures \ref{fig:minmaxavg_level1} we show, for each grid and agglomeration strategy, the minimum, maximum, and average value of each metric.
The metrics obtained after performing a second round of agglomeration are reported in Table \ref{tab:metrics2} and Figure \ref{fig:minmaxavg_level2}. We observe that, on average, the metrics related to the R-tree approach are better than with \textsc{METIS}. However, looking at the second row of Tables~\ref{tab:metrics1} and~\ref{tab:metrics2} it can be seen that the minimum values of the CR metric can be slightly worse for the R-tree approach.
For what concerns the other metrics, \textsc{BR} is showing values always close to $1$ with the R-tree and high gaps compared to \textsc{METIS}, confirming the fact that shapes are preserved while increasing the levels. Furthermore, the overlap factor does not deteriorate using the R-tree, meaning that the global percentage of overlap of the bounding boxes is very low and that the resulting polygonal grid is very close to a global refinement of the coarser agglomerates.

  \bgroup
\def\arraystretch{1.5}%
  \begin{table}[H]
    \centering
    \resizebox{.8\textwidth}{!}{%
    \begin{tabular}{|c|c|cc|cc|cc|cc|c}
    \cline{1-10}
    Grid       & \# polygons & \multicolumn{2}{c|}{UF}              & \multicolumn{2}{c|}{CR}      & \multicolumn{2}{c|}{BR}     & \multicolumn{2}{c|}{OF}             &  \\ \cline{1-10}
               &  & \multicolumn{1}{c|}{\textbf{R-tree}} & \textbf{METIS}  & \multicolumn{1}{c|}{\textbf{R-tree}}   & \textbf{METIS} & \multicolumn{1}{c|}{\textbf{R-tree}} & \textbf{METIS} & \multicolumn{1}{c|}{\textbf{R-tree}} & \textbf{METIS} &  \\ \cline{1-10}
    $\Omega_1$ & 16 & \multicolumn{1}{c|}{1.0}      & 0.8441 & \multicolumn{1}{c|}{0.7071} & 0.5054 & \multicolumn{1}{c|}{1.0}     & 0.7713 & \multicolumn{1}{c|}{1.00} & 0.75 &  \\ \cline{1-10}
    $\Omega_2$ & 20 & \multicolumn{1}{c|}{0.7317} & 0.6193 & \multicolumn{1}{c|}{0.4432} & 0.3690 & \multicolumn{1}{c|}{0.8541}    & 0.5485 & \multicolumn{1}{c|}{0.8} &  0.50 &  \\ \cline{1-10}
    $\Omega_3$ & 91 & \multicolumn{1}{c|}{0.7622} & 0.7500 & \multicolumn{1}{c|}{0.7622} & 0.4070 & \multicolumn{1}{c|}{0.9235}    & 0.5965 & \multicolumn{1}{c|}{0.93} & 0.58 &  \\ \cline{1-10}
    \end{tabular}%
    }
    \caption[]{Average Values for the Uniformity Factor (UF), Circle Ratio (CR), and Box Ratio (BR), and global Overlap Factor (OF) for the \emph{first} level of agglomerated meshes with the two different agglomeration strategies. %
    }
    \label{tab:metrics1}
  \end{table}
  \egroup

\pgfplotsset{every error bar/.style={line width=1mm}}

\pgfplotstableread{
x         y             y-max             y-min
R-tree     1             0                 0
METIS     0.844103      0.155897          0.073219
}{\squareuf}
\pgfplotstableread{
x         y             y-max             y-min
R-tree     0.707107      0                 0
METIS     0.505435      0.07865           0.18177499999999996
}{\squarecr}
\pgfplotstableread{
x         y             y-max             y-min
R-tree     1             0                 0
METIS     0.77          0.1572            0.2462
}{\squarebr}

\pgfplotstableread{
x         y             y-max             y-min
R-tree     0.73172       0.26827           0.269369
METIS     0.619338      0.38              0.197  
}{\balluf}

\pgfplotstableread{
x         y             y-max             y-min
R-tree     0.443221     0.160              0.248
METIS     0.36906      0.111236           0.14
}{\ballcr}

\pgfplotstableread{
x         y             y-max             y-min
R-tree     0.854118      0.0924            0.2339
METIS     0.548554      0.1463            0.2423
}{\ballbr}

\pgfplotstableread{
x         y             y-max             y-min
R-tree     0.762249      0.23              0.29
METIS     0.750007      0.25              0.30  
}{\squareunstructureduf}

\pgfplotstableread{
x         y             y-max             y-min
R-tree     0.502087      0.1736            0.286913
METIS     0.407028      0.1482            0.155009 
}{\squareunstructuredcr}

\pgfplotstableread{
x         y             y-max             y-min
R-tree     0.923544      0.0364            0.0421
METIS     0.596511      0.1758            0.2323
}{\squareunstructuredbr}

\begin{figure}[h] 
  \centering 
\begin{tikzpicture}[scale=0.9] 
\begin{axis} [
  title={$\Omega_1$ (Structured square)},
  ylabel={Uniformity Factor},
  width  = 0.4*\textwidth,
  height = 5cm,
  ymin=0,
  ymax=1.1,  
  symbolic x coords={R-tree, METIS},
  minor ytick={0,0.25,0.5,0.75,1.},
  yminorgrids,
  xtick=data,
  ticklabel style = {font=\tiny},
  x tick label style={rotate=45,anchor=east},
  legend style={at={(0.4,0.5)},anchor=north west,cells={anchor=west},column
  sep=1ex}
]
\addplot+[black, very thick, forget plot,only marks] 
  plot[very thick, error bars/.cd, y dir=plus, y explicit]
  table[x=x,y=y,y error expr=\thisrow{y-max}] {\squareuf};
\addplot+[black, very thick, only marks,xticklabels=\empty] 
  plot[very thick, error bars/.cd, y dir=minus, y explicit]
  table[x=x,y=y,y error expr=\thisrow{y-min}] {\squareuf};
  \addplot[only marks,mark=square*,color=black] 
  table[x=x,y expr=\thisrow{y}+\thisrow{y-max}] {\squareuf};
  \addplot[only marks,mark=square*,color=black] 
  table[x=x,y expr=\thisrow{y}-\thisrow{y-min}] {\squareuf};
\end{axis} 
\end{tikzpicture}
\begin{tikzpicture}[scale=0.9] 
  \begin{axis} [
    title={$\Omega_2$ (Structured ball)},
    width  = 0.4*\textwidth,
    height = 5cm,
    ymin=0,
    ymax=1.1,  
    symbolic x coords={R-tree, METIS},
    minor ytick={0,0.25,0.5,0.75,1.},
    yminorgrids,
    xtick=data,
    ticklabel style = {font=\tiny},
    x tick label style={rotate=45,anchor=east},
    legend style={at={(0.4,0.5)},anchor=north west,cells={anchor=west},column
    sep=1ex}
  ]
  \addplot+[black, very thick, forget plot,only marks] 
    plot[very thick, error bars/.cd, y dir=plus, y explicit]
    table[x=x,y=y,y error expr=\thisrow{y-max}] {\balluf};
  \addplot+[black, very thick, only marks,xticklabels=\empty] 
    plot[very thick, error bars/.cd, y dir=minus, y explicit]
    table[x=x,y=y,y error expr=\thisrow{y-min}] {\balluf};
    \addplot[only marks,mark=square*,color=black] 
    table[x=x,y expr=\thisrow{y}+\thisrow{y-max}] {\balluf};
    \addplot[only marks,mark=square*,color=black] 
    table[x=x,y expr=\thisrow{y}-\thisrow{y-min}] {\balluf};
  \end{axis} 
\end{tikzpicture}
\begin{tikzpicture}[scale=0.9] 
  \begin{axis} [
    title={$\Omega_3$ (Unstructured square)},
    width  = 0.4*\textwidth,
    height = 5cm,
    ymin=0,
    ymax=1.1,  
    symbolic x coords={R-tree, METIS},
    minor ytick={0,0.25,0.5,0.75,1.},
yminorgrids,
xtick=data,
ticklabel style = {font=\tiny},
x tick label style={rotate=45,anchor=east},
legend style={at={(0.4,0.5)},anchor=north west,cells={anchor=west},column
sep=1ex}
]
\addplot+[black, very thick, forget plot,only marks] 
plot[very thick, error bars/.cd, y dir=plus, y explicit]
table[x=x,y=y,y error expr=\thisrow{y-max}] {\squareunstructureduf};
\addplot+[black, very thick, only marks,xticklabels=\empty] 
plot[very thick, error bars/.cd, y dir=minus, y explicit]
table[x=x,y=y,y error expr=\thisrow{y-min}] {\squareunstructureduf};
\addplot[only marks,mark=square*,color=black] 
table[x=x,y expr=\thisrow{y}+\thisrow{y-max}] {\squareunstructureduf};
\addplot[only marks,mark=square*,color=black] 
table[x=x,y expr=\thisrow{y}-\thisrow{y-min}] {\squareunstructureduf};
\end{axis} 
\end{tikzpicture}
\begin{tikzpicture}[scale=0.9] 
\begin{axis} [
  ylabel={Circle Ratio},
  width  = 0.4*\textwidth,
  height = 5cm,
  ymin=0,
  ymax=1.1,  
  symbolic x coords={R-tree, METIS},
  minor ytick={0,0.25,0.5,0.75,1.},
  yminorgrids,
  xtick=data,
  ticklabel style = {font=\tiny},
  x tick label style={rotate=45,anchor=east},
  legend style={at={(0.4,0.5)},anchor=north west,cells={anchor=west},column
  sep=1ex}
]
\addplot+[black, very thick, forget plot,only marks] 
  plot[very thick, error bars/.cd, y dir=plus, y explicit]
  table[x=x,y=y,y error expr=\thisrow{y-max}] {\squarecr};
\addplot+[black, very thick, only marks,xticklabels=\empty] 
  plot[very thick, error bars/.cd, y dir=minus, y explicit]
  table[x=x,y=y,y error expr=\thisrow{y-min}] {\squarecr};
  \addplot[only marks,mark=square*,color=black] 
  table[x=x,y expr=\thisrow{y}+\thisrow{y-max}] {\squarecr};
  \addplot[only marks,mark=square*,color=black] 
  table[x=x,y expr=\thisrow{y}-\thisrow{y-min}] {\squarecr};
\end{axis} 
\end{tikzpicture}
\begin{tikzpicture}[scale=0.9] 
  \begin{axis} [
    width  = 0.4*\textwidth,
    height = 5cm,
    ymin=0,
    ymax=1.1,  
    symbolic x coords={R-tree, METIS},
    minor ytick={0,0.25,0.5,0.75,1.},
    yminorgrids,
    xtick=data,
    ticklabel style = {font=\tiny},
    x tick label style={rotate=45,anchor=east},
    legend style={at={(0.4,0.5)},anchor=north west,cells={anchor=west},column
    sep=1ex}
  ]
  \addplot+[black, very thick, forget plot,only marks] 
    plot[very thick, error bars/.cd, y dir=plus, y explicit]
    table[x=x,y=y,y error expr=\thisrow{y-max}] {\ballcr};
  \addplot+[black, very thick, only marks,xticklabels=\empty] 
    plot[very thick, error bars/.cd, y dir=minus, y explicit]
    table[x=x,y=y,y error expr=\thisrow{y-min}] {\ballcr};
    \addplot[only marks,mark=square*,color=black] 
    table[x=x,y expr=\thisrow{y}+\thisrow{y-max}] {\ballcr};
    \addplot[only marks,mark=square*,color=black] 
    table[x=x,y expr=\thisrow{y}-\thisrow{y-min}] {\ballcr};
  \end{axis}
\end{tikzpicture}
\begin{tikzpicture}[scale=0.9] 
    \begin{axis} [
    width  = 0.4*\textwidth,
    height = 5cm,
    ymin=0,
    ymax=1.1,  
    symbolic x coords={R-tree, METIS},
    minor ytick={0,0.25,0.5,0.75,1.},
    yminorgrids,
    xtick=data,
    ticklabel style = {font=\tiny},
    x tick label style={rotate=45,anchor=east},
    legend style={at={(0.4,0.5)},anchor=north west,cells={anchor=west},column
    sep=1ex}
    ]
    \addplot+[black, very thick, forget plot,only marks] 
    plot[very thick, error bars/.cd, y dir=plus, y explicit]
    table[x=x,y=y,y error expr=\thisrow{y-max}] {\squareunstructuredcr};
    \addplot+[black, very thick, only marks,xticklabels=\empty] 
    plot[very thick, error bars/.cd, y dir=minus, y explicit]
    table[x=x,y=y,y error expr=\thisrow{y-min}] {\squareunstructuredcr};
    \addplot[only marks,mark=square*,color=black] 
    table[x=x,y expr=\thisrow{y}+\thisrow{y-max}] {\squareunstructuredcr};
    \addplot[only marks,mark=square*,color=black] 
    table[x=x,y expr=\thisrow{y}-\thisrow{y-min}] {\squareunstructuredcr};
  \end{axis} 
\end{tikzpicture}
\begin{tikzpicture}[scale=0.9] 
\begin{axis} [
  ylabel={Box Ratio},
  width  = 0.4*\textwidth,
  height = 5cm,
  ymin=0.0,
  ymax=1.1,  
  symbolic x coords={R-tree, METIS},
  minor ytick={0,0.25,0.5,0.75,1.},
  yminorgrids,
  xtick=data,
  ticklabel style = {font=\tiny},
  x tick label style={rotate=45,anchor=east},
  legend style={at={(0.4,0.5)},anchor=north west,cells={anchor=west},column
  sep=1ex}
]
\addplot+[black, very thick, forget plot,only marks] 
  plot[very thick, error bars/.cd, y dir=plus, y explicit]
  table[x=x,y=y,y error expr=\thisrow{y-max}] {\squarebr};
\addplot+[black, very thick, only marks,xticklabels=\empty] 
  plot[very thick, error bars/.cd, y dir=minus, y explicit]
  table[x=x,y=y,y error expr=\thisrow{y-min}] {\squarebr};
  \addplot[only marks,mark=square*,color=black] 
  table[x=x,y expr=\thisrow{y}+\thisrow{y-max}] {\squarebr};
  \addplot[only marks,mark=square*,color=black] 
  table[x=x,y expr=\thisrow{y}-\thisrow{y-min}] {\squarebr};
\end{axis} 
\end{tikzpicture}
\begin{tikzpicture}[scale=0.9] 
  \begin{axis} [
    width  = 0.4*\textwidth,
    height = 5cm,
    ymin=0.0,
    ymax=1.1,  
    symbolic x coords={R-tree, METIS},
    minor ytick={0,0.25,0.5,0.75,1.},
yminorgrids,
xtick=data,
ticklabel style = {font=\tiny},
x tick label style={rotate=45,anchor=east},
legend style={at={(0.4,0.5)},anchor=north west,cells={anchor=west},column
sep=1ex}
]
\addplot+[black, very thick, forget plot,only marks] 
plot[very thick, error bars/.cd, y dir=plus, y explicit]
table[x=x,y=y,y error expr=\thisrow{y-max}] {\ballbr};
\addplot+[black, very thick, only marks,xticklabels=\empty] 
plot[very thick, error bars/.cd, y dir=minus, y explicit]
table[x=x,y=y,y error expr=\thisrow{y-min}] {\ballbr};
\addplot[only marks,mark=square*,color=black] 
table[x=x,y expr=\thisrow{y}+\thisrow{y-max}] {\ballbr};
\addplot[only marks,mark=square*,color=black] 
table[x=x,y expr=\thisrow{y}-\thisrow{y-min}] {\ballbr};
\end{axis} 
\end{tikzpicture}
\begin{tikzpicture}[scale=0.9] 
  \begin{axis} [
    width  = 0.4*\textwidth,
    height = 5cm,
    ymin=0.0,
    ymax=1.1,    
    symbolic x coords={R-tree, METIS},
    minor ytick={0,0.25,0.5,0.75,1.},
yminorgrids,
xtick=data,
ticklabel style = {font=\tiny},
x tick label style={rotate=45,anchor=east},
legend style={at={(0.4,0.5)},anchor=north west,cells={anchor=west},column
sep=1ex}
]
\addplot+[black, very thick, forget plot,only marks] 
plot[very thick, error bars/.cd, y dir=plus, y explicit]
table[x=x,y=y,y error expr=\thisrow{y-max}] {\squareunstructuredbr};
\addplot+[black, very thick, only marks,xticklabels=\empty] 
plot[very thick, error bars/.cd, y dir=minus, y explicit]
table[x=x,y=y,y error expr=\thisrow{y-min}] {\squareunstructuredbr};
\addplot[only marks,mark=square*,color=black] 
table[x=x,y expr=\thisrow{y}+\thisrow{y-max}] {\squareunstructuredbr};
\addplot[only marks,mark=square*,color=black] 
table[x=x,y expr=\thisrow{y}-\thisrow{y-min}] {\squareunstructuredbr};
\end{axis} 
\end{tikzpicture}
\caption{Min, Max and Average chart of metrics UF, CR, and BR for grids $\Omega_1,\Omega_2,\Omega_3$ and different agglomeration strategies. In the plots on the left, metrics related to the R-tree collapse to a single dot since all mesh elements have the same value for that particular metric.}
\label{fig:minmaxavg_level1}
\end{figure}
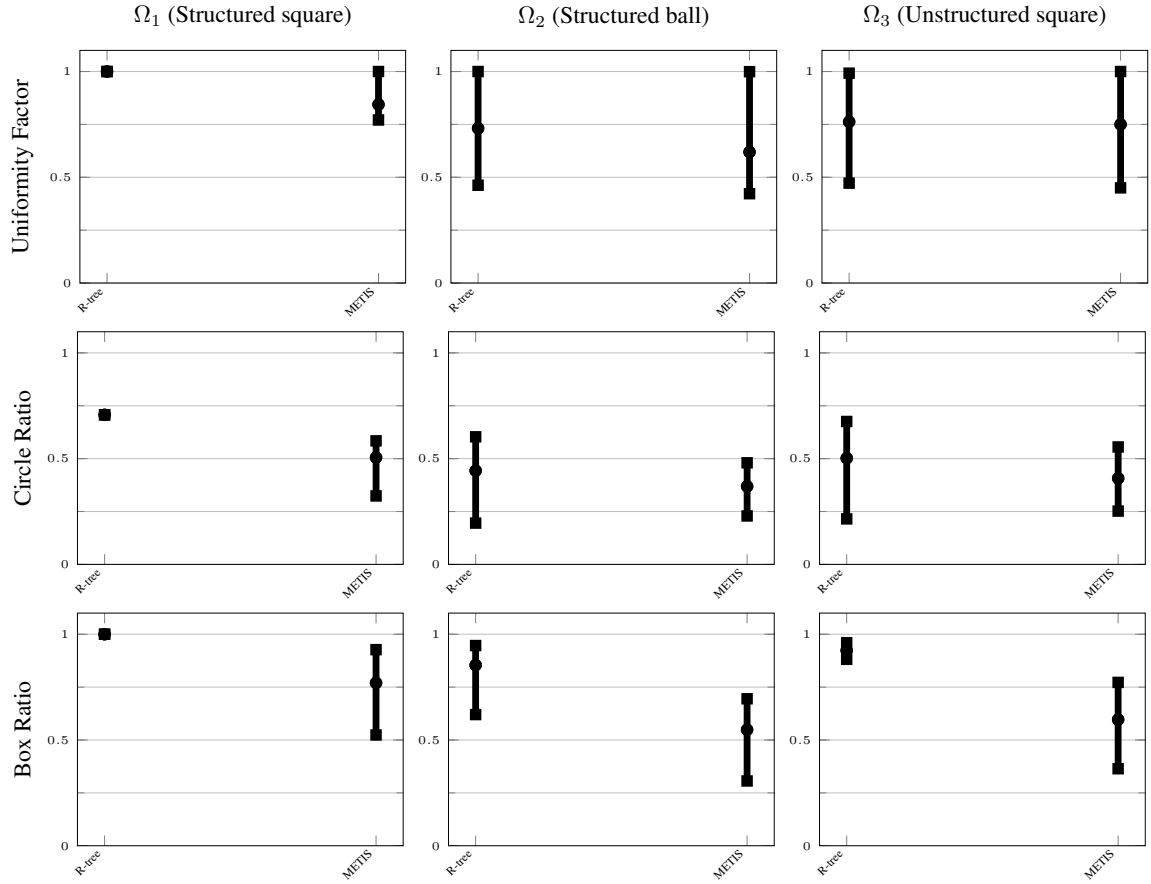

\bgroup
\def\arraystretch{1.5}%
  \begin{table}[H]
    \centering
    \resizebox{.8\textwidth}{!}{%
    \begin{tabular}{|c|c|cc|cc|cc|cc|c}
    \cline{1-10}
    Grid       &  \# polygons  & \multicolumn{2}{c|}{UF}              & \multicolumn{2}{c|}{CR}             & \multicolumn{2}{c|}{BR}                & \multicolumn{2}{c|}{OF}             &  \\ \cline{1-10}
               &    & \multicolumn{1}{c|}{\textbf{R-tree}} & \textbf{METIS}  & \multicolumn{1}{c|}{\textbf{R-tree}}   & \textbf{METIS}  & \multicolumn{1}{c|}{\textbf{R-tree}} & \textbf{METIS} & \multicolumn{1}{c|}{\textbf{R-tree}} & \textbf{METIS} &  \\ \cline{1-10}
    $\Omega_1$ & 64 & \multicolumn{1}{c|}{1.0}        & 0.6382 & \multicolumn{1}{c|}{0.7071} & 0.3368 & \multicolumn{1}{c|}{1.0}         & 0.5605   & \multicolumn{1}{c|}{1.00}& 0.5   &  \\ \cline{1-10}
    $\Omega_2$ & 80 & \multicolumn{1}{c|}{0.599658} & 0.5457 & \multicolumn{1}{c|}{0.4462} & 0.4002 & \multicolumn{1}{c|}{0.8356}    & 0.5890   & \multicolumn{1}{c|}{0.81}    & 0.56   &  \\ \cline{1-10}
    $\Omega_3$ & 364& \multicolumn{1}{c|}{0.698383} & 0.6111 & \multicolumn{1}{c|}{0.4799} & 0.3970 & \multicolumn{1}{c|}{0.8630}    & 0.5815   & \multicolumn{1}{c|}{0.87}    & 0.57   &  \\ \cline{1-10}
    \end{tabular}%
    }
    \caption[]{Average Values for the Uniformity Factor (UF), Circle Ratio (CR), Box Ratio (BR), and global Overlap Factor (OF) for the \emph{second} level of agglomerated meshes with the two different agglomeration strategies.%
    }
    \label{tab:metrics2}
    \end{table}
    \egroup

\pgfplotsset{every error bar/.style={line width=1mm}}

\pgfplotstableread{
x         y             y-max             y-min
R-tree     1             0                 0
METIS     0.726721      0.2732            0.1243
}{\refinedsquareuf}

\pgfplotstableread{
x         y             y-max             y-min
R-tree     0.707107      0                 0
METIS     0.357932      0.1909            0.1697
}{\refinedsquarecr}

\pgfplotstableread{
x         y             y-max             y-min
R-tree     1             0                 0
METIS     0.560528      0.4394            0.3605
}{\refinedsquarebr}

\pgfplotstableread{
x         y             y-max             y-min
R-tree     0.599658      0.4003            0.2204
METIS     0.547057      0.4529            0.2013
}{\refinedballuf}

\pgfplotstableread{
x         y             y-max             y-min
R-tree     0.446243      0.1842            0.3060
METIS     0.401347      0.1595            0.2524
}{\refinedballcr}

\pgfplotstableread{
x         y             y-max             y-min
R-tree     0.835658      0.1124            0.29682
METIS     0.589042      0.1605            0.2815
}{\refinedballbr}

\pgfplotstableread{
x         y             y-max             y-min
R-tree     0.698383      0.3016            0.3105
METIS     0.708752      0.2912            0.3303
}{\refinedsquareunstructureduf}

\pgfplotstableread{
x         y             y-max             y-min
R-tree     0.479989      0.1815            0.2903
METIS     0.394476      0.2013            0.1737
}{\refinedsquareunstructuredcr}

\pgfplotstableread{
x         y             y-max             y-min
R-tree     0.863073      0.0845            0.0846
METIS     0.581566      0.2517            0.2636
}{\refinedsquareunstructuredbr}

\begin{figure}[h] 
  \centering 
\begin{tikzpicture}[scale=0.9] 
\begin{axis} [
  title={$\Omega_1$ (Square)},
  ylabel={Uniformity Factor},
  width  = 0.4*\textwidth,
  height = 5cm,
  ymin=0,
  ymax=1.1,  
  symbolic x coords={R-tree, METIS},
  minor ytick={0,0.25,0.5,0.75,1.},
  yminorgrids,
  xtick=data,
  ticklabel style = {font=\tiny},
  x tick label style={rotate=45,anchor=east},
  legend style={at={(0.4,0.5)},anchor=north west,cells={anchor=west},column
  sep=1ex}
]
\addplot+[black, very thick, forget plot,only marks] 
  plot[very thick, error bars/.cd, y dir=plus, y explicit]
  table[x=x,y=y,y error expr=\thisrow{y-max}] {\refinedsquareuf};
\addplot+[black, very thick, only marks,xticklabels=\empty] 
  plot[very thick, error bars/.cd, y dir=minus, y explicit]
  table[x=x,y=y,y error expr=\thisrow{y-min}] {\refinedsquareuf};
  \addplot[only marks,mark=square*,color=black] 
  table[x=x,y expr=\thisrow{y}+\thisrow{y-max}] {\refinedsquareuf};
  \addplot[only marks,mark=square*,color=black] 
  table[x=x,y expr=\thisrow{y}-\thisrow{y-min}] {\refinedsquareuf};
\end{axis} 
\end{tikzpicture}
\begin{tikzpicture}[scale=0.9] 
  \begin{axis} [
    title={$\Omega_2$ (Ball)},
    width  = 0.4*\textwidth,
    height = 5cm,
    ymin=0,
    ymax=1.1,  
    symbolic x coords={R-tree, METIS},
    minor ytick={0,0.25,0.5,0.75,1.},
    yminorgrids,
    xtick=data,
    ticklabel style = {font=\tiny},
    x tick label style={rotate=45,anchor=east},
    legend style={at={(0.4,0.5)},anchor=north west,cells={anchor=west},column
    sep=1ex}
  ]
  \addplot+[black, very thick, forget plot,only marks] 
    plot[very thick, error bars/.cd, y dir=plus, y explicit]
    table[x=x,y=y,y error expr=\thisrow{y-max}] {\refinedballuf};
  \addplot+[black, very thick, only marks,xticklabels=\empty] 
    plot[very thick, error bars/.cd, y dir=minus, y explicit]
    table[x=x,y=y,y error expr=\thisrow{y-min}] {\refinedballuf};
    \addplot[only marks,mark=square*,color=black] 
    table[x=x,y expr=\thisrow{y}+\thisrow{y-max}] {\refinedballuf};
    \addplot[only marks,mark=square*,color=black] 
    table[x=x,y expr=\thisrow{y}-\thisrow{y-min}] {\refinedballuf};
  \end{axis} 
\end{tikzpicture}
\begin{tikzpicture}[scale=0.9] 
  \begin{axis} [
    title={$\Omega_3$ (Unstructured square)},
    width  = 0.4*\textwidth,
    height = 5cm,
    ymin=0,
    ymax=1.1,  
    symbolic x coords={R-tree, METIS},
    minor ytick={0,0.25,0.5,0.75,1.},
yminorgrids,
xtick=data,
ticklabel style = {font=\tiny},
x tick label style={rotate=45,anchor=east},
legend style={at={(0.4,0.5)},anchor=north west,cells={anchor=west},column
sep=1ex}
]
\addplot+[black, very thick, forget plot,only marks] 
plot[very thick, error bars/.cd, y dir=plus, y explicit]
table[x=x,y=y,y error expr=\thisrow{y-max}] {\refinedsquareunstructureduf};
\addplot+[black, very thick, only marks,xticklabels=\empty] 
plot[very thick, error bars/.cd, y dir=minus, y explicit]
table[x=x,y=y,y error expr=\thisrow{y-min}] {\refinedsquareunstructureduf};
\addplot[only marks,mark=square*,color=black] 
table[x=x,y expr=\thisrow{y}+\thisrow{y-max}] {\refinedsquareunstructureduf};
\addplot[only marks,mark=square*,color=black] 
table[x=x,y expr=\thisrow{y}-\thisrow{y-min}] {\refinedsquareunstructureduf};
\end{axis} 
\end{tikzpicture}
\begin{tikzpicture}[scale=0.9] 
\begin{axis} [
  ylabel={Circle Ratio},
  width  = 0.4*\textwidth,
  height = 5cm,
  ymin=0,
  ymax=1.1,  
  symbolic x coords={R-tree, METIS},
  minor ytick={0,0.25,0.5,0.75,1.},
  yminorgrids,
  xtick=data,
  ticklabel style = {font=\tiny},
  x tick label style={rotate=45,anchor=east},
  legend style={at={(0.4,0.5)},anchor=north west,cells={anchor=west},column
  sep=1ex}
]
\addplot+[black, very thick, forget plot,only marks] 
  plot[very thick, error bars/.cd, y dir=plus, y explicit]
  table[x=x,y=y,y error expr=\thisrow{y-max}] {\refinedsquarecr};
\addplot+[black, very thick, only marks,xticklabels=\empty] 
  plot[very thick, error bars/.cd, y dir=minus, y explicit]
  table[x=x,y=y,y error expr=\thisrow{y-min}] {\refinedsquarecr};
  \addplot[only marks,mark=square*,color=black] 
  table[x=x,y expr=\thisrow{y}+\thisrow{y-max}] {\refinedsquarecr};
  \addplot[only marks,mark=square*,color=black] 
  table[x=x,y expr=\thisrow{y}-\thisrow{y-min}] {\refinedsquarecr};
\end{axis} 
\end{tikzpicture}
\begin{tikzpicture}[scale=0.9] 
  \begin{axis} [
    width  = 0.4*\textwidth,
    height = 5cm,
    ymin=0,
    ymax=1.1,  
    symbolic x coords={R-tree, METIS},
    minor ytick={0,0.25,0.5,0.75,1.},
    yminorgrids,
    xtick=data,
    ticklabel style = {font=\tiny},
    x tick label style={rotate=45,anchor=east},
    legend style={at={(0.4,0.5)},anchor=north west,cells={anchor=west},column
    sep=1ex}
  ]
  \addplot+[black, very thick, forget plot,only marks] 
    plot[very thick, error bars/.cd, y dir=plus, y explicit]
    table[x=x,y=y,y error expr=\thisrow{y-max}] {\refinedballcr};
  \addplot+[black, very thick, only marks,xticklabels=\empty] 
    plot[very thick, error bars/.cd, y dir=minus, y explicit]
    table[x=x,y=y,y error expr=\thisrow{y-min}] {\refinedballcr};
    \addplot[only marks,mark=square*,color=black] 
    table[x=x,y expr=\thisrow{y}+\thisrow{y-max}] {\refinedballcr};
    \addplot[only marks,mark=square*,color=black] 
    table[x=x,y expr=\thisrow{y}-\thisrow{y-min}] {\refinedballcr};
  \end{axis}
\end{tikzpicture}
\begin{tikzpicture}[scale=0.9] 
    \begin{axis} [
    width  = 0.4*\textwidth,
    height = 5cm,
    ymin=0,
    ymax=1.1,  
    symbolic x coords={R-tree, METIS},
    minor ytick={0,0.25,0.5,0.75,1.},
    yminorgrids,
    xtick=data,
    ticklabel style = {font=\tiny},
    x tick label style={rotate=45,anchor=east},
    legend style={at={(0.4,0.5)},anchor=north west,cells={anchor=west},column
    sep=1ex}
    ]
    \addplot+[black, very thick, forget plot,only marks] 
    plot[very thick, error bars/.cd, y dir=plus, y explicit]
    table[x=x,y=y,y error expr=\thisrow{y-max}] {\refinedsquareunstructuredcr};
    \addplot+[black, very thick, only marks,xticklabels=\empty] 
    plot[very thick, error bars/.cd, y dir=minus, y explicit]
    table[x=x,y=y,y error expr=\thisrow{y-min}] {\refinedsquareunstructuredcr};
    \addplot[only marks,mark=square*,color=black] 
    table[x=x,y expr=\thisrow{y}+\thisrow{y-max}] {\refinedsquareunstructuredcr};
    \addplot[only marks,mark=square*,color=black] 
    table[x=x,y expr=\thisrow{y}-\thisrow{y-min}] {\refinedsquareunstructuredcr};
  \end{axis} 
\end{tikzpicture}
\begin{tikzpicture}[scale=0.9] 
\begin{axis} [
  ylabel={Box Ratio},
  width  = 0.4*\textwidth,
  height = 5cm,
  ymin=0.0,
  ymax=1.1,  
  symbolic x coords={R-tree, METIS},
  minor ytick={0,0.25,0.5,0.75,1.},
  yminorgrids,
  xtick=data,
  ticklabel style = {font=\tiny},
  x tick label style={rotate=45,anchor=east},
  legend style={at={(0.4,0.5)},anchor=north west,cells={anchor=west},column
  sep=1ex}
]
\addplot+[black, very thick, forget plot,only marks] 
  plot[very thick, error bars/.cd, y dir=plus, y explicit]
  table[x=x,y=y,y error expr=\thisrow{y-max}] {\squarebr};
\addplot+[black, very thick, only marks,xticklabels=\empty] 
  plot[very thick, error bars/.cd, y dir=minus, y explicit]
  table[x=x,y=y,y error expr=\thisrow{y-min}] {\squarebr};
  \addplot[only marks,mark=square*,color=black] 
  table[x=x,y expr=\thisrow{y}+\thisrow{y-max}] {\squarebr};
  \addplot[only marks,mark=square*,color=black] 
  table[x=x,y expr=\thisrow{y}-\thisrow{y-min}] {\squarebr};
\end{axis} 
\end{tikzpicture}
\begin{tikzpicture}[scale=0.9] 
  \begin{axis} [
    width  = 0.4*\textwidth,
    height = 5cm,
    ymin=0.0,
    ymax=1.1,  
    symbolic x coords={R-tree, METIS},
    minor ytick={0,0.25,0.5,0.75,1.},
yminorgrids,
xtick=data,
ticklabel style = {font=\tiny},
x tick label style={rotate=45,anchor=east},
legend style={at={(0.4,0.5)},anchor=north west,cells={anchor=west},column
sep=1ex}
]
\addplot+[black, very thick, forget plot,only marks] 
plot[very thick, error bars/.cd, y dir=plus, y explicit]
table[x=x,y=y,y error expr=\thisrow{y-max}] {\ballbr};
\addplot+[black, very thick, only marks,xticklabels=\empty] 
plot[very thick, error bars/.cd, y dir=minus, y explicit]
table[x=x,y=y,y error expr=\thisrow{y-min}] {\ballbr};
\addplot[only marks,mark=square*,color=black] 
table[x=x,y expr=\thisrow{y}+\thisrow{y-max}] {\ballbr};
\addplot[only marks,mark=square*,color=black] 
table[x=x,y expr=\thisrow{y}-\thisrow{y-min}] {\ballbr};
\end{axis} 
\end{tikzpicture}
\begin{tikzpicture}[scale=0.9] 
  \begin{axis} [
    width  = 0.4*\textwidth,
    height = 5cm,
    ymin=0.0,
    ymax=1.1,  
    symbolic x coords={R-tree, METIS},
    minor ytick={0,0.25,0.5,0.75,1.},
yminorgrids,
xtick=data,
ticklabel style = {font=\tiny},
x tick label style={rotate=45,anchor=east},
legend style={at={(0.4,0.5)},anchor=north west,cells={anchor=west},column
sep=1ex}
]
\addplot+[black, very thick, forget plot,only marks] 
plot[very thick, error bars/.cd, y dir=plus, y explicit]
table[x=x,y=y,y error expr=\thisrow{y-max}] {\squareunstructuredbr};
\addplot+[black, very thick, only marks,xticklabels=\empty] 
plot[very thick, error bars/.cd, y dir=minus, y explicit]
table[x=x,y=y,y error expr=\thisrow{y-min}] {\squareunstructuredbr};
\addplot[only marks,mark=square*,color=black] 
table[x=x,y expr=\thisrow{y}+\thisrow{y-max}] {\squareunstructuredbr};
\addplot[only marks,mark=square*,color=black] 
table[x=x,y expr=\thisrow{y}-\thisrow{y-min}] {\squareunstructuredbr};
\end{axis} 
\end{tikzpicture}
\caption{Min, Max and Average chart of metrics UF, CR, and BR for grids $\Omega_1,\Omega_2,\Omega_3$ and different agglomeration strategies, after creation of one more level. In the plots on the left, metrics related to the R-tree collapse to a single dot since all mesh elements have the same value for that particular metric.}
\label{fig:minmaxavg_level2}
\end{figure}
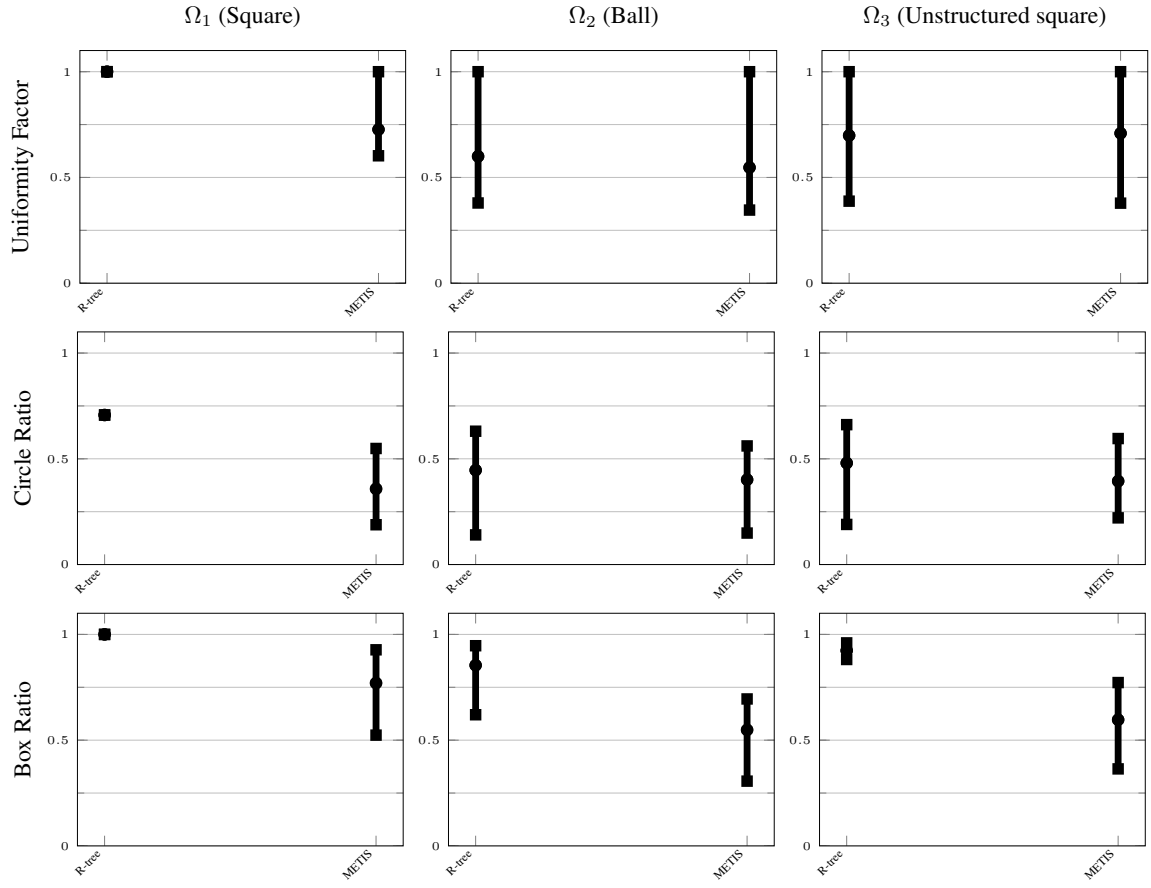

\subsection{Performance validation}
\label{sec:performance}
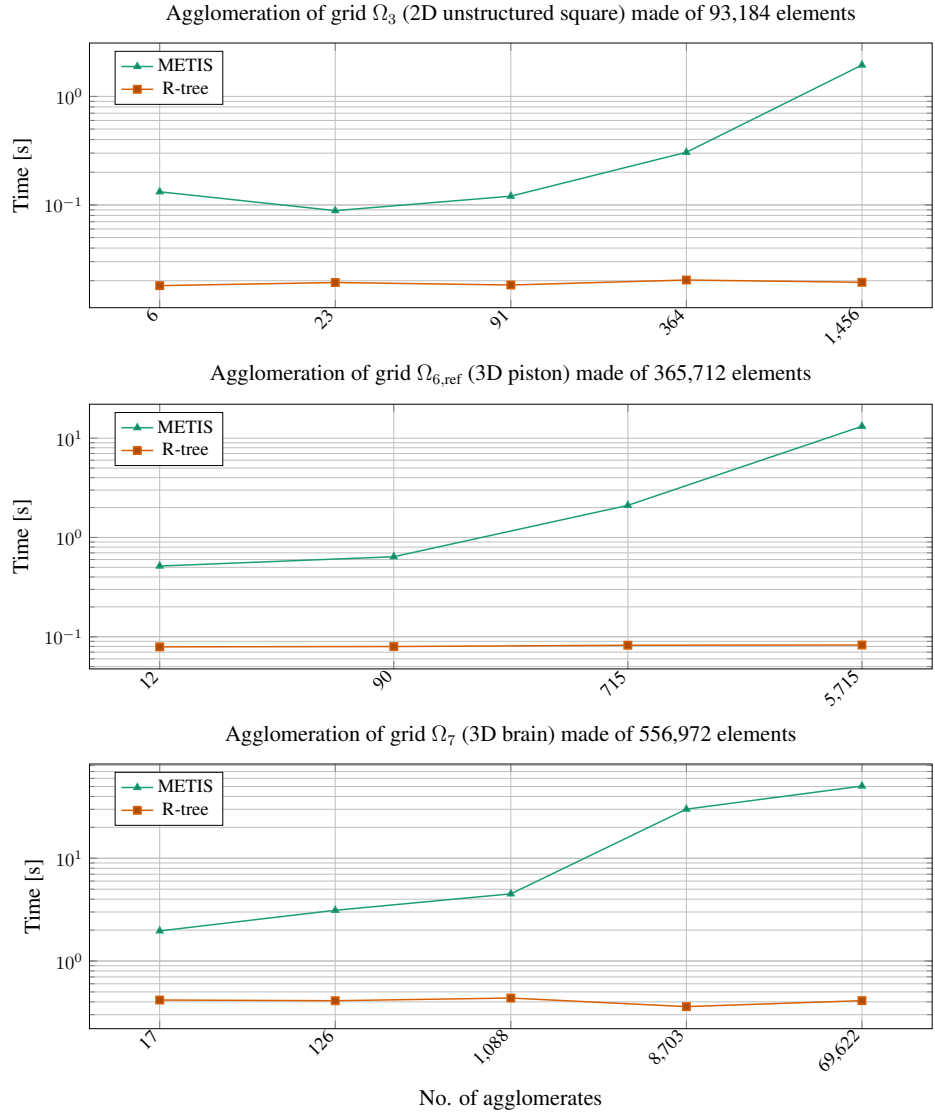
\begin{figure}[h]
  \centering
\begin{tikzpicture}[trim axis left,scale=0.7]
  \begin{semilogyaxis}[
      scale only axis, %
      height=5cm,
      width=\textwidth, %
      grid=both,
      title= \large{Agglomeration of grid $\Omega_3$ (2D unstructured square) made of 93,184 elements},
      ylabel= \large{Time [s]},
      legend pos=north west,
      x tick label style={rotate=45, anchor=east, align=center},
      xtick={0,.2,.4,.6,.8},
      xticklabels={6,23,91,364,1\text{,}456},
      ]
      \addplot+[solid, every mark/.append style={solid}, mark=triangle*]
      coordinates {
        (0,0.132089)(.2,0.0885683)(.4,0.120304 )(.6,0.305496)(.8,1.94425 )
        };
        \addlegendentry{METIS};
      \addplot+[solid, every mark/.append style={solid}, mark=square*]
      coordinates {
        (0,0.0180248)(.2,0.0192636)(.4,0.0182785)(.6,0.0203172)(.8,0.0193534)
        };
      \addlegendentry{R-tree};
 
  \end{semilogyaxis}
\end{tikzpicture}

\begin{tikzpicture}[trim axis left,scale=0.7]

  \begin{semilogyaxis}[
      scale only axis, %
      height=5cm,
      width=\textwidth, %
      grid=both,
      title= \large{Agglomeration of grid $\Omega_{6,\text{ref}}$ (3D piston) made of 365,712 elements},
      ylabel=\large{Time [s]},
      legend pos=north west,
      x tick label style={rotate=45, anchor=east, align=center},
      xtick={0,.2,.4,.6,.8},
      xticklabels={12,90,715,5\text{,}715},
      ]
      \addplot+[solid, every mark/.append style={solid}, mark=triangle*]
      coordinates {
        (0,0.515566)(.2,0.639702)(.4,2.10608)(.6,13.1708)
        };
        \addlegendentry{METIS};
        \addplot+[solid, every mark/.append style={solid}, mark=square*]
        coordinates {
        (0,0.0791843)(.2,0.0797387)(.4,0.0821672)(.6,0.0826353)
        };
      \addlegendentry{R-tree};
 
  \end{semilogyaxis}
\end{tikzpicture}

\begin{tikzpicture}[trim axis left,scale=0.7]
  \begin{semilogyaxis}[
      scale only axis, %
      height=5cm,
      width=\textwidth, %
      grid=both,
      title= \large{Agglomeration of grid $\Omega_7$ (3D brain) made of 556,972 elements},
      xlabel= \large{No. of  agglomerates},
      ylabel= \large{Time [s]},
      legend pos=north west,
      x tick label style={rotate=45, anchor=east, align=center},
      xtick={0,.2,.4,.6,.8},
      xticklabels={17,126,1\text{,}088,8\text{,}703,69\text{,}622},
      ]
      \addplot+[solid, every mark/.append style={solid}, mark=triangle*]
      coordinates {
        (0.,1.9634)(.2,3.1164)(.4,4.50613)(.6,30.0712)(.8,50.4756)
        };
        \addlegendentry{METIS};
        \addplot+[solid, every mark/.append style={solid}, mark=square*]
        coordinates {
          (0.,0.416847)(.2,0.410571)(.4,0.435945)(.6,0.359528)(.8,0.411066)
        };
      \addlegendentry{R-tree}; 
  \end{semilogyaxis}
\end{tikzpicture}
\caption{Wall-clock time (in seconds) needed to build polytopal grids with R-tree and METIS.}
\label{fig:benchmark}
\end{figure}

To further assess our method, a breakdown of the computing times required by the agglomeration process is tracked both in 2D
and in 3D test cases. The grid used for the 2D case is $\Omega_3$ (the unstructured square $[0,1]^2$), while for 3D meshes we consider
a globally refined version of the piston grid in Figure \ref{fig:piston_view}, referred to as $\Omega_{6,\text{ref}}$ and the brain mesh $\Omega_7$. %
The elapsed time we are interested in measuring is the time to build the target polytopal mesh. This means that different
components must be measured depending on the agglomeration strategy. In the case of METIS, we measure the wall-clock time (in seconds) associated with the call to the
METIS function \texttt{METIS\_PartGraphKway()}. For the R-tree-based strategy, we time cumulatively the following phases:
\begin{itemize}
  \item Build the R-tree;
  \item Visit the hierarchy and store data structures needed to generate agglomerates;
  \item Flag elements of the underlying triangulation $\Omega$.
\end{itemize}
We perform 10 runs for each agglomeration strategy and average the recorded wall-clock times. All the experiments have been performed on a 2.60GHz Intel Xeon processor. We observe from Figure~\ref{fig:benchmark} that the R-tree based approach keeps
constant timings which are independent of the extraction level. On the other hand, a graph partitioner
requires increasing computational times which in all cases are orders of magnitude larger than those required by the R-tree-based strategy. Moreover, the parallel extension of our approach in the case of distributed grids is conceptually straightforward. In an MPI-distributed framework, each process stores only a local part of the computational domain. And, within each locally owned partition, agglomeration can be performed locally. Since the original fine grid is already distributed, it must
be noted that the ghost polytopal elements do not live on the layer of standard ghost cells as is usual with classical FEMs implementations, but are owned by a different partition. Our algorithmic realization
carefully avoids calls to potentially expensive collective MPI routines and performs the exchange of such information in a setup phase.

%% file: poly_dg.tex
\section{Polygonal discontinuous Galerkin}
\label{sec:polydg}
\subsection{Notation and model problem}

For a Lipschitz domain $\omega \subset {\mathbb R}^d$, $d=1,2,3$,  we denote by $H^s(\omega)$ the Hilbertian Sobolev space of  index $s\ge 0$ of real--valued functions defined on
$\omega$, endowed with the seminorm $|\cdot |_{H^s(\omega)}$ and norm $\|\cdot\|_{H^s(\omega)}$. Furthermore, we let~$L^p(\omega)$, $p\in[1,\infty]$,  be the standard Lebesgue space on $\omega$, equipped with the norm~$\|\cdot\|_{L^p(\omega)}$. In the case $p=2$, we shall simply write $\|\cdot\|_{\omega}$  to denote the $L^2$-norm over $\omega$ and simplify this further to $\|\cdot\|$ when $\omega=\Omega$, the physical domain. We denote with $|\omega|$ the $d$-dimensional Hausdorff measure of $\omega$.

Let $\Omega$ be a bounded, simply connected, and open polygonal/polyhedral domain in $\mathbb{R}^d$, $d=2,3$. We indicate with $\Gamma_{\ddd} \coloneqq \partial\Omega$ the boundary $\Omega$. We consider the linear elliptic problem: find $u \in H^1(\Omega)$, such that

{
\begin{equation}\label{eq:problem}
  \begin{cases}
    -\Delta u &= f \qquad \text{in   } \Omega, \\
    u &= g \qquad \text{on }  \Gamma_{\ddd}
  \end{cases}
  \end{equation}}
with data $f \in L^2(\Omega)$, $g \in H^{1/2}(\Gamma_{\ddd})$. Setting $H^1_\ddd:= \{v\in H^1(\Omega): v = 0  \text{ on } \Gamma_{\ddd}  \}$, %
the weak formulation of~\eqref{eq:problem} reads: find $u\in H^1(\Omega)$, with $ u = g$ on $\Gamma_{\ddd}$ such that
{\begin{equation}\label{eq:prob}
\int_\Omega  \nabla u\cdot \nabla v\ud \uu{x}  =  \int_\Omega f  v \ud \uu{x}
\end{equation}}
for all $v\in H^1_{\ddd}(\Omega)$.  The well-posedness of the weak problem~\eqref{eq:prob} is guaranteed by the Lax-Milgram Lemma.

\subsection{Finite element spaces and trace operators}
We consider meshes consisting of general polygonal (for $d=2$) or polyhedral (for $d=3$) mutually disjoint open elements $K \in \Omega$, henceforth termed collectively as \emph{polytopic},
with $\cup_{K \in \Omega}\bar{K} =\bar{\Omega}$.

Given $h_{K}:=\diam(K)$, the diameter of $K \in \Omega$, we define the mesh-function $\mbf{h}:\cup_{K \in \Omega} K\to\mathbb{R}_+$ by $\mbf{h}|_{ K}=h_{ K}$, $K \in \Omega$.
Further, we let $\Gamma:=\cup_{K \in \Omega}\partial K$ denote the mesh skeleton and set $\Gamma_{\dint}:=\Gamma\backslash\partial\Omega$. The mesh skeleton $\Gamma$ is decomposed into $(d-1)$--dimensional simplices $F$ denoting the mesh \emph{faces}, shared by at most two elements. These are distinct from elemental \emph{interfaces}, which are defined as the simply connected components of the intersection between the boundary of an element and either a neighbouring element or $\partial \Omega$.
As such,  an interface between two elements may consist of more than one face, separated by hanging nodes/edges shared by those two elements only.

Over $\Omega$ we introduce the discontinuous \emph{finite element space} defined by
\[
V_h :=\{u\in L^2(\Omega)
:u|_{ K}\in\mathcal{Q}_{p}( K),K \in \Omega\},
\]
for some $p\in\mathbb{N}$
with,
$\mathcal{Q}_{p}(K)$ denoting the space of tensor-product polynomials of degree $p$ on $K$. In our implementation we exploit once again the bounding boxes to construct $V_h$: we use readily available tensor-product basis for polynomials defined on the bounding boxes and then consider their restriction over the physical element $K$; see Section~\ref{sec:convergence} for more details.

Let $K_i$ and $K_j$ be two adjacent elements of $\Omega$ sharing a face $F\subset\partial K_i \cap \partial K_j \subset \Gamma_{\dint}$. The outward unit normal
vectors on $F$ of $\partial K_i$ and $\partial K_j$ are denoted by $\boldsymbol{n}_{K_i}$ and $\boldsymbol{n}_{K_j}$, respectively. For $v$ and $\mbf{q}$ element-wise continuous scalar- and vector-valued functions, respectively,  we define the \emph{average} across $F$ by $
\mean{v}|_F:=\frac{1}{2}(v|_{F\cap K_i}+v|_{F\cap K_j})$, $\mean{\mbf{q}}|_F:=\frac{1}{2}(\mbf{q}|_{F\cap K_i}+\mbf{q}|_{F\cap K_j})$,
respectively, and the \emph{jump} across $F$ by
$
\jump{v}|_F :=v|_{F\cap K_i}\boldsymbol{n}_{K_i} + v|_{F\cap K_j} \boldsymbol{n}_{K_j}$, $\jump{\mbf{q}}|_F := \mbf{q}|_{F\cap K_i} \cdot \boldsymbol{n}_{K_i} + \mbf{q}|_{F\cap K_j} \cdot \boldsymbol{n}_{K_j}$,
using the convention $i>j$ in the element numbering to determine the sign.
On a boundary face $F\subset  \Gamma_{\ddd}$, with $F \subset \partial K_i$, $K_i\in \Omega$,
we set
$
\mean{v}|_F:=v_i, $  $ \mean{\mbf{q}}|_F:=\mbf{q}_i \cdot \boldsymbol{n},  $
$\jump{v}|_F :=v_i ,$ and $ \jump{\mbf{q}}|_F := \mbf{q}_i \boldsymbol{n}$, respectively, with $\boldsymbol{n}$
the outward unit normal to the boundary $\partial \Omega$.

For $v\in V_h $  we denote by $\nabla_h v$ the element-wise gradient; namely, $(\nabla_h v)|_{K}:=\nabla (v|_K)$ for all $K \in \Omega$. Then,
the symmetric interior penalty discontinuous Galerkin method reads: find $u_h \in V_h$ such that

\begin{equation}
  \label{eqn:weak_form}
  B(u_h,v_h) = l(v_h), \qquad \forall v_h \in V_h,
\end{equation}
with 
{\begin{equation}
  B(u_h,v_h) = \int_\Omega  \nabla_h u_h \cdot \nabla_h v_h \ud \uu{x}
  - \int_{\Gamma} \Bigl( \mean{  \nabla u_h} \cdot \jump{v_h} + \mean{ \nabla v_h} \cdot \jump{u_h} \Bigr)   \ud s
  + \int_{\Gamma} \sigma \jump{u_h} \cdot \jump{v_h} \ud s,
\end{equation}}
and
{\begin{equation}
  l(v_h) = \int_\Omega f  v \ud \uu{x} + \int_{\Gamma_{\ddd}}  g(\sigma v_h - \nabla v_h \cdot \mbf{n}) \ud s,
\end{equation}}
where $\sigma: \Gamma \rightarrow \mathbb{R}$ is the so-called penalization function, which we fix with~\eqref{eq:sigma} below. The stability of the discontinuous Galerkin method is linked to the correct choice of $\sigma$. Furthermore, it can strongly affect convergence properties when high local variation of geometry or diffusion coefficients occur, as detailed in~\cite{dong2022robust}. A rigorous analysis showing that the method can be made stable and optimally convergent even on the rough grids considered herein is presented in~\cite{Cangiani_PolyDG,dgp,dgease}.

%% file: tests.tex
\subsection{Convergence tests}\label{sec:convergence}

In the following examples, we consider the  Poisson model problem~\ref{eq:problem} posed on domain $\Omega$, with Dirichlet boundary conditions $g$ forced by a smooth manufactured analytical solution $u(\boldsymbol{x}) = \Pi_{i=1}^{d} \sin(\pi x_i)$. We discretize it using the polygonal interior penalty discontinuous Galerkin method described in Section~\ref{sec:polydg}.

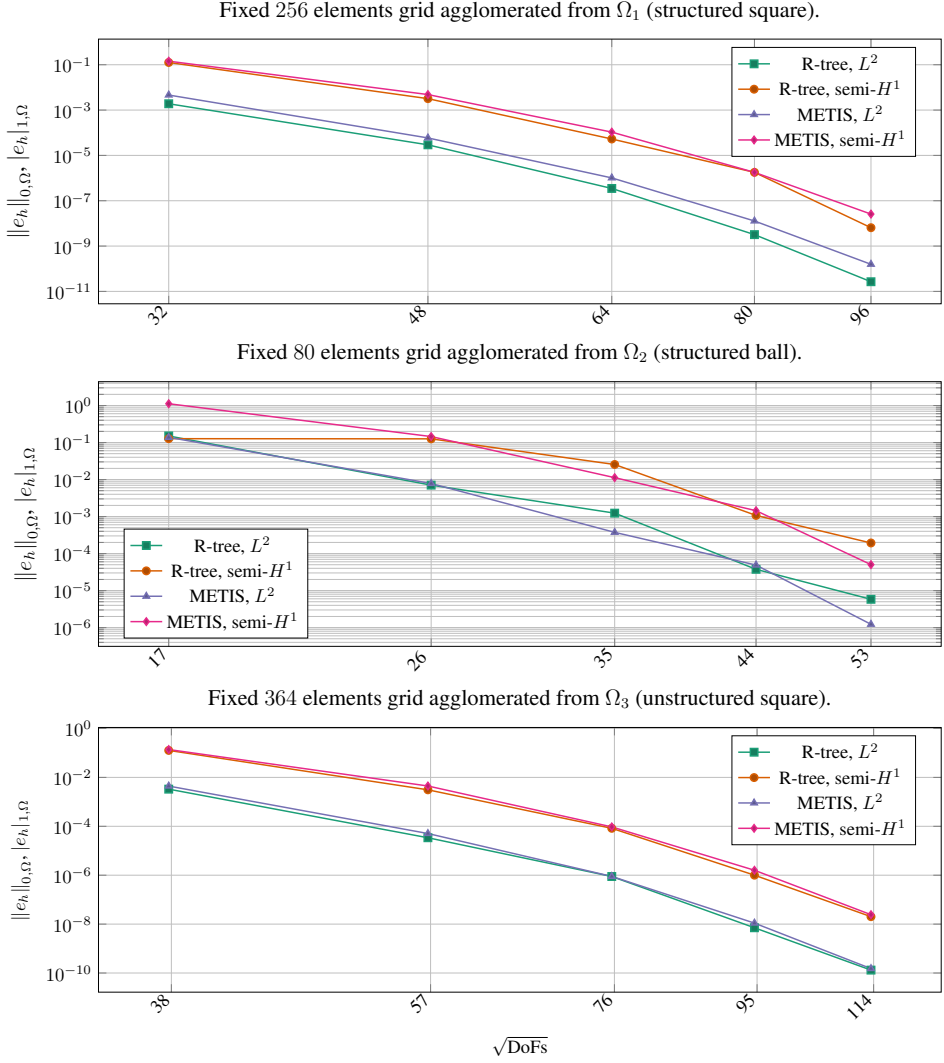
\begin{figure}[h]
  \centering
\begin{tikzpicture}[trim axis left,scale=0.7]
  \begin{axis}[
      scale only axis, %
      height=5cm,
      width=\textwidth, %
      xmode=log,ymode=log,
      grid=both,
      title= \large{Fixed $256$ elements grid agglomerated from $\Omega_1$ (structured square).},
      ylabel=\large{$\norm {e_h}_{0,\Omega}$, $\vert e_h \vert_{1,\Omega}$},
      legend pos=north east,
      x tick label style={rotate=45, anchor=east, align=center},
      xtick={32., 48., 64., 80., 96.},
      xticklabels={32, 48, 64, 80, 96},
      ]
      \addplot+[solid, every mark/.append style={solid}, mark=square*]
      coordinates {
        (32,0.00189974)(48.,2.93253e-05)(64.,3.48542e-07)(80.,3.1999e-09)(96,2.67161e-11)
        };
        \addlegendentry{R-tree, $L^2$};
        \addplot+[solid, every mark/.append style={solid}, mark=*]
        coordinates {
          (32,0.125874)(48,0.00319308)(64,5.29533e-05)(80,1.78721e-06)(96,6.46256e-09)
          };
          \addlegendentry{R-tree, $\text{semi-}H^1$};
          
      \addplot+[solid, every mark/.append style={solid}, mark=triangle*]
      coordinates {
        (32,0.00460829)(48,5.90358e-05)(64,1.01695e-06)(80,1.27898e-08)(96,1.57099e-10)
        };
      \addlegendentry{METIS, $L^2$};
      \addplot+[solid, every mark/.append style={solid}, mark=diamond*]
      coordinates {
        (32,0.143864)(48,0.0048376)(64,0.000106531)(80,1.78721e-06)(96,2.57866e-08)
        };
      \addlegendentry{METIS, $\text{semi-}H^1$};
 
  \end{axis}
\end{tikzpicture}

\begin{tikzpicture}[trim axis left,scale=0.7]
  \begin{axis}[
      scale only axis, %
      height=5cm,
      width=\textwidth, %
      xmode=log,ymode=log,
      grid=both,
      title= \large{Fixed $80$ elements grid agglomerated from $\Omega_2$ (structured ball).},
      ylabel=\large{$\norm {e_h}_{0,\Omega}$, $\vert e_h \vert_{1,\Omega}$},
      legend pos=south west,
      x tick label style={rotate=45, anchor=east, align=center},
      xtick={17,  26 ,35, 44, 53},
      xticklabels={17,  26 ,35, 44, 53},
      ]
      \addplot+[solid, every mark/.append style={solid}, mark=square*]
      coordinates {
        (17, 0.150468)(26, 0.00706559)(35.,0.00123692)(44.,3.75907e-05)(53,5.81712e-06)
        };
      \addlegendentry{R-tree, $L^2$};
      \addplot+[solid, every mark/.append style={solid}, mark=*]
        coordinates {
        (17,0.127004)(26,0.125721)(35,0.025485)(44,0.00107344)(53,0.000193528)
          };
      \addlegendentry{R-tree, $\text{semi-}H^1$};
          
      \addplot+[solid, every mark/.append style={solid}, mark=triangle*]
      coordinates {
        (17,0.138)(26,0.00781356)(35,0.000374667)(44,4.89077e-05)(53,1.22989e-06)
        };
      \addlegendentry{METIS, $L^2$};
      \addplot+[solid, every mark/.append style={solid}, mark=diamond*]
      coordinates {
        (17,1.11689)(26,0.144995)(35,0.0112497)(44,0.00144248)(53,5.03074e-05)
        };
      \addlegendentry{METIS, $\text{semi-}H^1$};
 
  \end{axis}
\end{tikzpicture}

\begin{tikzpicture}[trim axis left,scale=0.7]
  \begin{axis}[
      scale only axis, %
      height=5cm,
      width=\textwidth, %
      xmode=log,ymode=log,
      grid=both,
      title= \large{Fixed $364$ elements grid agglomerated from $\Omega_3$ (unstructured square).},
      xlabel={$\sqrt{\text{DoFs}}$},
      ylabel={$\norm {e_h}_{0,\Omega}$, $\vert e_h \vert_{1,\Omega}$},
      legend pos=north east,
      x tick label style={rotate=45, anchor=east, align=center},
      xtick={38.15756806,  57.23635209,  76.31513611,  95.39392014,
      114.47270417},
      xticklabels={38, 57, 76, 95, 114},
      ]
      \addplot+[solid, every mark/.append style={solid}, mark=square*]
      coordinates {
        (38,0.00330861)(57.,3.40476e-05)(76.,8.84921e-07)(95.,7.08998e-09)(114,1.31514e-10)
        };
      \addlegendentry{R-tree, $L^2$};
      \addplot+[solid, every mark/.append style={solid}, mark=*]
        coordinates {
          (38,0.127004)(57,0.00310559)(76,8.27386e-05)(95,1.00255e-06)(114,2.01561e-08)
          };
      \addlegendentry{R-tree, $\text{semi-}H^1$};
          
      \addplot+[solid, every mark/.append style={solid}, mark=triangle*]
      coordinates {
        (38,0.00442053)(57,5.06846e-05)(76,8.93129e-07)(95,1.09854e-08)(114,1.51489e-10)
        };
      \addlegendentry{METIS, $L^2$};
      \addplot+[solid, every mark/.append style={solid}, mark=diamond*]
      coordinates {
        (38,0.137959)(57,0.00445175)(76,9.45825e-05)(95,1.59371e-06)(114,2.42952e-08)
        };
      \addlegendentry{METIS, $\text{semi-}H^1$};
 
  \end{axis}
\end{tikzpicture}
\caption{Problem~\eqref{eq:problem} with $d=2$. Convergence under $p$-refinement for $p=1,2,3,4,5$.}
\label{fig:p_convergence}
\end{figure}

For any agglomerated polytopic element $K \in \Omega$, we define
on its bounding box $B_K$ the standard polynomial space $\mathcal{Q}^p(B_K)$ spanned by tensor-product Lagrange polynomials of degree $p$ in each variable, denoted with $\{ \phi_i \}_{i=1}^{N_p}$ , with local dimension $N_p=(p+1)^d$. Since $ \overline{K} \subset \overline{B_K}$, the basis
on $K$ is defined by restricting each basis function to $K$. Standard Gauß-Lobatto quadrature rules of order $2p+1$ are defined on the already available sub-tessellation $\{\tau_K\}$
of $K$ part of the underlying grid. Several other choices are possible, and we refer to the textbook \cite{Cangiani_PolyDG} for a discussion about relevant implementation details. 
We stress that in all comparisons presented below the only difference lies in the choice of the agglomeration strategy, while the number of degrees of freedom
is identical as the number of total partitions is the same. For each domain, we investigate convergence under $p$ refinement in both the $L^2$-norm and the $H^1$-seminorm. The penalization function in~\eqref{eqn:weak_form} is fixed as
\begin{equation}\label{eq:sigma}
\sigma(\boldsymbol{x}) = C_\sigma
\begin{cases}
\frac{p^2}{h_K} & \text{ on } f \in \Gamma \cap \partial \Omega, \\
\frac{p^2}{\min\{h_K^+,h_K^-\}} &  \text{ on } f \in \Gamma_{\dint},
\end{cases}
\end{equation}
where in the forthcoming experiments we have set $C_\sigma=10$. In Figure \ref{fig:p_convergence} we display errors in function of the polynomial degree $p$ for fixed grids obtained by agglomeration of the  two-dimensional grids $\Omega_1,\Omega_2$ and $\Omega_3$. We repeat the same procedure starting from the unit cube grid $\Omega_4$ consisting of 32,768 hexahedra, which is then agglomerated to generate coarse
polytopic grids $\mathcal{T}$ obtained with either METIS or R-tree and comprising only 72 polytopes. The respective findings are reported in Figure \ref{fig:p_convergence_cube}. In all cases, we observe the exponential convergence predicted by the theory~\cite{Cangiani_PolyDG,dgp,dgease}. METIS and R-tree results are always very close to each other, with the exclusion of the structure square case the perfectly square mesh produced by the R-tree yields marginally superior results. In any case our results confirm the robustness of the polygonal discontinuous Galerkin method over very general polytopal grids. 

The fact that each polytope is close to its enclosing bounding box motivates the usa of tensor product polynomials of degree $p$ in each coordinate direction.
Nevertheless, it is known that the usage of $\mathcal{Q}_p$ spaces on polytopic elements increases the number of degrees of freedom without substantial gain, and for this reason the $\mathcal{P}_p$ space of polynomials of \emph{total} degree $p$
can be employed. The effect of the use of $\mathcal{P}_p$ basis functions is shown in Figure~\ref{fig:comparison_dgq_dgp}, where convergence under $p$-refinement is investigated for the unstructured mesh $\Omega_3$. On the left, we can observe
the same patterns described in~\cite{CGH14}, highlighting the benefit of using $\mathcal{P}_p$ elements: higher accuracy is achieved for a comparable number of degrees of freedom with respect to tensor product
elements, thanks to the higher-order local polynomial space. Identical patterns have been found also for the $H^1$ seminorm, and are omitted for readability. On the right, we again observe exponential convergence under $p$-refinement also for $\mathcal{P}_p$ elements.

\begin{figure}[h]
  \centering
  \begin{tikzpicture}[trim axis left,scale=0.8]
    \begin{axis}[
        scale only axis, %
        height=7cm,
        width=\textwidth, %
        xmode=log,ymode=log,
        grid=both,
        xlabel={$\sqrt{\text{DoFs}}$},
        ylabel=\large{$\norm {e_h}_{0,\Omega}$, $\vert e_h \vert_{1,\Omega}$},
        legend pos=south west,
        x tick label style={rotate=45, anchor=east, align=center},
        xtick={ 9.59166305 , 19.07878403 , 38.15756806 ,76.31513611, 152.63027223, 305.26},
        xticklabels={9,  19 ,38, 76, 152, 305},
        ]
        \addplot+[dashed, every mark/.append style={solid}, mark=pentagon*]
        coordinates {
          (9, 4.447668e-02)
          (19,1.269757e-02)
          (38, 3.140411e-03)
          (76,7.887963e-04)
          (152,1.854817e-04)
          };
        \addlegendentry{R-tree, $L^2$, $\mathcal{Q}_1$};
        \addplot+[dashed, every mark/.append style={solid}, mark=square*]
        coordinates {
          (14, 1.797703e-03)
          (28,2.489369e-04)
          (57, 3.377351e-05)
          (114,4.405226e-06)
          (228,5.558035e-07)
          };
        \addlegendentry{R-tree, $L^2$, $\mathcal{Q}_2$};

        \addplot+[dashed, every mark/.append style={solid}, mark=oplus*]
        coordinates {
          (19,  1.429895e-04)
          (38, 1.021753e-05)
          (76, 7.321601e-07)
          (152, 5.126069e-08)
          (305,3.919937e-09)
          };
        \addlegendentry{R-tree, $L^2$, $\mathcal{Q}_3$};
            
        \node at (axis cs:154, 1.0e-4) [anchor=south west] {$\sim h^{-2}$};
        \addplot+[solid, every mark/.append style={solid}, mark=pentagon*]
        coordinates {
          (9,1.157225e-01)
          (19,1.581105e-02)
          (38,4.252653e-03)
          (76,1.142052e-03)
          (152,2.860103e-04)
          };
          \addlegendentry{METIS, $L^2$, $\mathcal{Q}_1$};
   
        \node at (axis cs:230, 4.0e-7) [anchor=south west] {$\sim h^{-3}$};
        \addplot+[solid, every mark/.append style={solid}, mark=square*]
        coordinates {
          (14,2.164270e-03)
          (28,3.984938e-04)
          (57,5.428619e-05)
          (114,7.033799e-06)
          (228,8.471648e-07)
          };
          \addlegendentry{METIS, $L^2$, $\mathcal{Q}_2$};

        \node at (axis cs:309, 2.4e-9) [anchor=south west] {$\sim h^{-4}$};
        \addplot+[solid, every mark/.append style={solid}, mark=oplus*]
        coordinates {
          (19,  1.710574e-04)
          (38, 9.541160e-06)
          (76, 7.302473e-07)
          (152, 4.946087e-08)
          (305,3.353402e-09)
          };
          \addlegendentry{METIS, $L^2$, $\mathcal{Q}_3$};
   
    \end{axis}
  \end{tikzpicture}
  \caption{Convergence under $h$-refinement for agglomerated mesh $\Omega_3$ with $\mathcal{Q}_p$ elements ($p=1,2,3$), and different agglomeration strategies.}\label{fig:h_refinement}
\end{figure}
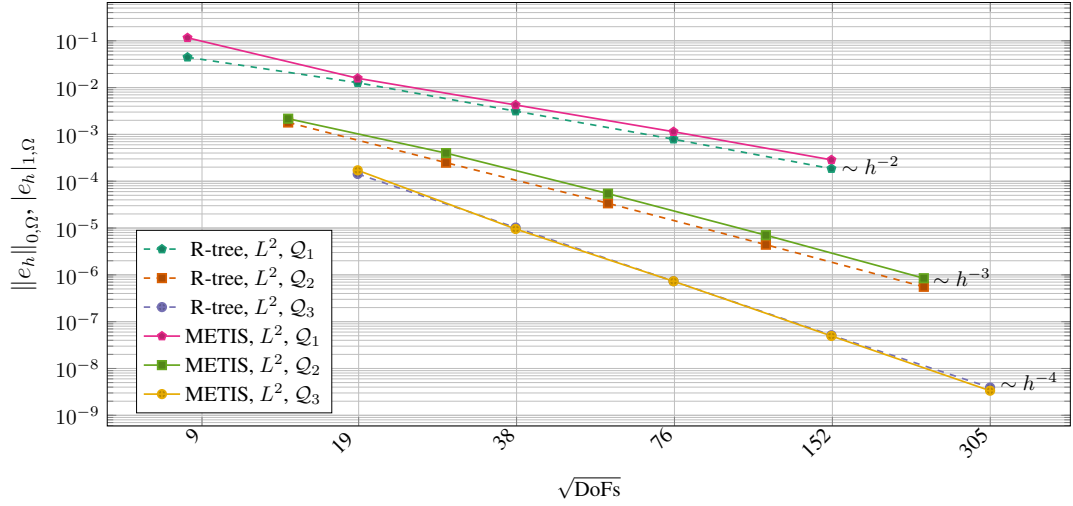

Finally, we assess convergence under $h$-refinement considering the mesh $\Omega_3$ and varying the polynomial degree $p$ and the agglomeration strategy. We display in Figure~\ref{fig:h_refinement} the $L^2$ and $H^1$ seminorm errors with respect to the same manufactured solution $u$ employed for the $p$-convergence test. We observe optimal rates for all the polynomial degrees and the chosen agglomeration procedures. Remarkably, the curves associated to the R-tree approach are always lower or equal to the ones associated with \textsc{METIS} partitioning.

\begin{figure}[h]
  \centering
\begin{tikzpicture}[trim axis left,scale=0.7]
  \begin{axis}[
      scale only axis, %
      height=5cm,
      width=\textwidth, %
      xmode=log,ymode=log,
      grid=both,
      title= \normalsize{Fixed $72$ elements grid agglomerated from $\Omega_5$ (structured cube).},
      xlabel={$\text{DoFs}^{1/3}$},
      ylabel={$\norm {e_h}_{0,\Omega}$, $\vert e_h \vert_{1,\Omega}$},
      legend pos=north east,
      x tick label style={rotate=45, anchor=east, align=center},
      xtick={8.32033529, 12.48050294, 16.64067058, 20.80083823},
      xticklabels={8, 12, 16, 20},
      ]
      \addplot+[solid, every mark/.append style={solid}, mark=square*]
      coordinates {
        (8,0.0339277)(12,0.00268563)(16.,0.000245186)(20.,1.13356e-05)
        };
      \addlegendentry{R-tree, $L^2$};
      \addplot+[solid, every mark/.append style={solid}, mark=*]
        coordinates {
          (8,0.507079)(12,0.068158)(16.,0.00706919)(20.,0.00043592)
          };
      \addlegendentry{R-tree, $\text{semi-}H^1$};
          
      \addplot+[solid, every mark/.append style={solid}, mark=triangle*]
      coordinates {
        (8,0.028252)(12,0.00245023)(16.,0.000190443)(20.,1.04804e-05)
        };
      \addlegendentry{METIS, $L^2$};
      \addplot+[solid, every mark/.append style={solid}, mark=diamond*]
      coordinates {
        (8,0.468394)(12,0.06474)(16.,0.00583337)(20.,0.000410002)
        };
      \addlegendentry{METIS, $\text{semi-}H^1$};
 
  \end{axis}
\end{tikzpicture}
\caption{Problem~\eqref{eq:problem} with $d=3$. Convergence under $p$-refinement for $p=1,2,3,4$.}
\label{fig:p_convergence_cube}
\end{figure}
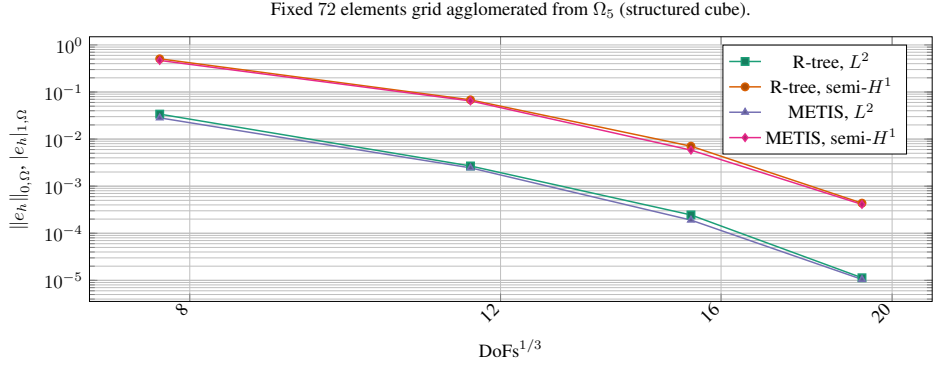

\begin{figure}[h]
  \centering
\begin{tikzpicture}[trim axis left,scale=0.5]
  \begin{axis}[
      scale only axis, %
      height=8cm,
      width=0.8\textwidth, %
      xmode=log,ymode=log,
      grid=both,
      title= \large{$p$-convergence: comparison between $Q_p$ and $\mathcal{P}_p$ elements.},
      xlabel={$\sqrt{\text{DoFs}}$},
      ylabel=\large{$\norm {e_h}_{0,\Omega}$, $\vert e_h \vert_{1,\Omega}$},
      legend pos=south west,
      x tick label style={rotate=45, anchor=east, align=center},
      xtick={38,57,76,95,114,133},
      xticklabels={38,57,76,95,114,133},
      ]
      \addplot+[solid, every mark/.append style={solid}, mark=square*]
      coordinates {
        (38, 3.305609e-03)
        (57.,3.404810e-05)
        (76.,8.849206e-07)
        (95.,7.089978e-09)
        (114, 1.315136e-10)
        (133, 9.782926e-13)
        };
        \addlegendentry{R-tree, $L^2 (Q_p)$};
        \addplot+[solid, every mark/.append style={solid}, mark=*]
        coordinates {
          (33, 4.132496e-02)
          (46, 1.474958e-04)
          (60.,3.253502e-06)
          (73.,6.479461e-08)
          (87, 8.542262e-10)
          (100,1.351629e-11)
        };
        \addlegendentry{R-tree, $L^2 (\mathcal{P}_p)$};

  \end{axis}
\end{tikzpicture}
\qquad %
\begin{tikzpicture}[trim axis left,scale=0.5]
  \begin{axis}[
      scale only axis, %
      height=8cm,
      width=0.8\textwidth, %
      xmode=log,ymode=log,
      grid=both,
      title= \large{$p$-convergence for $\mathcal{P}_p$ elements.},
      xlabel={$\sqrt{\text{DoFs}}$},
      legend pos=south west,
      x tick label style={rotate=45, anchor=east, align=center},
      xtick={33, 46, 60, 73, 87, 100},
      xticklabels={33, 46, 60, 73, 87, 100},
      ]
      \addplot+[solid, every mark/.append style={solid}, mark=square*]
      coordinates {
        (33, 4.132496e-02)(46, 1.474958e-04)(60.,3.253502e-06)(73.,6.479461e-08)(87,8.542262e-10)(100,1.351629e-11)
        };
      \addlegendentry{R-tree, $L^2 (\mathcal{P}_p$)};
      \addplot+[solid, every mark/.append style={solid}, mark=*]
        coordinates {
          (33, 2.541862e-01)
          (46, 1.048612e-02)
          (60.,3.360319e-04)
          (73.,7.606444e-06)
          (87, 1.360712e-07)
          (100,2.293088e-09)
          };
          \addlegendentry{R-tree, $\text{semi-}H^1 (\mathcal{P}_p)$};
          
          \addplot+[solid, every mark/.append style={solid}, mark=triangle*]
          coordinates {
            (33, 2.362214e-02)
            (46, 1.565300e-04)
            (60.,4.292490e-06)
            (73.,8.383538e-08)
            (87, 1.856890e-09)
            (100,3.499751e-11)
            };
            \addlegendentry{METIS, $L^2 (\mathcal{P}_p)$};
            \addplot+[solid, every mark/.append style={solid}, mark=diamond*]
            coordinates {
            (33, 2.536485e-01)
            (46, 1.011700e-02)
            (60.,3.996037e-04)
            (73.,9.252206e-06)
            (87, 2.446019e-07)
            (100,5.083452e-09)
        };
      \addlegendentry{METIS, $\text{semi-}H^1 (\mathcal{P}_p)$};
  \end{axis}
\end{tikzpicture}
\caption{Fixed $364$ elements grid agglomerated from $\Omega_3$ (unstructured square). Convergence under $p$-refinement for $p=1,2,3,4,5,6$ and different spaces.}
\label{fig:comparison_dgq_dgp}
\end{figure}
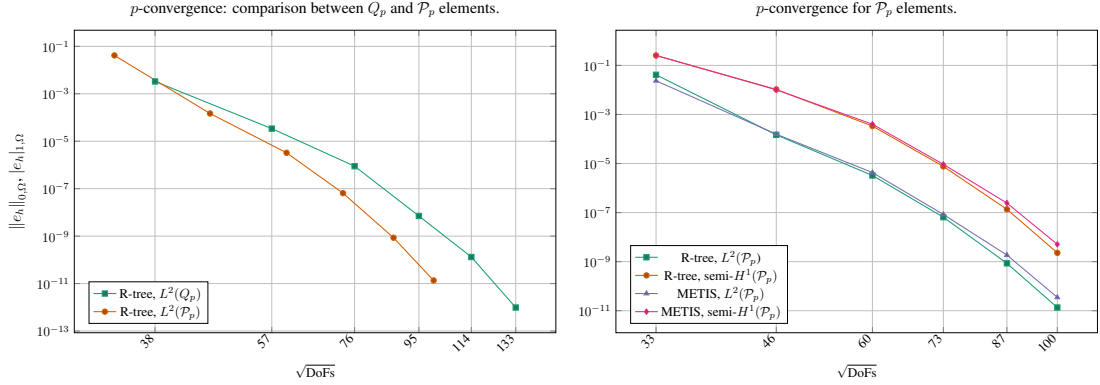

%% file: applications.tex
\label{sec:applications}

\subsubsection{3D geometry}\label{subsubsec:piston}
\begin{figure}[!htb]
  \centering
  \includegraphics[width=7.5cm]{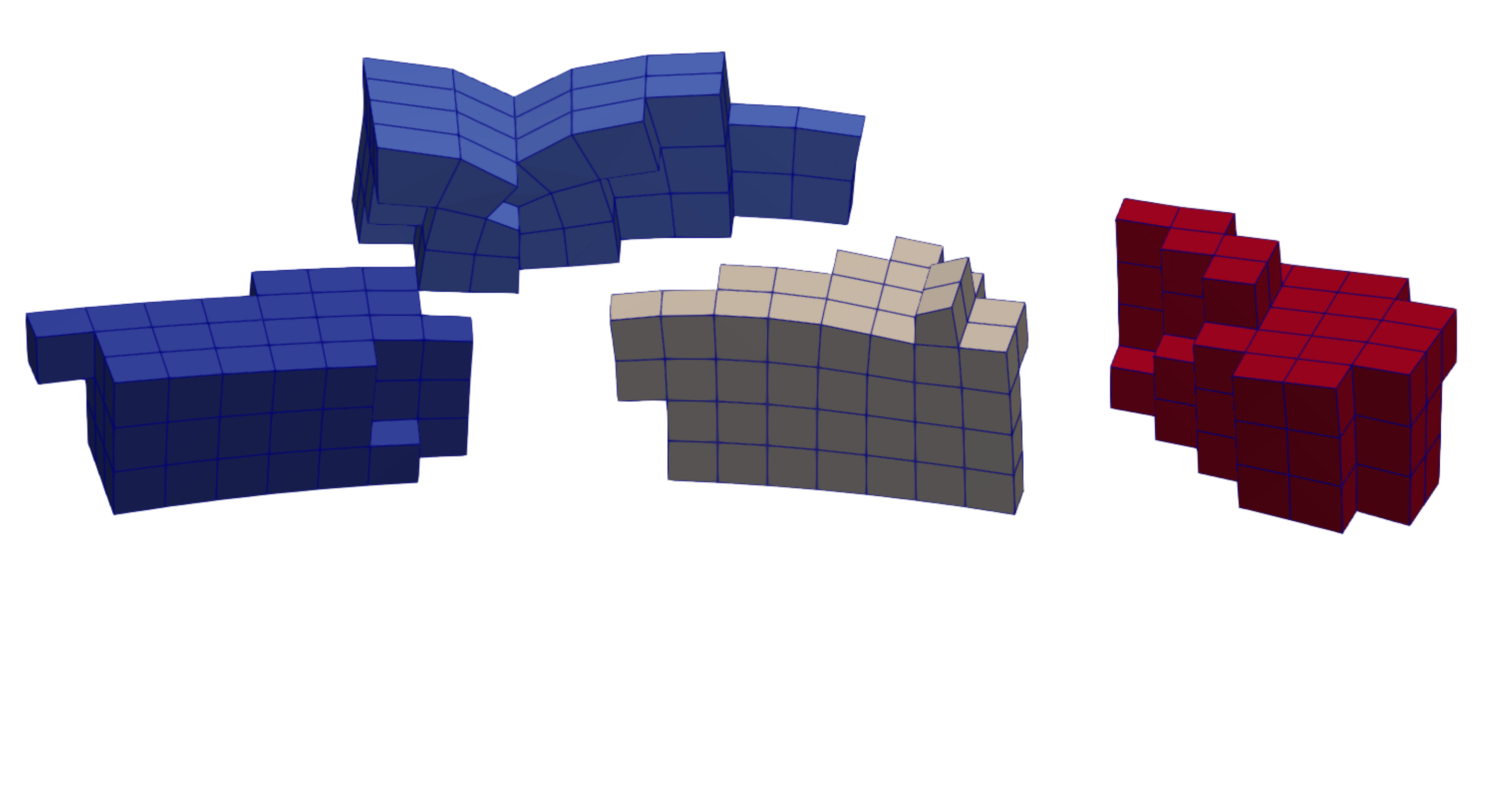}
  \hfill
  \includegraphics[width=8cm]{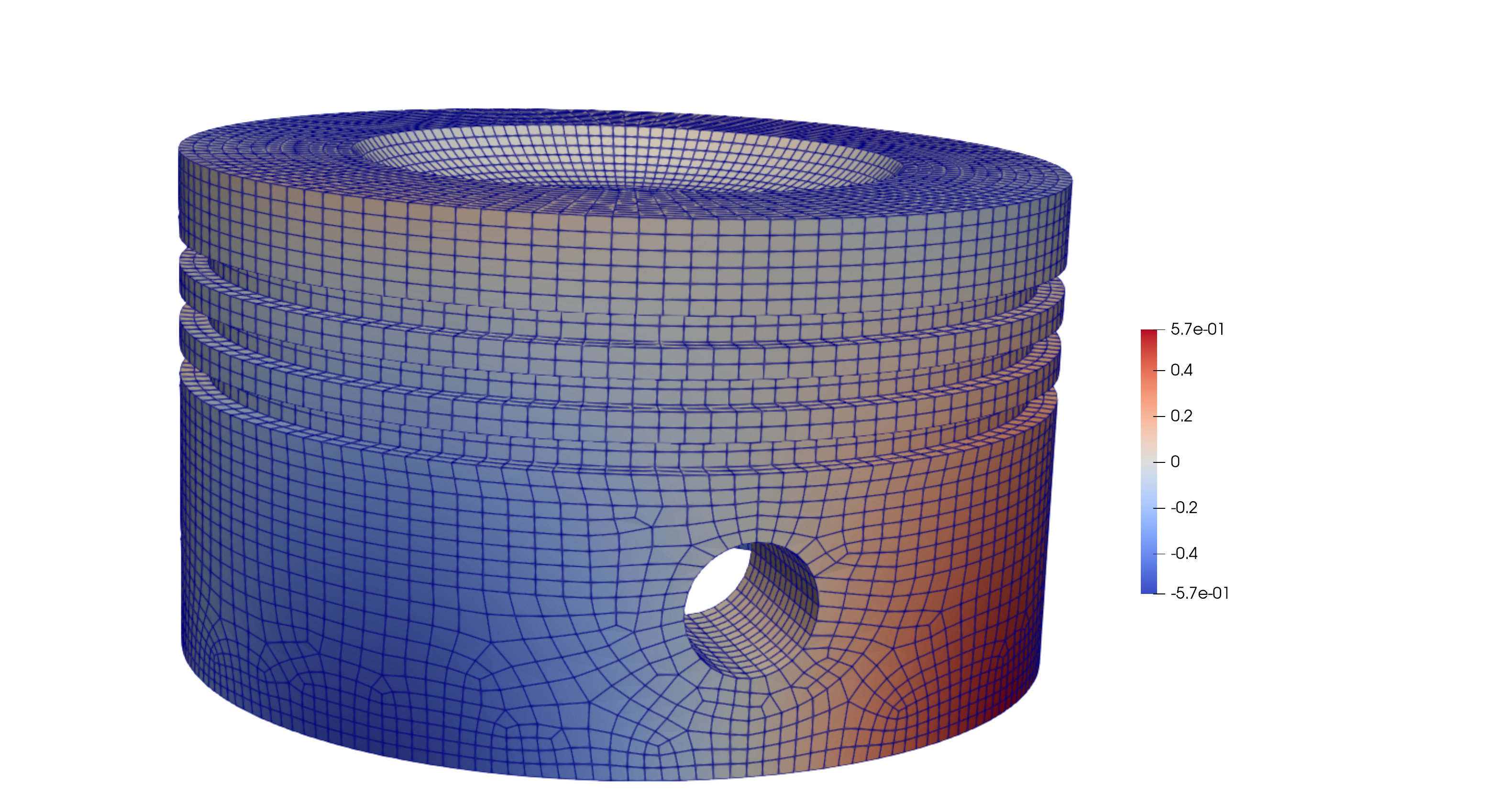}
  \caption{Left: sample agglomerates generated by the R-tree algorithm. Right: view of the solution.}
  \label{fig:piston_test}
\end{figure}

\begin{figure}[!htb]
  \centering
\begin{tikzpicture}[trim axis left,scale=0.82]
  \begin{axis}[
      scale only axis, %
      height=5cm,
      width=\textwidth, %
      xmode=log,ymode=log,
      grid=both,
      title= \normalsize{$\Omega_6$ ($3$D piston model)},
      xlabel={$\text{DoFs}^{1/3}$},
      ylabel={$\norm {e_h}_{0,\Omega}$, $\vert e_h \vert_{1,\Omega}$},
      legend pos=north east,
      x tick label style={rotate=45, anchor=east, align=center},
      xtick={18., 27., 36., 45.},
      xticklabels={18, 27, 36, 45},
      ]
      \addplot+[solid, every mark/.append style={solid}, mark=square*]
      coordinates {
        (18,0.0142143)(27,0.000500712)(36.,7.0761e-06)(45.,2.32691e-07)
        };
      \addlegendentry{R-tree, $L^2$};
      \addplot+[solid, every mark/.append style={solid}, mark=*]
        coordinates {
          (18,0.14325)(27,0.00483883)(36.,7.99999e-05)(45.,3.08876e-06)
          };
      \addlegendentry{R-tree, $\text{semi-}H^1$};
          
      \addplot+[solid, every mark/.append style={solid}, mark=triangle*]
      coordinates {
        (18,0.00870615)(27,0.000343421)(36.,3.63493e-06)(45.,1.45107e-07)
        };
      \addlegendentry{METIS, $L^2$};
      \addplot+[solid, every mark/.append style={solid}, mark=diamond*]
      coordinates {
        (18,0.0869164)(27,0.00411964)(36.,5.23425e-05)(45.,3.08053e-06)
        };
      \addlegendentry{METIS, $\text{semi-}H^1$};
 
  \end{axis}
\end{tikzpicture}
\caption{Convergence under $p$-refinement for the piston test case for $p=1,2,3,4$. Fixed polytopal grids made of $731$ agglomerates.}
\label{fig:p_convergence_piston}
\end{figure}
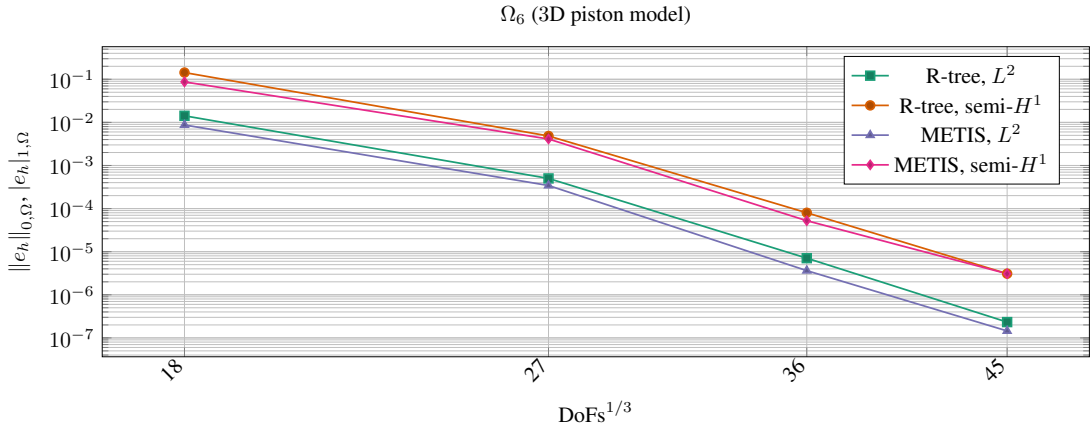

We finally compare the two strategies on the piston mesh $\Omega_6$, a three-dimensional example with a \emph{non-trivial} geometry which cannot be meshed with a low number of cells, due to the presence of local features requiring a fine mesh to be captured accurately. %
Such mesh comes from a real three-dimensional CAD model which has been first repaired and later meshed with hexahedra through the commercial software \textsc{CUBIT}. We adopt, on this
fixed mesh, both the METIS and R-tree agglomeration strategies, and investigate again convergence under refinement of the polynomial degree $p$, using as manufactured solution $u(\boldsymbol{x})=\Pi_{i=1}^3 \sin(\pi \frac{x_i}{10})$. Even
if agglomeration strategies strongly decrease the number of DoFs, sparse direct solvers quickly become prohibitive for non-trivial three-dimensional geometries with moderately high degrees, such as $p\geq 3$,
and modest mesh sizes. We report in Table \ref{tab:DoFs_per_degree_piston} the number of DoFs in function of the polynomial degree required by the discontinuous Galerkin method with the original hexahedral model and the agglomerated version. The level of agglomeration is fixed so as to keep the number of DoFs within a regime where sparse direct solvers are still effective. This allows us to solve the resulting linear system with the \textsc{MUMPS} solver \cite{MUMPS} through the \textsc{Trilinos} library \cite{Trilinos}. The original background grid $\Omega_5$ consists
of $45,714$ hexahedra. After generating the R-tree data structure as outlined in Section \ref{sec:rtree}, we extract its second level, resulting in a polytopal mesh with only $731$ elements, thereby reducing the size of the original problem by a factor of $64$. Following
the procedure of the previous experiments, we employ METIS by setting $731$ as the target number of partitions. A view of some elements of the coarse mesh $\mathcal{T}$ generated using
Algorithm \ref{alg:r_tree_algo}, as well as a view of the solution $u$ interpolated onto the finer mesh, are shown in Figure \ref{fig:piston_test}. The convergence history
is reported in Figure~ \ref{fig:p_convergence_piston}; it can be appreciated how both approaches yield comparable results in terms of accuracy also
for the present configuration, with METIS giving slightly more accurate results in the $L^2$-norm.

\begin{table}[!t]
  \begin{center}
      \centering
      \begin{tabular}{|r|r|r|}%
   \hline
   & \multicolumn{2}{|c|}{{Number of DoFs } } \\
   \hline
   $p$  &  Number of DoFs (original mesh) & Number of DoFs (agglomerated mesh)  \\ %
   \hline
    1   &  365,712    & 25,425  \\
    2   &  1,234,278  & 48,366  \\
    3   & 2,925,696   & 108,492 \\
    4   & 5,714,250   & 175,020 \\
   \hline
  \end{tabular}
\end{center}
      \caption{3D piston model. Total number of DoFs in function of the  polynomial degree $p$.}
      \label{tab:DoFs_per_degree_piston}
\end{table}

\subsection{Multigrid preconditioning}

We finally introduce R-tree based MultiGrid (\textsc{R3MG}) preconditioning. The R-tree agglomeration algorithm naturally produces \emph{nested} hierarchies of agglomerated grids. We exploit these grids and the flexibility of the discontinuous Galerkin framework to construct multigrid preconditioners. For an analysis of multigrid solvers applied to polygonal discontinuous Galerkin methods we refer to~\cite{AntoniettiPennesi,hpMGAntonietti}.

We consider again the model problem~\eqref{eq:problem} with $f=1$ and homogeneous Dirichlet boundary conditions. We use multigrid as preconditioner for the conjugate-gradient solver (CG) \cite{HestenesMCG} with one multigrid cycle per iteration
as this is known to be usually more robust than using multigrid as a solver. We stress that a sequence of \emph{nested} agglomerated grids $\{ \mathcal{T}_l \}_{l=1}^M$ can be directly generated thanks
to the structure of the tree. Denoting with $V_h^l$ the finite-dimensional discontinuous space defined on $\mathcal{T}_l$, the sequence of nested grids induces a nested sequence of spaces
$V_h^1 \subset V_h^2 \subset \ldots \subset V_h^M$. Thanks to this property, the intergrid transfer operators (\emph{restriction} and \emph{prolongation}) are more easily defined and cheaper to compute 
compared to a non-nested version. Indeed,  the latter requires the computation of expensive $L^2$ projections over arbitrarily intersecting meshes.

The prolongation operator between the spaces $V_{l-1}$ and $V_l$ is denoted by $\mathcal{P}_{l-1}^l$ and consists of the natural injection operator, $\mathcal{P}_{l-1}^l: V_{l-1} \xhookrightarrow{} V_l$, while as restriction operator we choose 
$\mathcal{R}_{l}^{l-1} \coloneqq \bigl( \mathcal{P}_{l-1}^l \bigr)^T$. The action of such operators can be performed matrix-free without
the need to explicitly store the sparse matrices associated to intergrid transfers. For each domain, the experiments are configured in the following way:
\begin{itemize}
    \item The conjugate-gradient solver is run with $\mathtt{abstol}=10^{-12}$ and $\mathtt{reltol}=10^{-9}$;
    \item The conjugate-gradient solver is preconditioned by a single V-cycle of multigrid;
    \item As pre- and post-smoothers,  $10$ steps of a Chebyshev smoother of degree $5$ is employed, using eigenvalue estimates computed with 20 iterations of the Lanczos iteration.
    \item As coarse grid solver, a direct solver is used.
\end{itemize}

We start by considering a refined version of the grids shown in Figures~\ref{fig:structured_square}, \ref{fig:structured_ball}, and \ref{fig:unstructured_square}, in order to produce a large enough number of agglomerated levels. We report the iteration counts when varying the number of levels $l$, and the polynomial degree $p$, and the number of degrees of freedom associated to each level of the hierarchy. To show the effectiveness of the multigrid preconditioner, we also report the number of iterations needed by the conjugate-gradient method without preconditioning applied to
the finest level of the multilevel hierarchy, which coincides with the original computational mesh.

In Tables~\ref{tab:2D_examples_MG} we report iteration counts obtained for the two-dimensional geometries in function of the number of levels $l$ and the polynomial degree $p$ from $1$ to $3$. We observe a roughly constant number of iterations when MG is employed as a preconditioner. The iterations with plain CG are significantly higher when no preconditioner is used. This is expected since CG is applied to the finest
level of the hierarchy. In this scenario, the benefit of a preconditioner is even more evident as the number of iterations required by plain CG becomes soon order of magnitudes larger compared to the preconditioned version. With $p=3$, in particular, plain CG does not achieve convergence within $10^5$ iterations with the mesh $\Omega_3$.

Finally, we test our approach on the three-dimensional grids $\Omega_5$ and $\Omega_6$, using the same experimental setup. The findings are reported in Table~\ref{tab:3D_examples_MG}.
We observe that even for a non-trivial geometry such as $\Omega_6$ we get a consistent reduction in the number of iterations required by CG for every polynomial degree $p$.

\begin{figure}[h]
 \begin{minipage}{0.55\textwidth}
    \begin{subfigure}{\linewidth}
  \vspace{2cm}
   \end{subfigure}
  \begin{subfigure}{\linewidth}
  \includegraphics[width=1.0\textwidth]{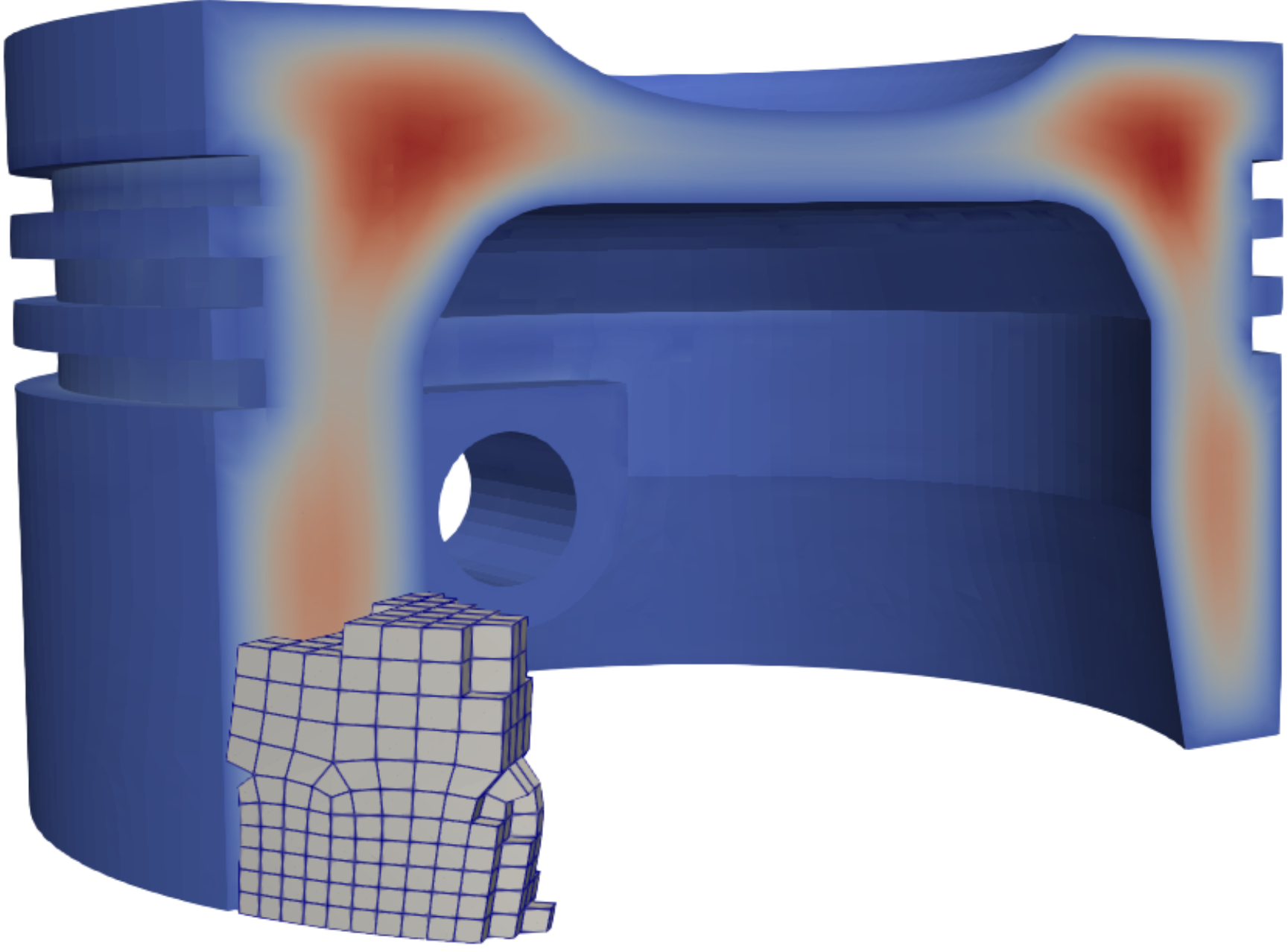}
  \end{subfigure}
  \end{minipage}
  \hspace*{3cm}
  \begin{minipage}{0.4\textwidth} 
  \begin{subfigure}{\textwidth}
  \includegraphics[width=0.6\linewidth]{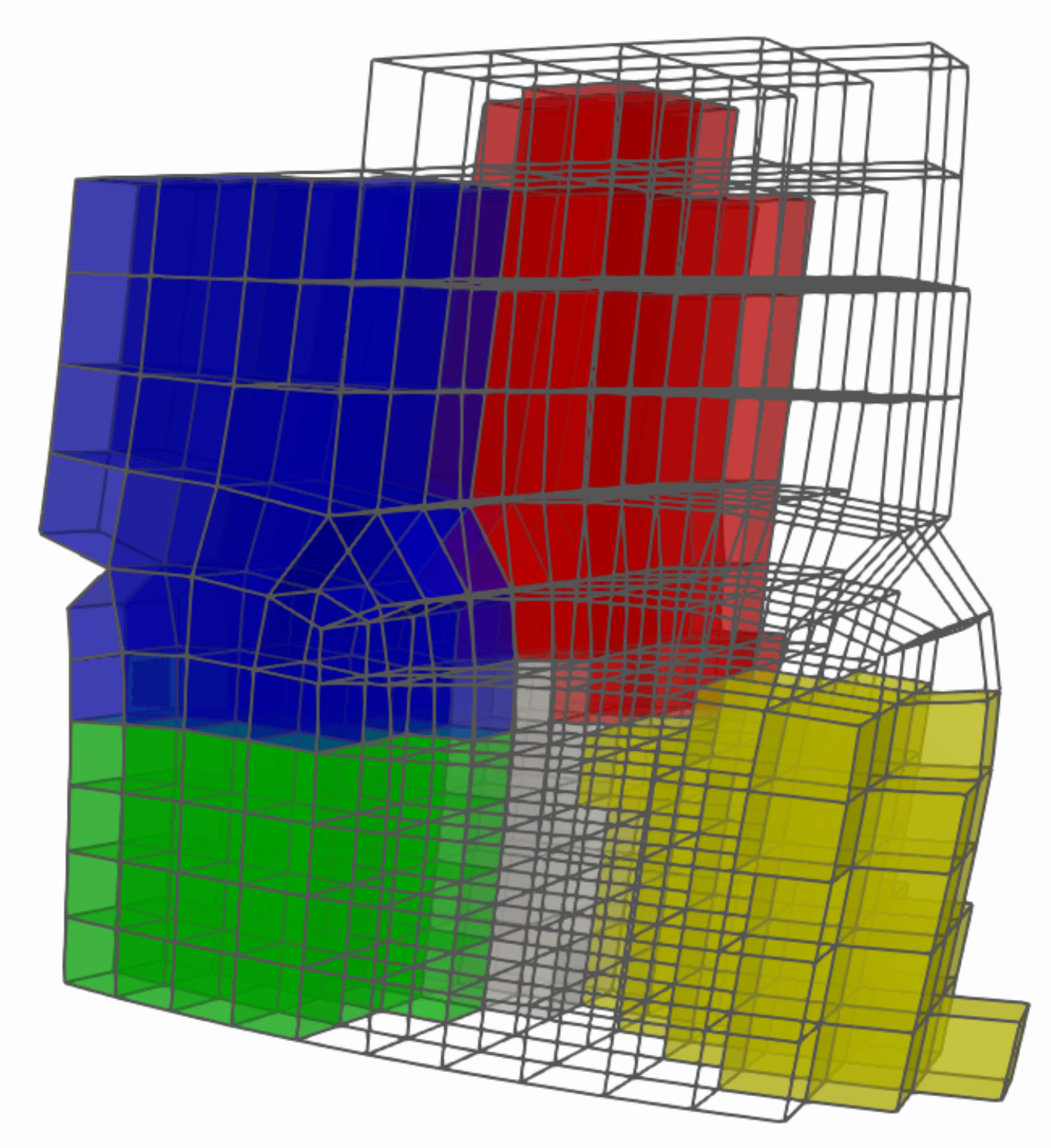}
  \end{subfigure}
  
  \vspace*{0.4cm}
  \begin{subfigure}{\textwidth}
  \includegraphics[width=0.55\linewidth]{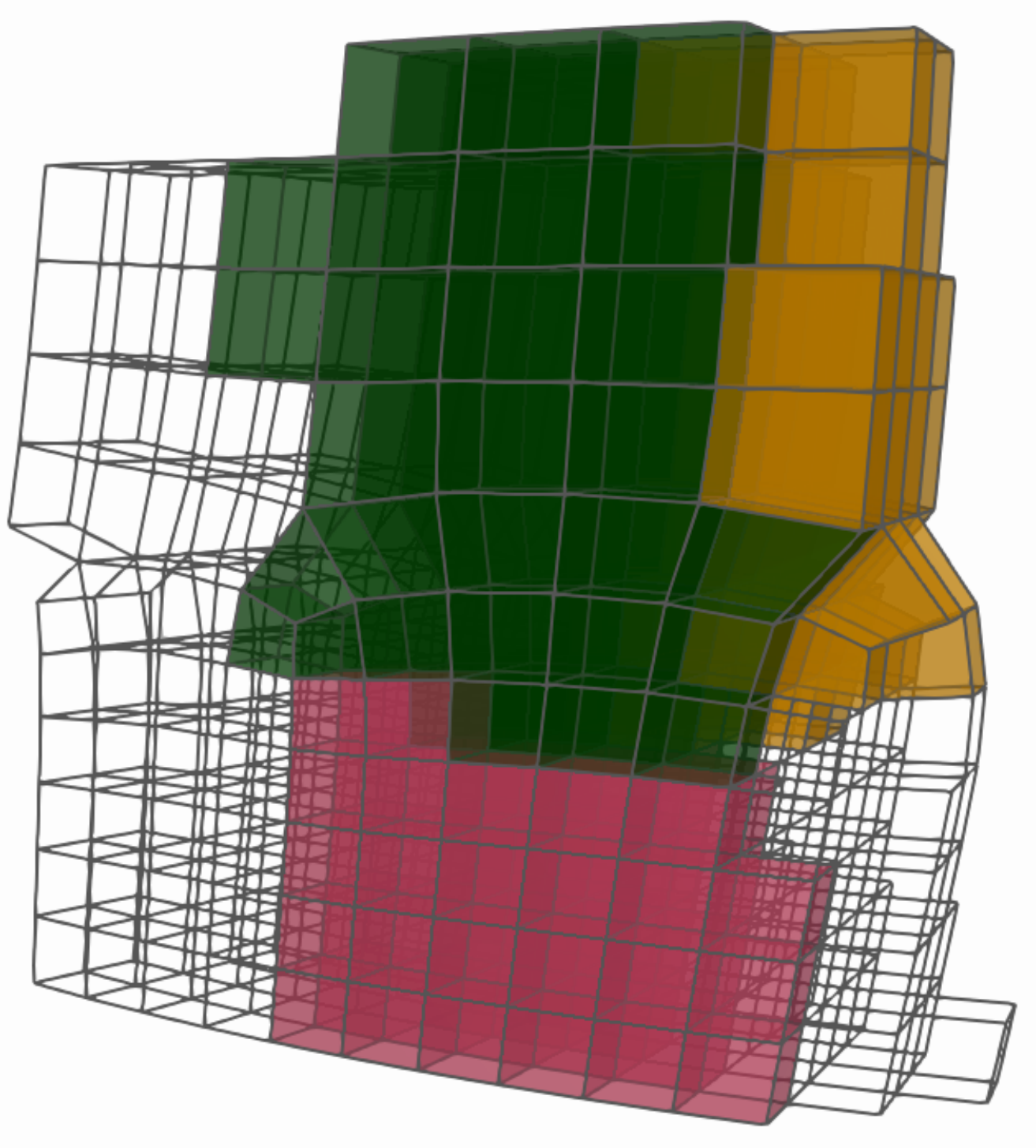}
  \end{subfigure}
  \end{minipage}
  \caption{Left: Clip of the piston geometry $\Omega_6$ with colour plot of the solution computed with $p=1$. In the foreground, highlighted in grey, a single element out of the $90$ elements of a coarse level grid. Right: detailed view of the highlighted element  (wireframe) and its $8$ sub-agglomerates belonging to the finer level. For better visualization, the sub-agglomerates, identified by different colors, are displayed in either of the two plots.}
  \label{fig:sub_agglomerates_piston}
  \end{figure}

We display in Figure~\ref{fig:sub_agglomerates_piston} a coarse three-dimensional polytopal element and its sub-agglomerates for the piston mesh $\Omega_6$, giving a pictorial representation of the capability of generating coarse elements even for complex geometries.

\begin{table}[!htb]
\tiny\setlength{\tabcolsep}{2pt}
\begin{minipage}{.3\textwidth}
\centering

\caption*{$\Omega_1$}%

\begin{tabular}{l|rrc|rrc|rrc}
  \toprule
  \emph{l}& \multicolumn{3}{|c}{$p=1$} & \multicolumn{3}{|c}{$p=2$} & \multicolumn{3}{|c}{$p=3$}  \\
  \cmidrule(lr){2-4}\cmidrule(lr){5-7}\cmidrule(lr){8-10}
  & \emph{\#i} & \#DoF & CG & \emph{\#i}  & \#DoF & CG & \emph{\#i} & \#DoF & CG   \\
  \midrule
  4 & 6 & 256   & -   & 5 & 576   & -   & 6 & 1024  & - \\
  5 & 6 & 1024  & -   & 6 & 2304  & -   & 7 & 4096  & - \\
  6 & 6 & 4096  & -   & 6 & 9216  & -   & 7 & 16384 & -  \\
  7 & 6 & 16384 & 253 & 6 & 36864 & 461 & 7 & 65536 & 682  \\
  \bottomrule
\end{tabular}
\end{minipage}\quad
\begin{minipage}{.3\textwidth}
\centering

\caption*{$\Omega_2$}%

\begin{tabular}{l|rrc|rrc|rrc}
  \toprule
  \emph{l}& \multicolumn{3}{|c}{$p=1$} & \multicolumn{3}{|c}{$p=2$} & \multicolumn{3}{|c}{$p=3$}  \\
  \cmidrule(lr){2-4}\cmidrule(lr){5-7}\cmidrule(lr){8-10}
  & \emph{\#i} & \#DoF & CG & \emph{\#i}  & \#DoF & CG & \emph{\#i} & \#DoF & CG   \\
  \midrule
  4 & 14 & 320   & -   & 12 & 720   & -    & 13    & 1280  & - \\
  5 & 14 & 1280  & -   & 12 & 2880  & -    & 14    & 5120  & - \\
  6 & 14 & 5120  & -   & 12 & 11520 & -    & 14    & 20480 & - \\
  7 & 14 & 20480 & 405 & 12 & 46080 & 1003 & 14 & 81920 & 1603 \\
  \bottomrule
\end{tabular}  
\end{minipage}
\begin{minipage}{.4\textwidth}

\centering

\caption*{$\Omega_3$}%

\begin{tabular}{l|rrc|rrc|rrc}
  \toprule
  \emph{l}& \multicolumn{3}{|c}{$p=1$} & \multicolumn{3}{|c}{$p=2$} & \multicolumn{3}{|c}{$p=3$}  \\
  \cmidrule(lr){2-4}\cmidrule(lr){5-7}\cmidrule(lr){8-10}
  & \emph{\#i} & \#DoF & CG & \emph{\#i}  & \#DoF & CG & \emph{\#i} & \#DoF & CG   \\
  \midrule
  6 & 17 & 5824   & -     & 12 & 13104  & -    & 14 & 23296   & - \\
  7 & 17 & 23296  & -     & 13 & 52416  & -    & 14 & 93184   & - \\
  8 & 17 & 93184  & -     & 13 & 209664 & -    & 14 & 372736  & -  \\
  9 & 17 & 372736 & 3046  & 13 & 838656 & 7405 & 14 & 1490944 & $>\num{1e5}$   \\
  \bottomrule
\end{tabular}  
\end{minipage} 
\caption{Number of required iterations, DoF per level, and iterations required by plain conjugate-gradient (CG) for grids $\Omega_1$, $\Omega_2$, and $\Omega_3$.}\label{tab:2D_examples_MG}
\end{table}

\begin{table}[!htb]
  \scriptsize\setlength{\tabcolsep}{3pt}
  \begin{minipage}{.5\textwidth}
  \centering
  
  \caption*{$\Omega_5$}%

  \begin{tabular}{l|rrc|rrc|rrc}
    \toprule
    \emph{l}& \multicolumn{3}{|c}{$p=1$} & \multicolumn{3}{|c}{$p=2$} & \multicolumn{3}{|c}{$p=3$}  \\
    \cmidrule(lr){2-4}\cmidrule(lr){5-7}\cmidrule(lr){8-10}
    & \emph{\#i} & \#DoF & CG & \emph{\#i}  & \#DoF & CG & \emph{\#i} & \#DoF & CG   \\
    \midrule
    2 & 4 & 64    & -  & 4 & 216    & -   & 4 & 512    & - \\
    3 & 5 & 512   & -  & 5 & 1728   & -   & 5 & 4096   & -  \\
    4 & 5 & 4096  & -  & 5 & 13824  & -   & 6 & 32768  & -   \\
    5 & 5 & 32768 & 88 & 6 & 110592 & 191 & 6 & 262144 & 259   \\
    \bottomrule
\end{tabular}
  \end{minipage}
  \begin{minipage}{.5\textwidth}
  
  \centering
  \caption*{$\Omega_6$}%

  \begin{tabular}{l|rrc|rrc|rrc}
    \toprule
    \emph{l}& \multicolumn{3}{|c}{$p=1$} & \multicolumn{3}{|c}{$p=2$} & \multicolumn{3}{|c}{$p=3$}  \\
    \cmidrule(lr){2-4}\cmidrule(lr){5-7}\cmidrule(lr){8-10}
    & \emph{\#i} & \#DoF & CG & \emph{\#i}  & \#DoF & CG & \emph{\#i} & \#DoF & CG   \\
    \midrule
    3 & 9 & 720    & -   & 11 & 2430    & -   & 12 & 5760    & -  \\
    4 & 9 & 5720   & -   & 12 & 19305   & -   & 12 & 45760   & -  \\
    5 & 9 & 45720  & -   & 12 & 154305  & -   & 13 & 365760  & -  \\
    6 & 9 & 365712 & 194 & 12 & 1234278 & 535 & 13 & 2925696 & 861\\
    \bottomrule
  \end{tabular}
  \end{minipage} 
  \caption{Number of required iterations, DoF per level, and iterations required by plain conjugate-gradient (CG) for grids $\Omega_5$, $\Omega_6$.}\label{tab:3D_examples_MG}
  \end{table}

